\documentclass[11pt]{amsart}
\pdfoutput=1
\usepackage{amscd,amssymb,amsmath,latexsym,enumerate}
\usepackage{mathrsfs}
\usepackage{epsfig}
\newtheorem{theo}{Theorem}
\newtheorem{defini}{Definition}
\newtheorem{proposi}{Proposition}
\newtheorem{lemma}{Lemma}
\newtheorem{coro}{Corollary}
\newtheorem{res}{Result}
\newtheorem{rem}{Remark}

\newtheorem{exam}{Example}

\newtheorem{prob}{Problem}

\newcommand{\Aa}{{\mathcal A}}
\newcommand{\Bb}{{\mathcal B}}
\newcommand{\Cc}{{\mathcal C}}

\newcommand{\Ee}{{\mathcal E}}
\newcommand{\Ff}{{\mathcal F}}

\newcommand{\Hh}{{\mathcal H}}

\newcommand{\Kk}{{\mathcal K}}

\newcommand{\Uu}{{\mathcal U}}

\newcommand{\Ww}{{\mathcal W}}

\newcommand{\id}{{\bf 1}}

\newcommand{\AG}{{\mathfrak A}}

\newcommand{\KG}{{\mathfrak K}}

\newcommand{\MG}{{\mathfrak M}}

\newcommand{\CM}{\mathbb C}
\newcommand{\EM}{\mathbb E}

\newcommand{\FM}{\mathbb F}
\newcommand{\KM}{\mathbb K}

\newcommand{\NM}{\mathbb N}
\newcommand{\RM}{\mathbb R}
\newcommand{\PM}{\mathbb P}
\newcommand{\QM}{\mathbb Q}

\newcommand{\TM}{\mathbb T}

\newcommand{\ZM}{\mathbb Z}

\newcommand{\hb}{{\hat{b}}}

\newcommand{\hn}{{\hat{n}}}

\newcommand{\hu}{{\hat{u}}}
\newcommand{\hy}{{\hat{y}}}

\newcommand{\hpi}{{\hat{\pi}}}
\newcommand{\hth}{{\hat{\theta}}}

\newcommand{\hD}{{\widehat{D}}}

\newcommand{\hK}{{\widehat{K}}}

\newcommand{\as}{{\mathscr A}}

\newcommand{\ms}{{\mathscr M}}

\newcommand{\ps}{{\mathscr P}}

\newcommand{\Cs}{$C^{\ast}$-algebra }         
\newcommand{\CS}{$C^{\ast}$-algebra}          
\newcommand{\Css}{$C^{\ast}$-algebras }       
\newcommand{\CsS}{$C^{\ast}$-algebras}        
\newcommand{\TR}{{\rm Tr\,}}                  
\newcommand{\Ker}{\mbox{\rm Ker}}             
\newcommand{\Aut}{\mbox{\rm Aut}}             
\newcommand{\Adj}{\mbox{\rm Ad}}              
\newcommand{\Orb}{{\mbox{\rm Orb}}}           
\newcommand{\Hull}{{\mbox{\rm Hull}}}         
\newcommand{\sign}{\mbox{\rm sign}}           
\newcommand{\Blip}{B_{\mbox{\tiny \it Lip}}}  
\newcommand{\Clip}{C_{\mbox{\tiny \it Lip}}}  
\newcommand{\Prob}{\mbox{\rm Prob}}           
\newcommand{\XP}{\Aa\rtimes_\alpha\ZM}        
\newcommand{\wu}{\mbox{\rm\sc wu}}            
\newcommand{\NMo}{{\mathbb N}_\ast}           

\def\SL{{\rm SL}}

\begin{document}

\title{Dynamical Systems on Spectral Metric Spaces}
\thanks{Work supported in part by NSF Grants No. 0901514, DMS-0651925 and DMS-1007207.}
\author{Jean V. Bellissard$^1$, Matilde Marcolli$^2$, Kamran Reihani$^3$ }

\address{$^1$Georgia Institute of Technology\\School of Mathematics\\
Atlanta GA 30332-0160}
\email{jeanbel@math.gatech.edu}

\address{$^2$California Institute of Technology\\Mathematics 253-37\\ Pasadena, CA 91125}
\email{matilde@caltech.edu}

\address{$^3$University of Kansas\\Department of Mathematics\\Lawrence, KS 66045-7594}
\email{reihani@math.ku.edu}

\begin{abstract}
Let $(\Aa,\Hh,D)$ be a spectral triple, namely: $\Aa$ is a \CS, $\Hh$ is a Hilbert space on which $\Aa$ acts and $D$ is a selfadjoint operator with compact resolvent such that the set of elements of $\Aa$ having a bounded commutator with $D$ is dense. A spectral metric space, the noncommutative analog of a complete metric space, is a spectral triple $(\Aa,\Hh,D)$ with additional properties which guaranty that the Connes metric induces the weak$^\ast$-topology on the state space of $\Aa$. A $\ast$-automorphism respecting the metric defined a dynamical system. This article gives various answers to the question: is there a canonical spectral triple based upon the crossed product algebra $\Aa\rtimes_{\alpha}\ZM$, characterizing the metric properties of the dynamical system ? If $\alpha$ is the noncommutative analog of an isometry the answer is yes. Otherwise, the metric bundle construction of Connes and Moscovici is used to replace $(\Aa,\alpha)$ by an equivalent dynamical system acting isometrically. The difficulties relating to the noncompactness of this new system are discussed. Applications, in number theory, in coding theory are given at the end.    
\end{abstract}

\maketitle


\pagestyle{myheadings}
\markboth{Dynamical Systems on Spectral Metric Spaces}{Dynamical Systems on Spectral Metric Spaces}


\tableofcontents


\vspace{.5cm}

\section{Introduction}
\label{met10.sect-intro}

\noindent Let $X=(\Aa, \Hh, D)$ be a spectral triple \cite{Co94} (also called a $K$-cycle), namely $\Aa$ is a separable \CS, $\Hh$ is a Hilbert space on which $\Aa$ is represented and $D$ is a selfadjoint operator with compact resolvent such that $\Cc^1(X)= \{a\in\Aa\,;\, \|[D,a]\| <\infty\}$ is a dense subset of $\Aa$. The condition $\|[D,a]\| <\infty$ means that $D$ admits a core invariant by multiplication by $a$ and that the operator $[D,a]=Da-aD$ extends as a bounded operator on $\Hh$. Unless otherwise stated, $\Aa$ will be identified with the image subalgebra of $\Bb(\Hh)$ by the representation. In some cases (see Section~\ref {met10.ssect-neqact} Example~\ref{met10.exam-shift}) it will be necessary to specify the representation, justifying the notation $(\Aa, \Hh,\pi, D)$, when $\pi$ denotes the representation.

\vspace{.1cm}

\noindent In view of the Gelfand Theorem, a \Cs can be seen as the noncommutative analog of the space of continuous functions on a locally compact space, vanishing at infinity. Spectral triples, under the name of {\em Fredholm modules} \cite{At67}, have been used initially to encode the cyclic cohomology of the noncommutative space defined by $\Aa$ \cite{Co94}. Eventually, Connes \cite{Co94} realized that they can also encode metric structures on this space. Such a structure extends as a metric on the state space of $\Aa$. It is remarkable that the extension of a metric from a space to its set of probability measures had been defined  and studied for a long time before, starting with the work of Kantorovich and Rubinstein \cite{KR57} for compact metric spaces, in the context of the mass transportation problem. It was later extended by Wasserstein \cite{Wa69} and Dobrushin \cite{Do70} to non compact complete metric spaces. Dobrushin proposed the name {\em Wasserstein distance} to such metrics, and the name has been used since then. In the present paper, though, we will refer to it as the {\em Connes distance} in the context of noncommutative metric spaces. It will be seen in Section~\ref{met10.ssect-defres} what conditions on a spectral triple are needed to encode a structure of locally compact complete metric space. The main problem investigated in this paper is the following

\begin{prob}
\label{met10.prob-main}
Let $\Aa$ be unital and $\alpha$ be a $\ast$-automorphism of $\Aa$. Given a spectral triple $X=(\Aa,\Hh,D)$, is there a canonical way of defining a spectral triple $Y=(\Aa\rtimes_{\alpha} \ZM, \Kk, \tilde{D})$ induced by $(\Aa, \Hh, D)$, for the crossed product of the original algebra $\Aa$ by the action of $\ZM$ defined by $\alpha$~? If so, what are the metric properties encoded by this triple~?
\hfill $\Box$
\end{prob}

\noindent It will be shown in this work, that the answer is not always positive: actually, it is positive if and only if $(\Aa, \Hh, D)$ represents the noncommutative analog of a metric space for which the group generated by $\alpha$ is an {\em equicontinuous} family of {\em quasi-isometries}. In particular, for uniformly hyperbolic actions on a compact metric space, such a construction is not possible. However, the present work offers a way to deal with this problem: following the philosophy introduced by Connes and Moscovici in the nineties \cite{CM95}, it is proposed to extend the algebra $\Aa$ to an algebra $\Bb$ representing the noncommutative analog of the set of continuous functions on the {\em metric bundle}, on which such action becomes isometric, so equicontinuous. This will permit to construct a spectral triple for the crossed product $\Bb\rtimes_{\alpha}\ZM$, in which $\Aa\rtimes_{\alpha}\ZM$ acts as a multiplier algebra. While it is possible to extend the construction to more general group actions, this work will deal only with the case of a single automorphism.

\vspace{.1cm}

\noindent The metric bundle construction comes with a price: it addresses the problem of spectral triples on {\em non unital} \CsS. This question was studied in particular by Latr\'emoli\`ere \cite{La07}. The main difficulty in this case is that the state space is not closed for the weak$^\ast$ topology. Thanks to the notion of {\em weak-uniform topology} on $\Aa$ it became possible to characterizes those Connes metrics on the state space giving rise to the weak$^\ast$ topology. Because of the importance of this result for the purpose of the present work, it will be stated explicitly (see Result~\ref{met10.res-latr}, in Section~\ref{met10.ssect-defres}). 

\vspace{.1cm}

\noindent It will be seen, however, that the work by Latr\'emoli\`ere does not solve all problems encountered with the present construction. For indeed, while it is sufficient whenever the metric is bounded, namely when the state space has a finite diameter in the Connes metric, it turns out to be insufficient for the examples that will be presented here at the end of this paper. Something analogous occurs with the Wasserstein distance when the metric space is complete but not compact and the metric is unbounded. In the latter case this distance encodes the $w^\ast$-topology on the set of probabilities only on compact subsets of finite diameters. Such subsets generate a dense subset of the probability space. For a spectral triple, this situation occurs whenever the {\em Lipschitz ball}, namely the set of $a\in\Aa$ for which $\|[D,a]\|\leq 1$, is unbounded. If so, it is expected that indeed there are, in the state space, enough compact subsets with finite diameter in the Connes metric on which this metric generates the $w^\ast$-topology. This problem will not be addressed in the present work though. 

\vspace{.3cm}

 \subsection{Definitions and Previous Results}
 \label{met10.ssect-defres}

\noindent Given a spectral triple $(\Aa,\Hh,D)$, let $d_C$ be the {\em Connes} pseudo-distance defined on the state space of $\Aa$ as follows \cite{Co94}:

$$d_C(\rho, \omega) =
   \sup \{|\rho(a)-\omega(a)|\,;\, \|[D,a]\|\leq 1\}\,
\hspace{2cm}
    \rho,\omega \in\MG_1(\Aa)\,.
$$

\noindent Here, $\MG_1(\Aa)$ will denote the state space of $\Aa$, using a notation coming from probability theory. The first problem is to decide whether this is a metric on $\MG_1(\Aa)$ and whether this metric defines the $w^\ast$-topology. The present answers to these questions are described here.

\begin{res}[\cite{RV06}]
\label{met10.res-rv}
Let $(\Aa,\Hh,D)$ be a spectral triple with $\Aa$ unital. The Connes pseudo-metric is a metric on the state space if the following two conditions hold

\noindent (i) the representation of $\Aa$ in $\Hh$ is non-degenerate, namely that $\Aa\Hh=\Hh$,

\noindent (ii) the metric commutant $\Aa_D'=\{a\in\Aa\,;\, [D,a]=0\}$ is trivial, namely it is reduced to multiples of the unit.
\end{res}

\noindent In particular this condition requires that the representation of $\Aa$ on the Hilbert space $\Hh$ be {\em faithful}.

\begin{res}[\cite{Pa98,Ri98,Ri99,Ri04,Oz-Ri}]
\label{met10.res-PavRie}
Let $X=(\Aa,\Hh,D)$ be a spectral triple with $\Aa$ unital. In addition, let $X$ be such that the representation of $\Aa$ in $\Hh$ is non-degenerate and the metric commutant is trivial. Then the Connes metric defines the $w^\ast$-topology if and only if the Lipschitz ball, defined by

$$\Blip  = \Blip(X)=
   \{a\in \Aa\,;\, \|[D,a]\|\leq 1\}\,,
$$

\noindent is pre-compact in the quotient space $\Aa/\Aa_D'$. Equivalently, $\Ker(\phi)\cap\Blip$ is pre-compact in $\Aa$, for some (hence for all) $\phi\in\MG_1(\Aa)$.
\end{res}

\noindent As proved in Lemma~\ref{met10.lem-blipclosed} below, $\Blip$ is closed in norm. If $\Aa$ is not unital, but separable, and if the Lipschitz ball is norm bounded, the main result of Latr\'emoli\`ere is the following (see \cite{La07} Th. 2.6 \&~4.1)

\begin{res}[\cite{La07}]
\label{met10.res-latr}
Let $X=(\Aa, \Hh,D)$ be a spectral triple with $\Aa$ a nonunital, separable \CS. If the $D$-commutant is trivial and the representation of $\Aa$ is nondegenerate, the Connes distance defines the $w^\ast$-topology on the state space if and only if there is a strictly positive element $h\in\Aa$ such that $h\Blip h$ is compact. In such a case the Connes metric is a path metric and the state space is complete for this metric. 
\end{res}

\noindent It is worth reminding the reader that a {\em strictly positive element} $h$ is a positive element such that $h\Aa h$ is norm dense in $\Aa$. This leads to the following definition, used in the present paper

\begin{defini}
\label{met10.def-qms}
A bounded spectral metric space is a spectral triple $X=(\Aa,\Hh,D)$ such that

(i) the representation of $\Aa$ on $\Hh$ is non-degenerate,

(ii) the metric commutant $\Aa_D'=\{a\in\Aa\,;\, [D,a]=0\}$ is trivial,

(iii) there is a strictly positive element $h\in\Aa$ such that the $h$-compressed Lipschitz Ball $h\Blip h$ has precompact image in the normed space $\Aa/\Aa_D'$.

\vspace{.1cm}

\noindent A spectral metric space $X=(\Aa,\Hh,D)$ will be called compact if $\Aa$ is unital.
\end{defini}

\noindent The term {\em spectral metric space} is appropriate since it is based on a spectral triple\footnote{The authors thank M.Rieffel for suggesting this name.}. 
The previous definition is a bit more restrictive than the term {\em quantum metric space} defined by Rieffel \cite{Ri04}. For indeed, here $\Aa$ is restricted to be a \Cs instead of an order-unit space. 

\vspace{.1cm}

\noindent It will be seen in Sections~\ref{met10.sect-NCtori} \& \ref{met10.sect-codes}, that there are important examples of spectral triple that do not satisfy the criterion of Result ~\ref{met10.res-latr}. This is because the corresponding Lipschitz balls are unbounded in norm. The problem of extending Result ~\ref{met10.res-latr} to the unbounded case will not be addressed in the present work. But partial results can be found in Section~\ref{met10.ssect-zm} in which the following important elementary example is considered: let $\ZM$ be endowed with a metric $d_\ZM$. Then let $c_0(\ZM)$ be the algebra of complex valued sequences vanishing at infinity. It is not unital. It acts on the Hilbert space $\Hh=\ell^2(\ZM)$ by pointwise multiplication. In Section~\ref{met10.ssect-zm}, several Dirac operators will be proposed to describe this metric. Each construction requires to enlarge $\Hh$ by tensoring it with some finite dimensional Hilbert space $\Ee$ on which some Clifford algebra acts. It will be shown that the corresponding spectral triple is a compact metric spectral space if and only if the distance $d_\ZM$ is summable, namely $\sum_n d_\ZM(n,n+1) <\infty$. A second construction allows to include translation invariant metrics (thus nonsummable): then it satisfies the criterion of Result~\ref{met10.res-latr} if and only if the metric $d_\ZM$ is bounded. Both results exclude the usual metric on $\ZM$. This is because the usual metric gives $\ZM$ an infinite diameter. Therefore the Connes metric associated with it is unbounded on the state space. However, the following result, valid in the commutative case, applies

\begin{res}[\cite{Do70}]
\label{met10.res-dob}
Let $(X,d)$ be a complete separable metric space and let $\MG_1(X)$ be the set of probability measures on $X$. A subset $F\subset\MG_1(X)$ is called $d$-tight if there is $x_0\in X$ such that for all $\epsilon >0$ there is $r>0$ for which $\int_{d(x_0,x)\geq r} \mu(dx) \;d(x_0,x) \leq \epsilon$ for all $\mu\in F$. Then the Wasserstein distance generates the $w^\ast$-topology on each weak$^\ast$-closed $d$-tight set.
\end{res}

\noindent The extension of such a result to noncommutative spectral triple is an interesting but difficult open problem. It will not be considered in the present work.

\vspace{.3cm}

 \subsection{Main Results}
 \label{met10.ssect-main}

\noindent The set $\Aut (\Aa)$ of $\ast$-automorphisms of $\Aa$ will be endowed with the norm pointwise topology, namely a basis of neighborhoods of an element $\alpha\in\Aut(\Aa)$ is the family of sets of the form $\Uu(\alpha;a,\epsilon) = \{\beta\in\Aut(\Aa)\,;\, \|\beta(a)-\alpha(a)\| <\epsilon\}$, where $a\in\Aa$ and $\epsilon >0$. 

\begin{defini}
\label{met10.def-qiso}
For $X=(\Aa,\Hh,D)$ a spectral metric space, a $\ast$-automorphism $\alpha\in\Aut(\Aa)$ will be called a quasi-isometry if $a\in \Cc^1(X)\,\Leftrightarrow \alpha(a)\in \Cc^1(X)$. 
\end{defini}

\noindent For $G$ a subgroup of $\Aut(\Aa)$, $\Cc^1(G,X)$ will denote the set of elements $a\in\Aa$ such that (i) $\|[D,g(a)]\|<\infty$ for all $g\in G$ and (ii) the map $g\in G \mapsto [D, g(a)]\in B(\Hh)$ is continuous in norm. In particular $\Cc^1(G,X)\subset \Cc^1(X)$. Similarly, $\Cc_b^1(G,X)\subset \Cc^1(G,X)$ will denote the set of elements $a\in\Aa$ such that $\sup_{g\in G}\|[D, g(a)]\|<\infty$. It is easy to check that $\Cc^1(X), \Cc^1(G,X),\Cc_b^1(G,X)$ are $\ast$-subalgebras of $\Aa$ invariant by holomorphic functional calculus (see Lemma~\ref{met10.lem-lipnorm} for instance). Transposed in the classical world,  the algebra $\Cc^1(G,X)$ contains $\Cc^1(X)$ if and only if $G$ acts by bi-Lipschitz transformations, namely by quasi-isometries. Similarly, if $\Cc_b^1(G,X)$ contains $\Cc^1(X)$ means that $G$ is an equicontinuous family of quasi-isometries. This justifies the following definitions used in the present work

\begin{defini}
\label{met10.def-smoothact}
Let $X=(\Aa,\Hh,D)$ be a spectral metric space. Then a subgroup $G\subset \Aut(\Aa)$ will be called quasi-isometric, if $\Cc^1(X)= \Cc^1(G,X)$. It will be called equicontinuous if $\Cc^1(X)= \Cc_b^1(G,X)$. It will be called isometric if it is quasi-isometric and $\|[D,g(a)]\|=\|[D,a]\|$ for all $g \in G$ and $a \in \Cc^1(X)=\Cc^1(G,X)$. The same definition applies to a single automorphism provided $G$ is the group it generates. The pair $(X,G)$ will be called a spectral dynamical system.
\end{defini}

\begin{rem}
\label{met10.rem-qiso}
{\em  A quasi-isometric automorphism of a {\em compact} spectral metric space $X=(\Aa,\Hh,D)$ is an isometry if and only if it preserves the Connes distance on the state space $\MG_1(\Aa)$ of $\Aa$. For indeed the ``if'' part follows from the fact that $a\mapsto\|[D,a]\|$ defines a lower semicontinuous Lipschitz seminorm on $\Cc^1(X)$ (\cite{Ri99}, Proposition 3.7), therefore, following \cite{Ri99}, Theorem 4.1, it can be recovered via Connes distance by 

$$\|[D,a]\|=
   \sup\left\{
     \frac{|\rho(a)-\omega(a)|}{d_C(\rho,\omega)}\,;\, 
      \rho\neq\omega\in \MG_1(\Aa)
       \right\}\,.
$$

\noindent The ''only if'' part is obvious from the definition of an isometry.
}
\hfill $\Box$
\end{rem}

\noindent In the classical world, equicontinuity implies compactness through the Arzel\`a-Ascoli theorem. An extension was considered in \cite{ACh04} in the noncommutative world, for discrete groups of rapid decay. The result in \cite{ACh04} leads to

\begin{theo}[Arzel\`a-Ascoli theorem]
\label{met10.th-qi}
Let $X=(\Aa,\Hh,D)$ be a compact spectral metric space. Let $G\subset \Aut(\Aa)$ be a quasi-isometric subgroup. Then $G$ is equicontinuous if and only if it has a compact closure.
\end{theo}

\noindent In a compact metric space, given a compact group of quasi-isometries, it is possible to find an equivalent metric that is invariant by the group. Such is the case also for spectral metric spaces, as shown by the next result. To express it, the simplest case is considered, namely the group generated by a single automorphism $\alpha\in\Aut(\Aa)$. This group defines a $\ZM$-action on the spectral metric space. Let $\Aa\rtimes_\alpha\ZM$ denote the crossed-product algebra, namely the universal \Cs generated by the elements of $\Aa$ and by a unitary $u$ such that $uau^{-1} = \alpha(a)$. The crossed product is an algebraic version of the existence of a unitary operator {\em implementing} the automorphism $\alpha$ in a representation. In addition, the dual group $\TM$ of $\ZM$ acts on the crossed product \cite{Co73,TaM73,TaH75} in the following way: $\eta_k(a)=a$ if $a\in\Aa$ and $\eta_k(u)=e^{\imath k}u$ for $k\in\TM$. In such a case, the following holds

\begin{theo}
\label{met10.th-qiCross}
Let $\Aa$ be a unital, separable \Cs and $\alpha\in\Aut(\Aa)$. Let $u$ denote the unitary implementing $\alpha$ in the crossed product algebra $\XP$. Then there is a compact spectral metric space $X=(\Aa,\Hh,D)$ based on $\Aa$ with $\alpha$ implementing an equicontinuous $\ZM$-action if and only if there is a compact spectral metric space $Y=(\XP,\Kk,\hD)$, based on the crossed product algebra, such that 

(i) the dual action on $\Aa\rtimes_\alpha\ZM$ is quasi-isometric, 

(ii) $u^{-1}[\hD,u]$ is bounded and commutes with the elements of $\Aa$, 

(iii) the Connes metrics induced on the state space of $\Aa$ associated with the two spectral metric spaces $X$ and $Y$ are equivalent.
\end{theo}

\noindent The dual action is required to identify $u$ within the algebra $\XP$. It will be proved that it acts isometrically on $Y$. If $u$ satisfies the condition (ii), then $[D,\alpha^n(a)]=u^n[D,a]u^{-n}$, meaning that the Connes metric defined by the spectral triple $X=(\Aa,\Kk,\hD)$, obtained by restricting $Y$ to $\Aa$, is $\alpha$-invariant. Conversely the construction of the spectral metric space $Y$ involves the use of the socalled {\em left regular representation} of $(\Aa,\alpha)$. 

\vspace{.2cm} 

\noindent Even in the commutative case there are examples of spectral metric spaces for which certain actions are quasi-isometric in the previous sense without being equicontinuous. A simple example is provided by the Arnold cat map (see Section~\ref{met10.ssect-neqact}). More generally, this is the case for uniformly hyperbolic actions. In such a case, building a spectral triple on the crossed product algebra becomes impossible using the construction leading to the proof of the Theorem~\ref{met10.th-qiCross}. Thanks to the seminal work of Connes and Moscovici \cite{CM95}, it is possible to overcome this difficulty. The main idea is the noncommutative analog of replacing the manifold by its {\em metric bundle}. As it turns out, in the classical case, this metric bundle is a quotient of the frame bundle that allows to build a diffeomorphism invariant measure on any manifold. The extension of such a construction to the noncommutative case is not yet available for spectral metric spaces in general. The automorphism groups might be too big for such a construction to be available in general. However, if $G$ is a subgroup of $\Aut(\Aa)$, as the image of a locally compact group by a continuous group homomorphism, it is sufficient to restrict this bundle to an orbit of the group, so that the fiber can be replaced by the group itself. In the case of a $\ZM$-action, this gives the following construction

\begin{theo}
\label{met10.th-bundle}
Let $X=(\Aa,\Hh,D)$ be a spectral metric space and let $\alpha\in\Aut(\Aa)$ be quasi-isometric. If $\alpha$ is not equicontinuous, there is a spectral metric space $Y$ based on $\Aa\otimes c_0(\ZM)$, induced by $X$, such that 

(i) $\Aa$ acts as a multiplier algebra on $Y$

(ii) $\alpha$ implements an equicontinuous action $\alpha_\ast$ on $Y$

(iii) there is a spectral metric space based on $\Aa\otimes c_0(\ZM)\rtimes_{\alpha_\ast}\ZM$, on which $\Aa\rtimes_{\alpha}\ZM$ acts as a multiplier algebra.
\end{theo}

\noindent It is important to note that, even if $X$ is compact, namely $\Aa$ is unital, the space $Y$ is non-compact in general. In the present paper, it will not be investigated whether the construction leads to a spectral triple that generates the $w^\ast$-topology on the noncommutative analog of $d$-tight subsets of the state space. 

\vspace{.1cm}

\noindent To illustrate the potentiality of this construction, three examples have been added at the end of the present paper. The metric bundle construction permits to complete prior results. The first example concerns an interpretation of some number theoretic quantities, the Shimizu $L$-functions, occurring as topological invariant through the Index Theorem in computing the signature of Hilbert's modular surfaces associated with totally real number fields \cite{ADS}. It was shown that the only contribution was coming from cusps that are seen as the boundary of a manifold, which turns out to be a solvmanifold, namely it comes from a quotient of a Lie algebra, all element of which are strictly hyperbolic. In a previous work \cite{Mar-solv}, these $L$-functions were shown to comes from a spectral triple over a non-commutative torus, in the case of a real quadratic number field. The corresponding solvmanifold, can be interpreted, modulo homotopy, as the crossed product of the noncommutative torus by the group of totally positive units of the number field. Due to the hyperbolicity of the action, building a spectral triple representing the metric on the solvmanifold was a problem that the present paper helps to solve. The other example, the algebro-geometric codes \cite{ManMar}, is coming form coding theory in Computer Science. This example is actually very close to the problem of finding a spectral triple on a tiling space and is a sophisticated extension of the example of a two-sided shift acting on the Cantor set and treated in Section~\ref{met10.ssect-neqact}, Example~\ref{met10.exam-shift}. In both cases, however, the natural metric that is behind the construction corresponds to a noncompact unbounded metric space. This is an incentive to extend the theory of unbounded complete metric spaces to the noncommutative world.

\vspace{.3cm}

 \subsection{A Bit of History}
 \label{met10.ssect-rev}

\noindent The notion of spectral triple, in view of representing the noncommutative analog of Riemannian manifolds, was proposed by Connes in the eighties \cite{Co94} (see Chapter 6). It came out of the construction of Fredholm modules used in the construction of cyclic cohomology. Immediately after understanding the conceptual foundation in this part of the theory \cite{Co88} Connes gave few important examples: the case of Riemannian manifolds, the Cantor set and the Julia sets \cite{Co94}. He also addressed the inverse problem, namely, when is a spectral triple coming from a smooth Riemannian manifold \cite{Co95}. However, the subject did not develop until the late nineties when several authors proposed using them for different purposes. The earliest work recorded in the literature came in 1994, when Lapidus proposed to extend the Connes example of the Cantor set to fractals \cite{La94,La97}. Then several works appeared, like the important contributions of Buyalo \cite{Bu99,Bu00} or the ones by Guido \& Isola \cite{GI01,GI03,GI05} in the analysis of metric spaces and fractals. On the more conceptual side were the works by Pavlovi\'c \cite{Pa98} and by Rieffel \cite{Ri98,Ri99,Ri04} and the extension to noncompact metric spaces by Latr\'emoli\`ere \cite{La07} that are at the origin of the Definition~\ref{met10.def-qms} in the present paper. The reconstruction problem for Riemannian manifold has been considered by Rennie \cite{Re01}, Rennie and Varilly \cite{RV06} and Connes \cite{Co08}. Several important classes of noncommutative spectral triples have been constructed by Rieffel \cite{Ri04} and by Christensen, Antonescu-Ivan and various collaborators \cite{ChI06,ChI07,CIL,CIS}. Recently, the Riemannian geometry of ultrametric Cantor sets was described in \cite{PB09} leading to the definition of an analog to the Laplace-Beltrami operator and some developments in the theory of tilings \cite{JS09}. Very recently the extension to compact metric space has been considered. In \cite{FS10} the author use a Christensen-Ivan spectral triple to encode various metric invariant on a multifractal space. In \cite{Pa10} an extension of \cite{PB09} to all compact spaces is treated. The problem of group actions on noncommutative metric spaces was considered by Connes in the context of gauge theory and the notion of Yang-Mills action principle in \cite{Co94}. However, the notion developed in the present paper have not been considered so far, in a systematic way. The present paper is therefore an introduction to this topic.

\vspace{.3cm}

 \subsection{Organization of the Paper}
 \label{met10.ssect-org}

\noindent The rest of the paper is organized as follows. In Section~\ref{met10.sect-compap}, the case of equicontinuous action is studied and is shown to correspond to almost periodicity in the sense of H.~Bohr. Then Section~\ref{met10.sect-stEqui} is dedicated to the construction of a canonical spectral metric space over the crossed product algebra whenever the action is equicontinuous. As it turns out there are some non trivial technicalities, in particular, the proof of the Proposition~\ref{met10.prop-lipcomp}, concerning the compactness of the Lipschitz ball, which is postponed to the Appendix~\ref{met10.sect-Xprod}. Then the Section~\ref{met10.sect-smb} is dedicated to the notion of spectral metric bundle. Two examples are given of non equicontinuous actions, the Arnold cat map acting on the $2$-torus and the shift acting on $\{0,1\}^\ZM$, the most elementary model of a Smale space. The study of the space of Euclidean metrics on a linear space gives a clue to build the spectral metric bundle. Then a long subsection is dedicated to the study of various metric structures on $\ZM$, the most elementary model of noncompact metric space, leading to introduce various technicalities like the $\wu$-topology \cite{La07}. The last three Sections, \ref{met10.sect-NCtori}, \ref{met10.sect-CKAF} and \ref{met10.sect-codes}, treat examples already developed previously by one of the present authors, but which were not completely treated before by lack of the notion of spectral metric bundle. 

\vspace{.3cm}

\noindent {\bf Acknowledgments: } This work was supported in part by NSF Grant No. 0901514 (J.B.), DMS-0651925 and DMS-1007207 (M.M). J.B. thanks the SFB 701 (Universit\"at Bielefeld, Germany) and the MAPMO and the F\'ed\'eration Denis Poisson (Orl\'eans University, France), for providing office space during the writing of this paper. He also wants to thank Jean Savinien and Ian Palmer for discussions and for providing references used in this work. K.R. thanks Marc Rieffel for various comments and for pointing the reference \cite{La07}.

\vspace{.5cm}

\section{Compactness and Almost Periodicity}
\label{met10.sect-compap}

\noindent This section is devoted to the notion of almost periodicity and to proving Theorem~\ref{met10.th-qi}. The standard references for the theory of almost periodic functions are the books by H. Bohr \cite{Bo47}, Besicovich \cite{Be32}, and the extension to all groups in the articles by von Neumann \cite{vN34,BvN35}.

\vspace{.3cm}

 \subsection{The Lipschitz Algebra}
 \label{met10.ssect-lipalg}

\noindent Some technical tools, presented here will be useful. In this Section $\Aa$ is a separable \Cs but it is not required to be unital.

\begin{lemma}
\label{met10.lem-lipnorm}
Let $X=(\Aa,\Hh,D)$ be a spectral metric space. The space $\Cc^1(X)$, endowed with the norm

$$\|a\|_1 = \|a\| + \|[D,a]\|\,,
$$

\noindent is a Banach $\ast$-algebra invariant by holomorphic functional calculus. The injection map $i:\Cc^1(X)\to \Aa$ is a compact homomorphism of Banach $\ast$-algebras. 
\end{lemma}

\noindent  {\bf Proof: } (i) That $\|\cdot \|_1$ defines a seminorm is obvious. Since $\|a\|_1=0\,\Rightarrow \,\|a\|=0\,\Rightarrow a=0$, it is a norm. That this norm is $\ast$-invariant is also clear, since $D$ is selfadjoint. Thanks to the Leibniz formula, $[D,ab]= [D,a]b+ a[D,b]$ it follows that

$$\|ab\|_1 \leq 
    \|a\|\|b\| + \|[D,a]\|\|b\| + \|a\|\|[D,b]\|\leq
     \|a\|_1\, \|b\|_1\,,
$$

\noindent showing that this is an algebraic norm.

\vspace{.1cm}

\noindent (ii) Let now $(a_n)_{n\in\NM}$ be a Cauchy sequence in $\Cc^1(X)$. Then it is Cauchy for the norm in $\Aa$ and therefore there is $a\in\Aa$ such that $\lim_{n\rightarrow \infty} \|a-a_n\| =0$. Moreover it follows that, for all $n\in\NM$, the sequence  $(\|a_n-a_m\|)_{m\in\NM}$ is Cauchy and thus converges in $\CM$ to a limit $c_n$ that converges to zero as $n\rightarrow \infty$. Let now $f,g\in\Hh$ belonging to the domain of $D$. It follows that $\langle Df|(a-a_n) g\rangle - \langle f|(a-a_n) Dg\rangle$ converges to zero as $n\rightarrow \infty$. Moreover

\begin{eqnarray*}
|\langle Df|(a-a_n) g\rangle -
    \langle f|(a-a_n) Dg\rangle| &=& 
\lim_{m\rightarrow\infty}
      |\langle Df|(a_m-a_n) g\rangle -
       \langle f|(a_m-a_n) Dg\rangle| \\
&\leq &\|f\|\|g\| \limsup_{m\rightarrow\infty}\|a_m-a_n\|_1
\end{eqnarray*}

\noindent Hence $\|a-a_n\|_1$ is finite and $\lim_{n\rightarrow \infty }\|a-a_n\|_1=0$, showing that $a\in\Cc^1(X)$ and that $a_n$ converges to $a$ in $\Cc^1(X)$. Therefore $\Cc^1(X)$ is a Banach $\ast$-algebra.

\vspace{.1cm}

\noindent (iii) Let $a\in\Cc^1(X)$. Then, its resolvent $(z\id-a)^{-1}$ either belongs to $\Aa$, if $\Aa$ is unital, or $(z\id-a)^{-1}-\id/z\in\Aa$ if not. Moreover, thanks to the Leibniz rule,

$$[D,(z\id-a)^{-1}-\id/z] = (z\id-a)^{-1}[D,a](z\id-a)^{-1}\,,
$$

\noindent so that $(z\id-a)^{-1}-\id/z \in \Cc^1(X)$. In addition

$$\|[D,(z\id-a)^{-1}-\id/z]\|_1\leq 
    \left\{\|(z\id-a)^{-1}-\id/z\|+1/|z|\right\}^2\, \|a\|_1
$$

\noindent This implies that for any function $F(z)$ holomorphic in a neighborhood $U$ of the spectrum of $a$ (and vanishing at $z=0$ whenever $\Aa$ is not unital), and any Jordan path $\gamma$ contained in $U$, surrounding the spectrum of $a$, the Cauchy integral $\oint_\gamma F(z) \{(z\id-a)^{-1} -1/z\}dz/2\imath \pi$ converges in $\Cc^1(X)$ showing that $F(a)\in\Cc^1(X)$. Hence, $\Cc^1(X)$ is invariant by holomorphic functional calculus.

\vspace{.1cm}

\noindent (iv) Since the unit ball of $\Cc^1(X)$ maps to $\Blip$ in the quotient space $\Aa/\CM\id$, it follows that it is compact in the norm topology of $\Aa$ (see Section~\ref{met10.ssect-quot}). Hence the canonical injection maps bounded sets into compact ones, meaning it is compact.
\hfill $\Box$

\vspace{.2cm}

\noindent A very similar proof leads to

\begin{coro}
\label{met10.cor-cc1b}
Let $X=(\Aa,\Hh,D)$ be a spectral metric space and let $G\subset \Aut(\Aa)$ be an equicontinuous group of automorphisms. Then the space $\Cc_b^1(G,X)$ endowed with the norm

$$\|a\|_{1,G} = \|a\|+
   \sup_{g\in G} \|[D,ga]\|
$$

\noindent is a Banach $\ast$-algebra, invariant by holomorphic functional calculus and the canonical injection $i:\Cc_b^1(G,X)\to \Aa$ is compact.
\end{coro}

\noindent The only additional property needed in the proof is the continuity of the maps $g\in G\mapsto [D,ga]\in \Bb(\Hh)$ for each $a\in \Cc_b^1(G,X)$. But the uniformity of the norm over $G$ implies that the limit of a Cauchy sequence is necessarily continuous. 

\begin{lemma}
\label{met10.lem-qiso}
Let $X=(\Aa,\Hh,D)$ be a spectral metric space. A $\ast$-automorphism $\alpha\in\Aut(\Aa)$ is quasi-isometric if and only if one of the two following conditions hold

(i) it defines a bounded $\ast$-automorphism of $\Cc^1(X)$,

(ii) the action of $\alpha$ on the state space of $\Aa$ is bi-Lipschitz for the Connes metric.
\end{lemma}

\noindent  {\bf Proof: } (i) If $\alpha :\Cc^1(X)\to \Cc^1(X)$ is bounded, it follows, from the Definition~\ref{met10.def-qiso}, that $\alpha$ is quasi-isometric. Conversely if $\alpha$ is quasi-isometric, it is a linear map $\alpha :\Cc^1(X)\to \Cc^1(X)$ which is defined everywhere. By the closed graph theorem it is bounded. That it is a $\ast$-homomorphism is obvious from the definition. The same can be said for $\alpha^{-1}$ since $\alpha$ is one-to-one and since, thanks to Definition~\ref{met10.def-qiso}, it is onto as well. 

\vspace{.1cm}

\noindent (ii) Let $\alpha\in\Aut(\Aa)$ be quasi-isometric. Then by the previous argument, $\alpha: \Cc^1(X)\to \Cc^1(X)$ is a bounded $\ast$-isomorphism. Then, if $\rho, \omega$ are two states on $\Aa$

$$d_C(\rho\circ\alpha,\omega\circ\alpha )=
  \sup\{|\rho(a)-\omega(a)|\,;\, \|[D,\alpha^{-1}(a)]\|\leq 1\}
$$

\noindent Since $\alpha$ is bounded on $\Cc^1(X)$ it follows that there is $K>0$ such that $\|[D,\alpha^{-1}(a)]\|\leq 1\,\Rightarrow \|[D,a]\|\leq K$. Therefore $d_C(\rho\circ\alpha,\omega\circ\alpha )\leq Kd_C(\rho,\omega)$. Replacing $\alpha$ by its inverse leads to

\begin{equation}
\label{met10.eq-biLip}
\frac{1}{K} d_C(\rho,\omega) \leq
   d_C(\rho\circ\alpha,\omega\circ\alpha )\leq 
    Kd_C(\rho,\omega)\,.
\end{equation}

\noindent Conversely, if eq.~(\ref{met10.eq-biLip}) holds, it follows that both $\alpha, \alpha^{-1}$ leaves $\Cc^1(X)$ invariant, showing that $\alpha$ is quasi-isometric.
\hfill $\Box$

\vspace{.3cm}

 \subsection{Almost Periodicity: Definitions}
 \label{met10.ssect-ap}

\begin{defini}
\label{met10.def-ap}
If $\alpha\in\Aut(\Aa)$ an element $a\in\Aa$ is called $\alpha$-almost periodic if the norm closure $\Hull(a)$ of its orbit $\Orb(a)=\{\alpha^n(a)\,;\, n\in\ZM\}$ is compact. A subset $F\subset \Aa$ is called $\alpha$-almost periodic if any of its elements are $\alpha$-almost periodic. A $\ast$-automorphism $\alpha\in\Aut(\Aa)$ is called {\em almost periodic} whenever any $a\in\Aa$ is $\alpha$-almost periodic. 
\end{defini}

\noindent An equivalent way of defining almost periodicity is the following

\begin{coro}
\label{met10.cor-apcomp1}
Let $\Aa$ be a \CS. A $\ast$-automorphism $\alpha\in\Aut(\Aa)$ is almost periodic if and only if the norm pointwise closure of the group it generates is compact. 
\end{coro}

\noindent The proof will be left to the reader since it consists only of using the definition of the norm-pointwise convergence. In much the same way, if $F\subset \Aa$ is $\alpha$-almost periodic, then let $\Aa_F$ be the \Cs generated by the elements of the form $\alpha^n(a)$ for $n\in\ZM$ and $a\in F$. Then $\alpha$ is an almost periodic $\ast$-automorphism of $\Aa_F$. In particular it generates a compact group of automorphisms. However, this group may not be contained in $\Aut(\Aa)$.

\vspace{.3cm}

 \subsection{Almost Periodicity: Bohr's Approach}
 \label{met10.ssect-bohr}

\noindent In this section, a more classical approach to almost periodicity, close to Bohr's original one, will be presented for the sake of the reader. Let $a\in\Aa$ be $\alpha$-almost periodic. Then, its hull is $\alpha$-invariant and $(\Hull(a),\alpha)$ becomes a topological dynamical system. An $\alpha$-periodic element is an $a\in\Aa$ such that there is $p\in\NM$ ($p\neq 0$), called a {\em period}, such that $\alpha^p(a)=a$. The orbit of a periodic element is finite thus compact. In particular if $a$ is periodic then it is almost periodic.

\begin{proposi}
\label{met10.prop-hull}
If $a\in\Aa$ is $\alpha$-almost periodic, $\Hull(a)$ inherits a canonical structure of compact abelian group in which $\Orb(a)$ is a dense subgroup isomorphic to $\ZM$.
\end{proposi}

\noindent {\bf Proof: } If $a$ is periodic, the result is elementary. Let $a$ be non periodic then. Let $S(a)$ be the set of infinite sequences $\xi=(x_k)_{k\in\NM}$ such that (i) $x_k\in\ZM$ for all $k$'s and (ii) $\lim_{k\rightarrow \infty} \alpha^{x_k}(a)$ exists in $\Hull(a)$. By definition $S(a)$ is non-empty: if $n\in\ZM$, then the constant sequence $\hn=(n_k)_{k\in\ZM}$ is defined by $n_k=n$ for all $k$ and is obviously in $S(a)$. Then given two such sequences $\xi=(x_k)_{k\in\NM},\eta=(y_k)_{k\in\NM}$, their sum is defined by $\xi+\eta=(x_k+y_k)_{k\in\NM}$. It follows that $\xi+\eta$ belongs to $S(a)$. For indeed, if $k,l$ are two integers, $\|\alpha^{x_k+y_k}(a)-\alpha^{x_l+y_l}(a)\|\leq 
\|\alpha^{x_k}\{\alpha^{y_k}(a)-\alpha^{y_l}(a)\}\|+ 
\|\alpha^{x_k}(\alpha^{y_l}(a))-\alpha^{x_l}(\alpha^{y_l}(a))\|= \|\alpha^{y_k}(a)-\alpha^{y_l}(a)\|+ \| \alpha^{x_k}(a)-\alpha^{x_l}(a)\|$. Hence the sequence $\left(\alpha^{x_k+y_k}(a)\right)_{k\in\NM}$ is Cauchy in $\Aa$ and therefore it converges to an element of $\Hull(a)$. In much the same way, the opposite of $\xi$ is defined by $-\xi=(-x_k)_{k\in\ZM}$. Then again, it follows that $\| \alpha^{-x_k}(a)-\alpha^{-x_l}(a)\|=\| \alpha^{-x_k-x_l}\{\alpha^{x_l}(a)-\alpha^{x_k}(a)\}\|= \| \alpha^{x_k}(a)-\alpha^{x_l}(a)\|$. Hence by the same argument $-\xi\in S(a)$. This gives $S(a)$ the structure of a an abelian group.

\vspace{.1cm}

\noindent It is easy to check that the relation $\xi\stackrel{a}{\sim}\eta$ defined by $\lim_{k\rightarrow\infty} \|\alpha^{x_k}(a)-\alpha^{y_k}(a)\|=0$ is an equivalence relation compatible with the group structure. Thus, $S(a)/\stackrel{a}{\sim}=\Omega_a$ is also an abelian group. Then, the function $d(\xi,\eta)=\lim_{k\rightarrow\infty} \|\alpha^{x_k}(a)-\alpha^{y_k}(a)\|$, defined on $S(a)\times S(a)$ satisfies (i) $d(\xi,\eta) \geq 0$, (ii) $d(\xi,\eta)= d(\eta,\xi)$, (iii) $d(\xi,\eta) \leq d(\xi,\zeta)+d(\zeta,\eta)$ for all $\zeta\in S(a)$, (iv) $d(\xi,\eta)=0$ if and only if $\xi\stackrel{a}{\sim}\eta$. In particular, $d$ defines a metric on $\Omega_a$. 

\vspace{.1cm}

\noindent The map $\phi:S(a)\to \Hull(a) $ defined by $\phi(\xi) = \lim_{k\rightarrow \infty} \alpha^{x_k}(a)$ satisfies (i) it is onto, by definition of the Hull, (ii) $\phi(\xi)=\phi(\eta)$ if and only if $\xi\stackrel{a}{\sim}\eta$, (iii) $\|\phi(\xi)-\phi(\eta)\|= d(\xi,\eta)$. Hence $\phi$ defines an isometry from $\Omega_a$ onto $\Hull(a)$, showing that $\Omega_a$ is compact indeed. Moreover, if $n\in\ZM$ it follows from the definition that $\phi(\xi+\hn)= \alpha^n(\phi(\xi))$. In particular $d(\xi+\hn,\eta+\hn)=d(\xi,\eta)$ for all $n\in\ZM$. In addition if $\eta=(y_k)_{k\in\NM}\in S(a)$ it follows that $d(\eta,\hy_k)\rightarrow 0$ as $k\rightarrow\infty$. Consequently the metric $d$ is invariant by translation on $\Omega_a$ so that the addition in $\Omega_a$ is continuous. Moreover, it is simple to show that $d(-\xi,-\eta)= d(\xi,\eta)$. Hence $(\Omega_a,d)$ is a topological compact abelian group. Then the equivalence class of $\hn$ for $n=1$ defines an element $\omega_0\in\Omega_a$ such that the set $\{n\omega_0\,;\, n\in\ZM\}$ is dense in $\Omega_a$. At last, the relation $\phi(\xi+\hn)= \alpha^n(\phi(\xi))$ shows that the map $n\in\ZM\mapsto \alpha^n(a)\in\Hull(a)$ can be continued in a unique way as a map $\omega\in\Omega_a\mapsto \alpha^{\omega}(a)\in\Hull(a)$ throughout the quotient map defined by $\phi$.
\hfill $\Box$

\begin{coro}
\label{met10.cor-Fap}
If $F\in\Aa$ is a finite subset of $\alpha$-almost periodic elements. Then there is compact abelian group $\Omega_F$ with an element $\omega_F$, called the unit, such that (i) the $\ZM$-action defined by $\omega\in\Omega_F\mapsto \omega+\omega_F\in\Omega_F$ is minimal, (ii) for any $a\in F$, the map $n\in\ZM\mapsto \alpha^n(a)\in\Aa$ can be continued in a unique way as a continuous surjective map $\phi_{F,a}:\omega\in\Omega_F \mapsto \alpha^\omega(a)\in\Hull(a)$ such that $\phi_{F,a}(\omega_F)= \alpha(a)$. Moreover, if $F\subset G$ there is a canonical surjective continuous group homomorphism $\phi_{G,F}:\Omega_G\mapsto \Omega_F$ which sends $\omega_G$ onto $\omega_F$.
\end{coro}
\noindent  {\bf Proof: } The same construction can be done if $S(a)$ is replaced by $S(F)= \bigcap_{a\in F} S(a)$. Since $F$ is finite this intersection is not empty. For indeed, let choose $a,b\in F$ and $\xi_a\in S(a)$. Then the sequence $\alpha^{x_k}(a)$ converges in $\Hull(a)$. Since $\Hull(b)$ is compact there is a subsequence $\xi_{a,b}$ of $\xi_a$ such that $\alpha^{x_{k_l}}(b)$ is convergent, namely $\xi_{a,b}\in S(a)\cap S(b)$. By recursion on the set $F$, this leads to the existence of elements in $S(F)$. The rest of the proof is similar to the proof of Proposition~\ref {met10.prop-hull} and can be done by inspection. 
\hfill $\Box$

\begin{rem}
\label{met10.rem-countap}
{\em If $F$ is infinite but countable, the set $S(F)$ can be shown to be nonempty by using a diagonal procedure to extract subsequences. If $F$ is infinite, the same construction works provided the sets $S(F)$ are replaced by a family of nets.
}
\hfill $\Box$
\end{rem}

\begin{proposi}
\label{met10.prop-approxap}
The set of $\alpha$-almost periodic elements of $\Aa$ is norm closed.
\end{proposi}

\noindent {\bf Proof: } If $a =\lim_j a_j$ with $a_j$ $\alpha$-almost periodic for all $j$'s. Then $\|\alpha^n(a)-\alpha^n(a_j)\|=\|a-a_j\|$. Let $\epsilon >0$ and let $j$ be such that $\|a-a_j\|\leq \epsilon/3$. The family of open balls $B_n=\{b\in\Aa\,;\, \|b-\alpha^n(a_j)\| <\epsilon/3\}$ cover the orbit of $a_j$ thus they cover $\Hull(a_j)$. Since $a_j$ is almost periodic, it follows that there is a finite subset $J_j\subset \ZM$ such that the balls $(B_k)_{k\in J_j}$ cover $\Hull(a_j)$ as well. In particular, $\inf_{k\in J_j} \|\alpha^n(a_j)-\alpha^k(a_j)\|\leq \epsilon/3$ for all $n\in\ZM$. It follows that $\inf_{k\in J_j} \|\alpha^n(a)-\alpha^k(a)\|\leq \epsilon$ for all $n\in\ZM$. Therefore, the balls $C_n=\{b\in\Aa\,;\, \|b-\alpha^n(a)\| <\epsilon$ with $n\in J_j$ cover $\Hull(a)$ as well, showing that $\Hull(a)$ is compact, namely that $a$ is  $\alpha$-almost periodic.
\hfill $\Box$

\begin{coro}
\label{met10.cor-alphaap}
Let $\Aa$ be a separable \Cs and let $\alpha\in\Aut(\Aa)$ be almost periodic. Then there is an abelian compact group $\Omega_\alpha$ with an element $\omega_\alpha$, called the unit, such that (i) the $\ZM$-action defined by $\omega_\alpha\in\Omega_\alpha\mapsto \omega+\omega_\alpha\in\Omega_\alpha$ is minimal, (ii) for any finite subset $F\subset \Aa$ there is a surjective continuous group homomorphism $\phi_F:\Omega_\alpha\mapsto \Omega_F$ preserving the units, such that, if $F\subset G$ then $\phi_G\circ\phi_{G,F}=\phi_F$. In particular for all $a\in\Aa$, the map $\phi_{\{a\}}$ is the unique continuous extension of $n\in\ZM\mapsto \alpha^n(a)\in\Aa$.
\end{coro}

\noindent  {\bf Proof: } Since $\Aa$ is separable, there is a countable dense subset $F_\infty$. The family of groups $\{\Omega_F\,;\, F\subset F_\infty\; \mbox{\rm finite}\}$ defines an inverse limit, thanks to the restriction maps $\phi_{G,F}$. Set $\Omega_\alpha = \projlim_{F}(\Omega_F,\phi_{G,F})$. Using the Proposition~\ref{met10.prop-approxap} this group does not depend on the choice of $F_\infty$. The rest can be proved by inspection, using the properties of projective limits.
\hfill $\Box$

\vspace{.3cm}

 \subsection{Proof of Theorem~\ref{met10.th-qi}}
 \label{met10.ssect-th1}

\noindent The following result will be helpful

\begin{lemma}
\label{met10.lem-blipclosed}
Let $X=(\Aa,\Hh,D)$ be a spectral triple. Then the Lipschitz ball is norm closed.
\end{lemma}

\noindent  {\bf Proof: } Let $I$ be the set of unit vectors in the domain of $D$. For $f,g\in I$, let $p_{f,g}$ be the semi-norm on $\Aa$ defined by

$$p_{f,g}(a) =
   \left|
     \langle Df, a g\rangle - \langle f, a Dg\rangle
   \right|
$$

\noindent Clearly $p_{f,g}$ is norm-continuous. Thus $B_{f,g} = \{a\in \Aa\,;\, p_{f,g}(a) \leq 1\}$ is norm-closed. Since the domain of $D$ is dense, it follows that 

$$p(a)= \|[D,a]\| =
   \sup_{f,g\in I} p_{f,g}(a)\,.
$$

\noindent Hence $p$ is a lower semi-continuous function for the norm topology. Since $\Blip=\bigcap_{f,g\in I}B_{f,g}$ it is also norm-closed.
\hfill $\Box$ 

\vspace{.2cm}

\noindent The proof of Theorem~\ref{met10.th-qi} goes as follows. Let $G\subset \Aut(\Aa)$ have a compact closure. Then, by definition, it is quasi-isometric. In particular, for any $a\in\Cc^1(G,X)$ the map $g\in G\mapsto \|[D,g(a)]\|\in \RM_+$ is continuous. Thus it is bounded. This shows that $\Cc^1(G,X)=\Cc_b^1(G,X)$ and that $G$ is equicontinuous. 

\vspace{.1cm}

\noindent Conversely, let $G\subset \Aut(\Aa)$ be equicontinuous. Then let $\Blip(G,X)$ be defined by

$$\Blip(G,X) = 
   \{a\in \Cc_b^1(G,X)\,;\, 
      \sup_{g\in G}\|[D,g(a)]\|\leq 1
   \}\,.
$$

\noindent It follows that $\Blip(G,X)$ is $G$-invariant. By an argument similar to the proof of Lemma~\ref{met10.lem-blipclosed} above, it follows that it is norm-closed as well. By construction, for any element $a\in \Cc_b^1(G,X)$ there is $\lambda >0$ such that $\lambda a\in\Blip(G,X)$. In particular, by density, $\Blip(G,X)$ is total in $\Aa$. If $B= \Blip(G,X)$ let $B^{B}$ be equipped with the product topology. By definition of the norm-pointwise convergence, the map $\phi: g\in G\mapsto (g(a))_{a\in \Blip(G,X)}\in B^{B}$ is a homeomorphism from $G$ onto its image. For $g\in G$ let $S_g: B^{B}\mapsto B^{B}$ be the shift defined by 

$$(S_g \xi)(a) = g(\xi(a))\,,
\hspace{2cm}
  \xi\in B^{B}, a\in B=\Blip(G,X)\,.
$$

\noindent For indeed such a $\xi$ satisfies $\xi(a)\in B$ for all $a\in \Blip(G,X)$, therefore any automorphism $g\in G$ transforms $B$ into itself, so that $S_g$ is well defined and continuous. Then, by construction $\phi(gg') = S_g\circ\phi(g')$ whenever $g,g'\in G$.  

\vspace{.1cm}

\noindent Let $\theta:\Aa\mapsto \Aa/\Aa_D'$ be the quotient map. If $\Aa_D'$ is trivial, it is left invariant by $G$, so that $G$ also acts by isometries on $\Aa/\Aa_D'$. In particular $\theta$ maps $B$ into a closed subset $K$ of $\theta(\Blip)$, which, since $X=(\Aa,\Hh, D)$ is a spectral metric space, is compact in $\Aa/\Aa_D'$. Let then $\hth$ be the map induced by $\theta$ on the product space $B^{B}$, that is

$$\hth(\xi)(a) = \theta(\xi(a))\,,
\hspace{2cm}
  \xi\in B^{B}, a\in B=\Blip(G,X)\,.
$$

\noindent This map is continuous. Therefore, the map $\psi=\hth\circ\phi:G\to K^{B}$ is continuous. Thanks to the Tychonov Theorem, $K^{B}$ is compact. In addition, it will be proved that $\psi$ is one-to-one. If so, $\psi$ defines a homeomorphism from $G$ into its image and its closure is compact. 

\vspace{.1cm}

\noindent Thanks to the group property, it is enough to show that $\psi(g)=\psi(id)$ implies that $g=id$. Now $\psi(g)=\psi(id)$ if and only if $g(a) =a\,\bmod\CM\id_{\Aa}$, namely if and only if there is $\eta(a)\in \CM$ such that $g(a)=a +\eta(a)\id_{\Aa}$ for all $a\in\Blip(G,X)$. Since the map $a\mapsto (g(a)-a)$ is linear and bounded from $\Aa$ into $\Aa$, it follows that $\eta$ is a bounded linear form on $\Aa$. Since any $g\in G$ is a $\ast$-automorphism of $\Aa$, it follows that $\eta(a^\ast) = \overline{\eta(a)}$.  In addition, $g(\id_{\Aa})=\id_{\Aa}$ implies $\eta(\id_{\Aa})=0$. Moreover, $g(a^\ast a)= g(a)^\ast g(a)$ implies that

$$\eta(a)a^\ast + \overline{\eta(a)} a =
   \{\eta(a^\ast a) -|\eta(a)|^2\}\id_{\Aa}
$$

\noindent Applying $\eta$ on both sides shows that $2|\eta(a)|^2=0$ so that $\eta=0$ and $g=id$.

\vspace{.5cm}

\section{Spectral Triples for Equicontinuous Group Actions}
\label{met10.sect-stEqui}

\noindent This Section is devoted to the proof of Theorem~\ref{met10.th-qiCross}. As in the statement of this Theorem, $\Aa$ is a unital separable \Cs and $\alpha\in\Aut(\Aa)$. Moreover, $u$ will denote the generator of the $\ZM$-action defined by $\alpha$ in the crossed product algebra $\XP$. For convenience, $C_c(\ZM,\Aa)$ will denote the dense subalgebra of the crossed product $\Aa\rtimes_\alpha \ZM$ of elements of the form $\sum_{n\in\ZM} b_n u^n$ with only finitely many nonzero $b_n$'s.

\vspace{.3cm}

 \subsection{The Dual Action}
 \label{met10.ssect-dualact}

\noindent The dual action on a crossed product was defined by Connes \cite{Co73} and Takesaki \cite{TaM73}, in the case of von Neumann algebras, and by Takai \cite{TaH74,TaH75}, for \CsS. Let $\TM= \RM/2\pi\ZM$ denote the dual group of $\ZM$. For $k\in\TM$, let $\eta_k\in \Aut(\Aa\rtimes_\alpha \ZM)$ be the dual action of $\TM$, defined by

\begin{equation}
\label{met10.eq-duact}
\eta_k(a)=a\,,\;\; \mbox{\rm for $a\in\Aa$}\,,
   \hspace{2cm}
    \eta_k(u) = e^{\imath k}\, u\,.
\end{equation}

\noindent It will be convenient to introduce extra tools of calculations at this point. The generator of the dual action is defined, for $b\in C_c(\ZM,\Aa)$, by

$$\partial b= 
   \frac{d\eta_k(b)}{dk}\upharpoonright_{k=0} =
    \sum_{l\in\ZM} \imath l\,b_l\,u^l\,.
$$

\noindent Thanks to the group property of the dual action this is a $\ast$-derivation, in particular, it satisfies the Leibniz rule $\partial(bc) = \partial(b) c+ b\partial(c)$ and $\partial(b^\ast)= \left(\partial b\right)^\ast$. Another important tool is the {\em canonical conditional expectation} $\EM:\XP\mapsto \Aa$ defined by

\begin{equation}
\label{met10.eq-cexp}
\EM(b) = 
   \int_{\TM} \frac{dk}{2\pi}\; \eta_k(b)= b_0\,.
\end{equation}

\noindent It is a completely positive linear map such that $\|\EM(b)\|\leq \|b\|$ and satisfies $\EM(aba')= a\EM(b)a'$ if both $a,a'\in\Aa$. It follows immediately from the inverse Fourier transform formula, that, if $b\in C_c(\ZM,\Aa)$ then

$$b_l = \EM(bu^{-l})\,.
$$

\noindent This latter formula shows that $b_l$ is defined for any element in $\XP$, thanks to eq.~(\ref{met10.eq-cexp}). Additional important tools, such as the use of the Fejer kernel, can be found in the Appendix (Section~\ref{met10.sect-Xprod}).

\vspace{.3cm}

 \subsection{A Perturbation Result}
 \label{met10.ssect-pert}

\noindent Let $Y=(\XP,\Kk,D)$ be a (compact) spectral metric space based on the crossed product algebra $\XP$. As it is well-known the crossed product $\XP$ is not a total invariant characterizing the conjugacy class of $\alpha$. For instance, if $C$ is the Cantor set and if $\phi$ is a minimal homeomorphism of $C$, then, modulo isomorphism, the algebra $\Cc(C)\rtimes_\phi\ZM$ only gives the $\ZM$-actions orbit equivalent to $\phi$ \cite{GPS95,GPS10}. In a similar way, there is no reason why this spectral metric space could distinguish the unitary $u$ from other elements of $\XP$. So, in general, it  cannot be expected that $u\in\Cc^1(Y)$ (see Definition~\ref{met10.def-smoothact}). However, the following result shows that, provided the dual action is quasi-isometric and at the cost of perturbing $\alpha$ by an inner automorphism, it is possible to assume that $u\in\Cc^1(Y)$

\begin{lemma}
\label{met10.lem-eqsms}
Let $Y=(\XP,\Kk,D)$ be a spectral metric space based on $\XP$ so that $u$ denotes the generator of $\alpha$ in $\XP$. Let $Y$ be such that the dual action is quasi-isometric. Then 

\noindent (i) for every $0 <\epsilon <1$ there is a unitary element $w\in\Aa$ such that $\|w-\id\|< O(\epsilon)$ and such that $wu\in\Cc^1(\TM,Y)$;

\noindent (ii) if $\beta=\Adj(w)\circ\alpha$ then $\beta$ is quasi-isometric in $X=(\Aa,\Kk,\rho,D)$.

\noindent (iii) in particular, the two Banach spaces $\Aa\rtimes_{\alpha}\ZM$ and $\Aa\rtimes_{\beta}\ZM$ coincide and they are isomorphic as \CsS. Moreover, $Y'=(\Aa\rtimes_{\beta}\ZM, \Kk,D)$ is a spectral metric space.
\end{lemma}

\noindent {\bf Proof: } (i) Since $\Cc^1=\Cc^1(\TM,Y)$ is dense, by hypothesis, it follows that, given any $0<\epsilon<1$, there is $b\in\Cc^1$ such that $\|b-u\|\leq \epsilon$. In particular, $b$ is invertible in $\XP$. Therefore

$$(b_1-\id)u =
   \int_{\TM} \frac{dk}{2\pi}\;e^{-\imath k}\,
     \eta_k(b-u)
$$

\noindent showing that $\|b_1u-u\|=\|b_1-\id\|\leq \epsilon$. Hence $b_1$ is invertible in $\Aa$. Moreover, $b\in\Cc^1$ implies $b_1 u\in\Cc^1$, since the dual action is quasi-isometric. In particular $b_1b_1^\ast= b_1 uu^\ast b_1^\ast$ is a positive invertible element of $\Aa$ in $\Cc^1$, thus in $\Cc^1(X)$. Thanks to the invariance of $\Cc^1$ by holomorphic functional calculus, it follows that $w=(b_1b_1^\ast)^{-1/2}b_1$ is a unitary in $\Aa$ and that $\|w-1\|\leq C\epsilon$ for some $C>0$. In addition, $wu\in\Cc^1$ by construction.

\vspace{.1cm}

\noindent (ii) Define $\beta\in\Aut(\Aa)$ by $\beta(a)=w\alpha(a)w^{-1}$. Hence $wua(wu)^{-1}= w\alpha(a)w^{-1}=\beta(a)$ for $a\in\Aa$. Then, a dense subalgebra of $\Aa\rtimes_\beta\ZM$ is the set of elements of the form $b=\sum_{|n|\leq N} a_n (wu)^n$, for some $N\in\NM$, with $a_n\in\Cc^1(X)$. Since $(wu)^n= w\alpha(w)\cdots \alpha^{n-1}(w) u^n$, it follows that such elements belong to $\XP$. Since $w$ is invertible, the same argument shows that the two \Css $\XP$ and $\Aa\rtimes_\beta \ZM$ are equal as a set and as a Banach space. Since $wu\in\Cc^1(\TM,Y)$ it follows that such $b$'s are all in $\Cc^1$ so that $\beta$ leaves $\Cc^1(X)$ invariant.

\vspace{.1cm}

\noindent (iii) Let $\phi: \Aa\rtimes_{\beta}\ZM\to \Aa\rtimes_{\alpha}\ZM$ be defined by $\phi(a)_n= a_n w\alpha(w)\cdots \alpha^{n-1}(w)$. Then an elementary calculation shows that $\phi(a\ast_\beta b)= \phi(a)\ast_\alpha \phi(b)$ for $a,b\in\Aa\rtimes_\beta\ZM$, where $\ast_\alpha, \ast_\beta$ denote the product in each of these algebras. Similarly, $\phi(a^\ast) = \phi(a)^\ast$. Thus $\phi$ is a $\ast$-homomorphism. Exchanging the roles of $\alpha$ and $\beta$ shows that $\phi$ is also a $\ast$-isomorphism, thus an isometry. Since, as sets, the two algebras are the same, the spectral triple $Y'=(\Aa\rtimes_\beta, \Kk,D)$ is also a spectral metric space, but now, the generator of the $\ZM$ action belongs to $\Cc^1(Y')$.
\hfill $\Box$

\vspace{.3cm}

 \subsection{Invariant Metrics}
 \label{met10.ssect-invmet}

\noindent In this Section, $Y=(\XP,\Kk,D)$ is a (compact) spectral metric space based on the crossed product algebra $\XP$. It is assumed that (i) the dual action is quasi-isometric, (ii) $u^{-1}[D,u]$ is bounded and commutes with the elements of $\Aa$. Then

\begin{lemma}
\label{met10.lem-smsubspace}
The family $X=(\Aa,\Kk,D)$ defines a compact spectral metric space.
\end{lemma}

\noindent {\bf Proof: } (i) Since the dual action is equicontinuous, it follows that, for $b\in\Cc^1(\TM,Y)$, the integral $\int_{\TM} [D,\rho\circ\eta_k(b)]dk$  converges in norm, leading to $\|[D,\rho\circ\EM(b)]\|<\infty$. Namely the projection of the dense subalgebra $\Cc^1(\TM,Y)\subset \XP$ onto $\Aa$, which is dense, is contained in the algebra $\Cc^1(\Aa,D)$. Hence  $X=(\Aa,\Kk,D)$ is a spectral triple.

\vspace{.1cm}

\noindent (ii) If $a\in\Aa$ belongs to the $D$-commutant for $\Aa$ it belongs to the $D$-commutant for $\XP$, and thus it is a multiple of the unit. Since the $D$-commutants are the same, the Lipschitz ball for $\Aa$ is contained in the Lipschitz ball for $\XP$, so that it is compact when seen in the quotient space $\Aa/\Aa_{\hD}'$. 
\hfill $\Box$

\begin{lemma}
\label{met10.lem-smeq}
The automorphism group generated by $\alpha$ is equicontinuous.
\end{lemma}

\noindent {\bf Proof: } Since $\gamma_u=u^{-1}[D,u]$ is bounded it follows that $u\in\Cc^1(Y)$. Therefore, $u^{-1}Du$ is well defined and equal to $D+\gamma_u$. Thus, since $\gamma_u$ commutes with the elements of $\Aa$, $[u^{-1}Du,a]=[D,a]$ for $a\in\Cc^1(X)$. Therefore, $[D,\alpha^{-1}(a)]=u^{-1}[D,a]u$. By recursion on $n\in\NM$, it follows that $[D,\alpha^{-n}(a)]=u^{-n}[D,a]u^n$ for all $n\in\NM$. In much the same way, $\gamma_{u^{-1}}= u[D,u^{-1}]=u\gamma_u u^{-1}$ also commutes with the elements of $\Aa$, since $[\gamma_{u^{-1}},a]= u[\gamma_u,\alpha^{-1}(a)]u^{-1}=0$. Therefore, $[D,\alpha^{-n}(a)]=u^{-n}[D,a]u^n$ for all $n\in\ZM$. Hence $\|[D,\alpha^{-n}(a)]\|=\|[D,a]\|$ for all $n\in\NM$ and $\alpha$ defines an isometric action of $\ZM$ on $X$. Hence it is equicontinuous.
\hfill $\Box$

\vspace{.3cm}

 \subsection{The Regular Representation on $\Aa\rtimes_\alpha\ZM$}
 \label{met10.ssect-strcros}

\noindent In this section the direct part of  Theorem~\ref{met10.th-qiCross} is considered. Let $X=(\Aa,\Hh,D)$ be a compact spectral metric space. Let $\alpha\in \Aut(\Aa)$ generate an equicontinuous group. In particular, $\Cc_b^1(\ZM,X) = \{a\in\Aa\,;\, \sup_{n\in\NM}\|[D,\alpha^n(a)]\|<\infty\}$ is dense in $\Aa$. Then let $\Kk=\Hh\otimes \ell^2(\ZM)\otimes \CM^2$. Hence, an element $f\in\Kk$ can be written as $f= (f_n)_{n\in\ZM}$ with $f_n= (f_{n+},f_{n-})$ where $f_{n,\pm}\in\Hh$. It will be convenient to introduce the {\em Pauli matrices}

\begin{equation}
\label{met10.eq-Pauli}
\sigma_1 = \left[
              \begin{array}{cc}
              0 & 1\\
              1 & 0\\
              \end{array}
             \right]\,,
\hspace{2cm}
\sigma_2 = \left[
              \begin{array}{cc}
              0 & -\imath\\
              \imath & 0\\
              \end{array}
             \right]\,,
\hspace{2cm}
\sigma_3 = \left[
              \begin{array}{cc}
              1 & 0\\
              0 & -1\\
              \end{array}
             \right]\,. 
\end{equation}

\noindent For the purpose of this work, the main property of these matrices is their anticommutation rule

$$\sigma_i\sigma_j+\sigma_j\sigma_i= 2\delta_{ij}
$$

\noindent The left regular representation $\hpi$ of $\Aa$ induced by $\alpha$ is defined by

\begin{equation}
\label{met10.eq-lrr}
\left(\hpi(a) f\right)_n =
   \alpha^{-n}(a)f_{n}\,,
\hspace{2cm}
\left( \hu f\right)_n = f_{n-1}\,.
\end{equation}

\noindent Hence $\hu$ is a unitary operator on $\Kk$ which satisfies $\hu\hpi(a)\hu^{-1} = \hpi\circ \alpha(a)$. In particular this representation extends to the crossed product $\Aa\rtimes_\alpha \ZM$ by setting $\hpi(\sum_n a_n u^n)= \sum_n \hpi(a_n) \hu^n$. A new Dirac operator is defined as follows

\begin{equation}
\label{met10.eq-dirlrr}
\left(\hD f\right)_n =
 \left[
   \begin{array}{cc}
    0 & D-\imath n\\
    D+\imath n & 0
   \end{array}
 \right]\,f_n = \left(
  D\otimes \sigma_1 + n \id_{\Hh}\otimes \sigma_2
                 \right) f_n\,,
\end{equation}

\noindent This operator is selfadjoint by construction. Moreover it has compact resolvent, because if $D$ is written in its spectral decomposition as $D=\sum_{k\in\ZM} \lambda_k \Pi_k$ where $\{ \cdots, < \lambda_k <\lambda_{k+1}<\cdots\}$ is the set of eigenvalues and $\Pi_k$ is the corresponding eigenprojection, then

$$\hD\upharpoonright_n = \sum_{k\in\ZM}
   \left[
   \begin{array}{cc}
    0 & (\lambda_k-\imath n)\Pi_k\\
   (\lambda_k +\imath n)\Pi_k & 0
   \end{array}
 \right]\,,
$$

\noindent so that the eigenvalues of $\hD$ are given by the $\pm\sqrt{\lambda_k^2+n^2}$ and the eigenprojections by 

$$\Pi_{k,n}^\pm =\frac{1}{2}\,
   \left[
   \begin{array}{cc}
\Pi_k & \pm \frac{(\lambda_k-\imath n)}{\sqrt{\lambda_k^2+n^2}}\Pi_k\\
\pm \frac{(\lambda_k-\imath n)}{\sqrt{\lambda_k^2+n^2}}\Pi_k & \Pi_k
   \end{array}
 \right]\otimes P_n\,,
$$

\noindent if $P_n$ denotes the rank one projection acting on $\ell^2(\ZM)$ which selects the $n$-th component. The eigenprojections are indeed finite dimensional and the inverse of the eigenvalues converge to zero as either $k$ or $n$ go to infinity. Hence $\hD$ has compact resolvent. Consequently, $\hD$ is a Dirac operator on $\Kk$. On the other hand

\begin{equation}
\label{met10.eq-comlrr}
\left([\hD, \hu]f\right)_n =
   \left[
   \begin{array}{cc}
    0 & -\imath\\
   \imath & 0
   \end{array}
 \right] \left(\hu f\right)_n\,,
\hspace{1cm}
  \left([\hD, \hpi(a)]f\right)_n =
   \left[
   \begin{array}{cc}
    0 & 1\\
   1 & 0
   \end{array}
 \right]\,[D,\alpha^{-n}(a)]f_n\,.
\end{equation}

\noindent It follows that $\|[\hD,\hpi(a)]\|<\infty$ for $a\in \Cc_b^1(\ZM,X)$. Let $\Bb_c$ be the subalgebra of $\Aa\rtimes_\alpha \ZM$ made of elements of the form $\sum_n a_n u^n$ where $a_n\in\Cc_b^1(\ZM,X)$ and $a_n=0$ but for a finite number of indices. Then, since $\Cc_b^1(\ZM,X)$ is dense in $\Aa$, $\Bb_c$ is dense in the crossed product and, thanks to eq.~(\ref{met10.eq-comlrr}), $b\in\Bb_c \Rightarrow \|[\hD,\hpi(b)]\| <\infty$. Hence $Y=(\Aa\rtimes_\alpha\ZM,\Kk,\hD)$ is a spectral triple.

\begin{defini}
\label{met10.def-rrX}
The spectral triple $Y=(\Aa\rtimes_\alpha\ZM,\Kk,\hD)$ will be called the regular representation of the metric dynamical system $(X,\alpha)$. It will be denoted by $Y=X\rtimes_\alpha\ZM$.
\end{defini}

\begin{lemma}
\label{met10.lem-rrprop}
The regular representation $Y=X\rtimes_\alpha\ZM$ has the following properties

(i) the dual action is isometric,

(ii) $\hu^{-1}[\hD,\hu]$ commutes with the elements of $\Aa$.
\end{lemma}

\noindent  {\bf Proof: } For $k\in\TM$, let $v_k$ be the unitary operator defined on $\Kk$ by

$$(v_kf)_n = e^{\imath nk} \,f_n\,,
\hspace{2cm}
    f\in\Kk\,.
$$

\noindent It is a standard result that $k\in\TM\mapsto v_k\in\Bb(\Kk)$ defines a strongly continuous unitary representation of the torus $\TM$. Moreover, an elementary calculation shows that

$$v_k \hpi(a) v_k^{-1} = \hpi(a)\,,\; a\in\Aa\,,
\hspace{2cm}
   v_k\hu v_k^{-1}= e^{\imath k}\,\hu\,,
  \hspace{2cm}
    v_k\hD v_k^{-1}= \hD\,.
$$

\noindent It follows immediately that $v_k$ implements the dual action on $\XP$ and that it defines an isometry on $Y$.

\vspace{.1cm}

\noindent An elementary calculation gives

$$\hu^{-1}[\hD,\hu] = 
   \id_{\Kk}\otimes \id_{\ell^2(\ZM)}\otimes \sigma_2\,,
$$

\noindent showing that it commutes with the representation of $\XP$, thus with the elements of $\Aa$. 
\hfill $\Box$

\vspace{.3cm}

 \subsection{$\Aa\rtimes_\alpha\ZM$ as a Spectral Metric Space}
 \label{met10.ssect-qmetcross}

\noindent It remains to prove that $Y$ is a spectral metric space. Namely, it requires to prove that the $\hD$-commutant of $\Aa\rtimes_\alpha \ZM$ is trivial and that its Lipschitz ball is pre-compact. It is worth reminding \cite{Ped79} that, if $\Aa$ is unital with unit $\id_\Aa$, then so is $\Aa\rtimes_\alpha \ZM$ and its unit is given by $\id_{\Aa\rtimes_\alpha \ZM} = \sum b_lu^l$ with $b_l=\delta_{l,0}\id_\Aa$. On the other hand, if $\Aa$ is nonunital, then so is $\Aa\rtimes_\alpha \ZM$.

\begin{lemma}
\label{met10.lem-hDcommut}
The $\hD$-commutant of $\Aa\rtimes_\alpha \ZM$ is trivial.
\end{lemma}

\noindent  {\bf Proof: } Let $b$ be an element of $\Cc_b^1(\TM,Y)$. Thanks to Lemma~\ref{met10.lem-rrprop}, it follows that $[\hD, \hpi(b)]=0$ if and only if $[\hD, \hpi\circ\eta_k(b)]=0$ for all $k\in\TM$. Therefore, using the Fourier decomposition {\em w.r.t.} $k$ gives

$$[\hD,\hpi(b_l)] +
      l\hpi(b_l) \id_{\Hh}\otimes \sigma_2=0
 \hspace{2cm}
    \forall l\in\ZM\,.
$$

\noindent Using eq.~(\ref{met10.eq-comlrr}), this implies 

$$[D,\alpha^{-n}(b_l)]\otimes \sigma_1 +
      l\,\alpha^{-n}(b_l) \id_{\Hh}\otimes \sigma_2=0
 \hspace{2cm}
    \forall l,n\in\ZM\,.
$$

\noindent Since the matrices $\id_2,\sigma_1,\sigma_2,\sigma_3\in M_2(\CM)$ are linearly independent, it follows that

$$\alpha^{-n}(b_l)=0\,,
   \hspace{1cm}\mbox{\rm if $l\neq 0$}\hspace{2cm}
    [D,\alpha^{-n}(b_0)] =0
 \hspace{1cm}
    \forall n\in\ZM\,.
$$

\noindent Since $X=(\Aa,\Hh,\pi,D)$ is a spectral metric space, it follows that $b_l=0$ for $l\neq 0$ and $b_0\in\Aa_D$, namely $b$ is a multiple of the identity in $\Aa\rtimes_\alpha\ZM$.
\hfill $\Box$

\begin{lemma}
\label{met10.lem-eqmetric}
The Connes metric on the original spectral metric space $X$ and the one induced on $X$ by $Y$ are equivalent.
\end{lemma}

\noindent  {\bf Proof: } Since $\alpha$ generates an equicontinuous group, it follows that the Connes metric associated with $X=(\Aa,\Hh,D)$, denoted here by $d_X$, satisfies  eq.~\ref{met10.eq-biLip} uniformly, namely, there is $K>0$ such that, for any pair of states $\rho,\omega$ on $\Aa$,

$$\frac{1}{K} d_X(\rho,\omega) \leq
   d_X(\rho\circ\alpha^n,\omega\circ\alpha^n )\leq 
    Kd_X(\rho,\omega)\,,
\hspace{2cm}
 \forall n\in\ZM\,.
$$

\noindent On the other hand the Connes metric $d_Y$ on the state space of $\Aa$ induced by the regular representation is given by 

$$d_Y(\rho,\omega)=
   \sup \{\rho(a)-\omega(a)\,;\, \sup_{n\in\ZM}\|[D,\alpha^n(a)]\|\leq 1\}\,.
$$

\noindent It follows that $d_Y(\rho,\omega)\leq d_X(\rho,\omega)$. By equicontinuity there is $K>0$ such that $\|[D,a]\|\leq 1\,\Rightarrow\, \sup_n \|[D,\alpha^n(a)]\|\leq K$, showing that $d_X(\rho,\omega)\leq K\,d_Y(\rho,\omega)$.
\hfill $\Box$

\vspace{.2cm}

\noindent At last, the following holds, finishing the proof of Theorem~\ref{met10.th-qiCross}

\begin{proposi}
\label{met10.prop-lipcomp}
The Lipschitz ball of the spectral triple $Y=(\Aa\rtimes_\alpha\ZM,\Kk,\hpi,\hD)$ is compact modulo the $\hD$-commutant.
\end{proposi}

\noindent The proof of this Proposition is technically demanding and is the subject of Section~\ref{met10.sect-Xprod}.

\vspace{.5cm}

\section{The Spectral Metric Bundle}
\label{met10.sect-smb}

\noindent This Section is devoted to the construction of a spectral metric bundle over the crossed product algebra $\XP$ whenever the automorphism $\alpha$ is not equicontinuous. Some examples of such situations are described in Section~\ref{met10.ssect-neqact}. To do so, the construction of Connes-Moscovici will be the main intuitive guideline \cite{CM95}. Namely, let $M$ be a compact manifold of dimension $d$. If $G$ is a group of diffeomorphisms of $M$, if $G$ is not amenable, there may not be any $G$-invariant measure on $M$. This can be cured if $M$ is replaced by its {\em frame bundle} $\ps$, because there is always on $\ps$ a measure that is invariant by all diffeomorphisms. Then if the concept of {\em measure} is replaced by the one of {\em Riemannian metric}, something similar happens, namely there may be no $G$-invariant Riemannian metric on $M$. However if $M$ is replaced by the {\em metric bundle} $\ms$, made of pairs $(x,g)$ where $x\in M$ and $g$ is a Euclidean metric on the tangent space $T_x M$, then the group of diffeomorphisms of $M$ acts on $\ms$ and there is an invariant Riemannian metric on it. The construction of the Riemannian metric is tautological, provided the set of Euclidean metrics on $\RM^d$ is endowed itself with a Riemannian metric as well. In Section~\ref{met10.ssect-eucmetr} below, this problem is described in more details. If now the group $G$ is smaller than the family of all diffeomorphisms, and if $M$ is endowed with a Riemannian metric, it is sufficient to restrict the metric bundle to the $G$-orbit of this metric.

\vspace{.1cm}

\noindent When passing from a manifold to a noncommutative one, the manifold is replaced by a \Cs $\Aa$. If the manifold is Riemannian, the metric will be described by a spectral triple $(\Aa,\Hh, D)$. Let $G$ be a group of $\ast$-automorphism of $\Aa$. The first step will be to identify a new \Cs $\Bb$ which is liable to play the role of the set of continuous functions on the metric bundle, vanishing at infinity. Since the metric bundle admits itself a metric, the spectral triple over $\Aa$ should induce a spectral triple over $\Bb$. The construction will be made in Section~\ref{met10.ssect-crossp}. It should be such that $G$ acts on this new spectral triple in a way such that the regular representation described in Section~\ref{met10.ssect-strcros} should work.

\vspace{.3cm}

 \subsection{Examples of Non-Equicontinuous Actions}
 \label{met10.ssect-neqact}

\begin{exam}[{\bf Arnold's cat map} \cite{AA67} and a {\bf parabolic map}]
\label{krfeb10.exam-arnold}
{\em Let $\Aa=\Cc(\TM^2)$ and let $\Aa_1=\Cc^1(\TM^2)$ be the dense $\ast$-subalgebra of continuously differentiable functions on the two-torus. The Hilbert space will be $L^2(\TM^2)\otimes \CM^2$. $\Aa$ acts by multiplication on the diagonal. The Dirac operator is given by

$$D = \left[
       \begin{array}{cc}
        0 & -\imath\partial_1-\partial_2\\
        -\imath\partial_1+\partial_2 & 0\\
       \end{array}
      \right]\,.
$$

\noindent Then $D^2= -(\partial_1^2+\partial_2^2)$ is the Laplacian and so the resolvent of $D$ is compact. The commutant of $D$ in $\Aa$ are the constant function. The norm $\|a\|+\|[D,a]\|$ is nothing but the Lipschitz norm on $\Aa_1$ and the Lipschitz ball is indeed compact in $\Aa$ by the Arzel\`a-Ascoli theorem. Thus $\alpha$ is an automorphism of $\Aa$ such that $\sup_n \|[D,\alpha^n(a)\| <\infty$ if and only if $\{\alpha^n(a)\,;\, n\in\ZM\}$ has a compact uniform closure, namely if and only if $\alpha$ is almost periodic. In particular, if

\begin{equation}
\label{met10.eq-cat}
M= \left[
      \begin{array}{cc}
       2 & 1\\
       1 & 1\\
      \end{array}
     \right]\,,
\end{equation}

\noindent then the automorphism $\alpha (a)(x) = a(M^{-1}x)$ is the {\em Arnold cat map} which is uniformly hyperbolic. Then $M$ can be decomposed into

$$M = \frac{\sqrt{5}+1}{2} \pi_+ + \frac{\sqrt{5}-1}{2} \pi_-\,,
$$

\noindent where $\pi_\pm$ are two orthogonal one-dimensional projections in the set of $2\times 2$ matrices, orthogonal to each other. Then $\nabla a\circ M^{-n} = M^{-n}\nabla(a)\circ M^{-n}$. It follows that, if $\sigma = (\sigma_1,\sigma_2)$, 

$$[D,\alpha^n(a)] =
   \frac{\sqrt{5}+1}{2} \sigma\cdot \pi_+ \alpha^n(\nabla(a)) + 
    \frac{\sqrt{5}-1}{2} \sigma\cdot \pi_- \alpha^n(\nabla(a))
$$

\noindent so that $\|[D,\alpha^n(a)]\|$ diverges as $|n|\rightarrow \infty$, because $(\sqrt{5}-1)/2= ((\sqrt{5}+1)/2)^{-1}$. If now $M$ is replaced by

\begin{equation}
\label{met10.eq-parab}
N= \left[
      \begin{array}{cc}
       1 & 1\\
       0 & 1\\
      \end{array}
     \right]\,,
\end{equation}

\noindent it gives a parabolic diffeomorphism of the torus. This map occurs in some problems related to the discrete Heisenberg group. Similarly $N^n$ is unbounded, even though it increases only linearly in $n$. But the corresponding automorphism is not equicontinuous.
}
\hfill $\Box$
\end{exam}

\begin{exam}[{\bf Bilateral shift}]
\label{met10.exam-shift}
{\em This example is illustrative of what happens for tiling spaces \cite{PB09}. Let $\Xi=\{0,1\}^{\ZM}$ be the Cantor set\footnote{The disconnected Greek letter $\Xi$ is an appropriate choice to denote a completely disconnected set.}, 
represented as the set of sequences of $0$'s and $1$'s indexed by the integers $\ZM$. The {\em bilateral shift} is the map $S:\Xi\to\Xi$ defined by

$$(Sx)_n=x_{n-1}\,,
\hspace{2cm}
   x=(x_n)_{n\in\ZM}\in\Xi\,.
$$

\noindent The construction of a metric and of a spectral triple follows the strategy in \cite{PB09} that is appropriate for tiling spaces as well. Namely, let $W_n=\{0,1\}^{2n+1}$ denote the set of words of length $2n+1$ and let $W$ denote the {\em dictionary}, namely, the disjoint union of the $W_n$'s. Given $w\in W_n$ then $2n+1=|w|$ is called the length of the word $w$. Then $\Xi(w)$ denotes the {\em acceptance zone} of $w$, namely, the set of sequences $x=(x_n)_{n\in\ZM}$ with middle word coinciding with $w$

$$x_{-n+k}=w_k\,,
   \hspace{2cm}
k=0,1,\cdots ,2n 
$$

\noindent An ultrametric $d$ is defined by demanding that the distance between two sequences $x,y$ be given by $2^{-|w|}$ if $w$ is the longest middle word common to $x$ and $y$. This metric can be recovered from the following continuous family of spectral metric spaces $X_{\tau}= (\Cc(\Xi), \Hh,\pi_\tau, D)$, where $\Hh=\ell^2(W)\otimes \CM^2$ and $D$ is defined by

$$D\psi(w) =2^{|w|}
     \left[
      \begin{array}{cc}
       0 & 1\\
       1 & 0\\
      \end{array}
     \right]\psi(w)\,,
\hspace{2cm}
  \psi=(\psi(w))_{w\in W}\in \Hh\,,\;\; \psi(w)\in\CM^2\,.
$$

\noindent A family of representations of $\Aa=\Cc(\Xi)$ should be used instead of only one, in order to recover the metric from the Connes procedure \cite{PB09}. Namely, a {\em choice} is a map $\tau:W\to \Xi\times \Xi$ with the following property: (i) if $\tau(w)=(x_w,y_w)$ then both $x_w,y_w$ belong to $\Xi(w)$, (ii) $d(x_w,y_w)=2^{-|w|}$, namely $w$ is the largest middle word common to $x_w,y_w$. The set of choices is also a Cantor set denoted by $\Upsilon$. Given a choice $\tau$, a representation can be defined by

$$\pi_{\tau}(f)\psi(w)=
     \left[
      \begin{array}{cc}
       f(x_w) & 0\\
       0 & f(y_w)\\
      \end{array}
     \right]\psi(w)\,,
\hspace{2cm}
  \psi\in \Hh\,.
$$

\noindent The Connes metric of the field $X_\tau$ of such spectral triples coincides with the previous metric \cite{PB09}, namely

$$d_C(x,y) = 
   \sup\left\{
      |f(x)-f(y)|\,;\,
       \sup_{\tau\in\Upsilon} \|[D,\pi_\tau(f)]\|\leq 1
       \right\} =d(x,y)\,.
$$

\noindent Here, the family $X=\left(X_\tau\right)_{\tau\in\Upsilon}$, not the individual element, defines a spectral metric space. As a matter of fact, the bilateral shift is not equicontinuous, especially because the metric $d$ defined above is not shift invariant. This model is the fundamental model in ergodic theory describing what is called today a {\em Smale space}.

\vspace{.1cm}

\noindent A similar conclusion holds if $\Xi$ is replaced by the tiling space of the Fibonacci sequence, obtained from the substitution $0\mapsto 01\,,\, 1\mapsto 0$. The tiling space is the set of all sequences sharing the same dictionary as the Fibonacci sequence. A similar construction for a metric and a spectral triple can be made and leads to a same conclusion for the shift.
}
\hfill $\Box$
\end{exam}

\vspace{.3cm}

 \subsection{The Space of Euclidean Metrics}
 \label{met10.ssect-eucmetr}

\noindent A Euclidean metric on $\RM^d$ is given by a $d\times d$ positive invertible matrix $Q$ will real coefficients. Namely, the corresponding metric is given by $g(x,y) = \langle x| Q y\rangle$ if $ \langle x| y\rangle = \sum_{i=1}^d x_iy_i$ denotes the usual dot product. Hence the space of metrics $\ms_d$ on $\RM^d$ is the interior of the positive cone on $M_d(\RM)$. Such a matrix $Q$ can be parametrized by an element of $GL_d$ as follows

$$\Lambda\in GL_d(\RM) \mapsto \Lambda\,\Lambda^t \in \ms_d\,.
$$ 

\noindent This map is onto, but not one-to-one, because the stabilizer of the metric $Q=\id$ is the subgroup $O_d\subset GL_d$ so that the previous map becomes one-to-one if seen as defined on the symmetric space $GL_d/O_d$. Hence $\ms_d$ can be seen as the manifold $GL_d/O_d$. The advantage of such an identification is that, since it is a symmetric space, the action of $GL_d$ defines a Riemannian metric on $\ms_d$ once an $O_d$-invariant Euclidean metric is defined on the tangent space to $\ms_d$ at $Q=\id$. This tangent space is nothing but the set of $d\times d$ real symmetric matrices on which the Hilbert-Schmidt inner product defines a canonical metric. The transport of the metric on the tangent space at $Q$ can be made as follows: if $Q=\Lambda\,\Lambda^t$ then $dQ=\Lambda\,H\,\Lambda^t$ for some $H=H^t\in M_d(\RM)$. Thus

\begin{equation}
\label{met10.eq-metmet}
ds^2 = \frac{1}{d} \TR (H^2)=
    \frac{1}{d} \TR \left(Q^{-1}dQQ^{-1} dQ \right)\,.
\end{equation}

\noindent The following result is straightforward and the proof will be left to the reader

\begin{proposi}
\label{met10.prop-geo}
Let $\ms_d$ be endowed with the Riemannian metric given by eq.~(\ref{met10.eq-metmet}). Then a geodesic is given by a map $s\in [0,1]\mapsto Q(s)$ as a solution of the equation

\begin{equation}
\label{met10.eq-geodseq}
\frac{d^2Q}{ds^2} = 
   \frac{dQ}{ds}\,Q^{-1}\frac{dQ}{ds} +
    \frac{
    \TR\left(
        \left(Q^{-1}\frac{dQ}{ds}\right)^3
        \right)-
    \TR\left(Q^{-1}\frac{dQ}{ds}Q^{-1}\frac{d^2Q}{ds^2}\right)
     }
{\TR\left(\left(Q^{-1}\frac{dQ}{ds}\right)^2\right)}\;
   \frac{dQ}{ds}\,.
\end{equation}
\end{proposi}

\noindent A special case is provided by the curve $\gamma_H$ given by  $s\in [s_0,s_1]\mapsto Q(s) = e^{sH}\in\ms_d$ for $H=H^t$. Since $Q^{-1} dQ/ds= H$ and since $H$ commutes with $Q(s)$, a simple calculation shows that this curve satisfies the geodesic equation~(\ref{met10.eq-geodseq}). Moreover

$$dL^2 = \frac{1}{d} \TR \left(Q^{-1}dQQ^{-1} dQ \right)=
       \frac{1}{d} \TR \left(H^2 \right)\,ds^2\,.
$$

\noindent In particular,

$$L(\gamma_H) = |s_1-s_0|
    \frac{1}{\sqrt{d}}
     \left(\TR \left(H^2 \right)\right)^{1/2}\,.
$$

\noindent An example of such a situation consists in evaluating the Riemannian distance between the two metrics given by $Q^m$ and $Q^n$ whenever $m\neq n$. Writing $Q=e^H$ this gives

\begin{equation}
\label{met10.eq-distnm}
d(Q^m,Q^n) = |m-n|
    \frac{1}{\sqrt{d}}
     \left(\TR \left(\ln{Q}^2 \right)\right)^{1/2}\,.
\end{equation}

\noindent This case corresponds to the situation for the Arnold cat map (see eq.~\ref{met10.eq-cat}), namely, along the orbit, $Q^{(n)}=(M^\ast)^n M^n=M^{2n}$ so that the previous formula applies with $Q=M^2$.

\vspace{.3cm}

 \subsection{Spectral Metric Spaces Based on $\ZM$}
 \label{met10.ssect-zm}

\noindent Since the orbit of the original metric in the metric bundle is identified with $\ZM$ it makes sense to consider first the possible structures of spectral metric spaces based on $\ZM$. The natural \Cs associated with this space is $c_0(\ZM)$, namely the set of sequences of complex numbers, vanishing at infinity, endowed with the uniform norm

$$a\in c_0(\ZM)\,\;\;\; a=\left(a_n\right)_{n\in\ZM}
\hspace{1cm} \Rightarrow \hspace{1cm}
  \lim_{n\rightarrow \infty} a_n =0\,,
\hspace{2cm}
  \|a\|= \sup_{n\in\ZM} |a_n|\,.
$$

\noindent The main feature of this algebra is that it is not unital. Characterizing the metrics equivalent to the weak$^\ast$ topology on the state space of $c_0(\ZM)$ can be done by using the results obtained in \cite{La07} (see Result~\ref{met10.res-latr}). In this case it is well known that the dual space of $c_0(\ZM)$ can be identified with the space $\ell^1(\ZM)$ of absolutely summable sequences $\eta= \left(\eta_n\right)_{n\in\ZM}$ with $\|\eta\|_{\ell^1}=\sum_{n\in\ZM} |\eta_n| <\infty$. Then $\eta\in\ell^1(\ZM)$ is positive if and only $\eta_n \geq 0$ for all $n$'s. If $C_1$ denotes the set of positive elements with norm less than or equal to one, then $C_1$ is weak$^\ast$ closed and compact (Banach-Alaoglu Theorem). A positive linear form $\eta$  is a state if and only if $\sum_{n\in\ZM} \eta_n =1$, so that the state space is nothing but the set $\MG_1(\ZM)$ of probability measures on $\ZM$ equipped with the weak$^\ast$ topology. It is well-known, and easy to prove, that $\MG_1(\ZM)$ is a weak$^\ast$ dense $G_\delta$ subset of $C_1$. In addition, the Prokhorov Theorem \cite{Pr56} gives a characterization of compact subsets in $\MG_1(\ZM)$: a weak$^\ast$ closed subset $K\subset \MG_1(\ZM)$ is compact if and only if, for all $\epsilon >0$, there is an integer $N$ such that $\sum_{|n| >N} \eta_n \leq \epsilon$ for all $\eta\in K$. This condition is called {\em tightness} in probability theory \cite{Bi99} and expresses the fact that a sequence in $K$ cannot escape at infinity. 

\vspace{.1cm}

\noindent In order to characterize the metrics on $\MG_1(\ZM)$ defining the weak$^\ast$ topology, the notion of {\em weak-uniform} topology ($\wu$) on $c_0(\ZM)$  will be used. It was introduced by Latr\'emoli\`ere \cite{La07} on separable \Css with or without units. Namely, if $\KG$ denotes the set of compact subsets in $\MG_1(\ZM)$, given any compact subset $K\in \KG$, let $p_K$ be the seminorm on $c_0$ defined by

$$p_K(a) = 
   \sup_{\omega\in K} |\omega(a)|\,.
$$

\noindent The $\wu$-topology is the locally convex topology defined by the family of seminorms $(p_K)_{K\in\KG}$ (if $c_0(\ZM)$ where unital, this topology would be the same as the norm topology). A subset $F\subset c_0(\ZM)$ is called {\em totally bounded} whenever, for any $\epsilon >0$ and any compact subset $K\subset \MG_1(\ZM)$, there is a finite set $A\subset F$ such that the open $p_K$-balls centered at the elements of $A$ with radius $\epsilon$ cover $F$. A crucial result in \cite{La07} is that $F$ is totally bounded if and only if there is a strictly positive element $h\in c_0(\ZM)$ such that $hFh$ is pre-compact (see Result~\ref{met10.res-latr}). Note that here $c_0$ is abelian so that $hFh=Fh^2$. It follows that there is a need to characterize the compact subsets of $c_0(\ZM)$. The following Lemma shows that they are exactly the bounded closed subsets {\em equicontinuous at infinity}, more precisely

\begin{lemma}
\label{met10.lem-compc0}
A closed subset $F\subset c_0(\ZM)$ is compact if and only if (i) it is bounded, (ii) for all $\epsilon >0$ there is $N\in\NM$ such that for $|n|\geq N$ $|a_n|\leq \epsilon$ for all $a\in F$. 
\end{lemma}

\noindent  {\bf Proof: } If $F\subset c_0(\ZM)$ is compact, then, for any $\epsilon >0$ there is a finite family $\{a^{(1)},\cdots,a^{(m)}\}$ in $F$ so that the open balls $B(a^{(i)},\epsilon/2)$ cover $F$. In particular, if $a\in F$, then there is an $1\leq i\leq m$ such that $\|a-a^{(i)}\|\leq \epsilon/2$ so that $\|a\|\leq \max_i \|a^{(i)}\|+\epsilon/2$. Hence $F$ is bounded. Moreover, since $a^{(i)}\in c_0(\ZM)$, there is a natural integer $N_i$ such that for $|n|\geq N_i$, $|a^{(i)}_n| < \epsilon/2$. If $N=\max\{N_i\}$ it follows that for any $a\in F$, there is $1\leq i\leq m$ such that $\|a-a^{(i)}\|< \epsilon/2$, implying that $|a_n|<\epsilon$ if $|n|\geq N$.

\vspace{.1cm}

\noindent Conversely, let $F$ be a closed bounded subset of $c_0(\ZM)$ such that for all $\epsilon >0$ there is $N\in\NM$ such that for $|n|\geq N$ $|a_n|\leq \epsilon$ for all $a\in F$. Then let $K$ denotes the disc in $\CM$ centered at the origin with radius $\sup_{a\in F} \|a\|$. Hence any sequence $a\in F$ satisfies $a_n\in K$ for all $n$. Therefore $F$ can be identified with a subset of the infinite product $\hK= K^{\ZM}$. Let $F_w$ denote the image of $F$ in $\hK$, that will be equipped with the product topology. Therefore, thanks to the Tychonov theorem, $F_w$ is pre-compact. The injection $\phi: F\to F_w$ is onto, by construction and continuous. It is sufficient to show that the inverse is continuous to prove that $F$ is compact. Since $\hK$ is metrizable, it is sufficient to show that if $(a^{(j)})_{j\in\NM}$ denotes a convergent sequence in $F_w$, then it converges uniformly. The convergence in the product topology means that for each $n\in\ZM$ the limit $a_n=\lim_{j\rightarrow \infty} a^{(j)}_n$ exists. Since $a^{(j)}\in F_w$, let $N$ be such that $|a^{(j)}_n| <\epsilon/2$ for all $j$'s and $|n|\geq N$. It follows that (i) $|a_n|\leq \epsilon/2$ whenever $|n\geq N$, showing that $a\in c_0(\ZM)$, and that (ii) $ \sup_{|n|\geq N}|a_n-a^{(j)}_n| \leq \epsilon$. Since the sequence $a^{(j)}$ converges pointwise, it follows that there is $J\in\NM$ such that for $j\geq J$ and $|n|<N$ $|a_n-a^{(j)}_n|\leq\epsilon$. Hence $\|a-a^{(j)}\|\leq\epsilon$ as well. Thus $\lim_{j\rightarrow \infty} \|a-a^{(j)}\| =0$. Since $F$ is closed, it follows that $a\in F$. 
\hfill $\Box$

\begin{coro}
\label{met10.cor-unBalwucompact}
Given any strictly positive element $h$ in $c_0(\ZM)$, the unit ball $B$ satisfies $hBh$ is compact.
\end{coro}

\noindent  {\bf Proof: } (i) $h$ is a positive element of $\Aa= c_0(\ZM)$ if and only if $h=(h_n)_{n\in\ZM}$ with $h_n \geq 0$ for all $n$'s. It is strictly positive if and only if $h_n >0$ for all $n$. For indeed, if so, let $a\in\Aa$ be any element. For $N\in\NM$ let $a^N$ be the element obtained from $a$ by changing $a_n$ into zero if $|n|>N$. Then $\lim_{N\rightarrow\infty} \|a-a^N\|=0$. Since $h_n>0$ for all $n$ the element $b^N= a^N/h^2$ is well defined in $c_0(\ZM)$, therefore $a^N= h^2 b^N \in h\Aa h$. Conversely, if there is $n\in\ZM$ such that $h_n=0$ then any element $a\in h\Aa h$ satisfies $a_n=0$ so that $h\Aa h$ is not norm dense in $\Aa$.

\vspace{.1cm}

\noindent (ii) If $h\in c_0(\ZM)$ is strictly positive, then if $\|a\|\leq 1$ it follows that $hah$ satisfies $|(hah)_n| = |a_n| h_n^2 \leq h_n^2$. Hence $hBh$ is obviously closed, bounded and equicontinuous at infinity, thus, thanks to Lemma~\ref{met10.lem-compc0}, is norm compact.
\hfill $\Box$

\vspace{.2cm}

\noindent Given a metric $d_\ZM$ on $\ZM$, a spectral triple can be constructed as follows. A natural faithful representation of $c_0(\ZM)$ on $\ell^2(\ZM)$ is given by the pointwise multiplication

$$(af)_n= a_n\,f_n\,,
\hspace{2cm}
   a\in c_0(\ZM)\,,\, f\in \ell^2(\ZM)\,.
$$

\noindent The discrete analog of a derivation is the finite difference $(\partial f)_n= f_{n}-f_{n-1}$, which is not self-adjoint. If the metric changes from site to site, the corresponding operator should rather be

$$(\nabla f)_n =
   \frac{f_n-f_{n-1}}{\imath d_\ZM(n,n-1)}\,,
\hspace{1cm} \Rightarrow\hspace{1cm}
   (\nabla^\ast f)_n =
   \frac{f_{n+1}}{\imath d_\ZM(n+1,n)}-
    \frac{f_{n}}{\imath d_\ZM(n,n-1)} \,,
$$

\noindent In this definition only the distance between consecutive points is required. It is natural to define a Dirac operator by changing $\ell^2(\ZM)$ into the spinor space $\ell^2(\ZM)\otimes\CM^2$ and set

$$D_{\ZM} =\left[
      \begin{array}{cc}
      0 & \nabla\\
      \nabla^{\ast} & 0\\
     \end{array}
    \right]\,=
   \frac{\sigma_1+\imath \sigma_2}{2} \;\nabla +
    \frac{\sigma_1-\imath \sigma_2}{2} \;\nabla^\ast\,.
$$

\noindent This gives a symmetric operator, but its resolvent is not compact in general. However, adding to $D_\ZM$ an operator commuting to $c_0(\ZM)$ may transform it into a compact one as seen below. Before stating the result, it will be convenient to introduce the following operators on $\ell^2(\ZM)$ and on $\ell^2(\ZM)\otimes \CM^2$:

\begin{equation}
\label{met10.eq-X}
(Xf)_n =\left[
      \begin{array}{cc}
      c_n c_{n+1} & 0\\
      0 & -c_n c_{n+1}\\
     \end{array}
    \right]f_n=
     c_n c_{n+1}\sigma_3 f_n\,,
\end{equation}

\noindent where the sequence $c= (c_n)_{n\in\ZM}$ is required to satisfy

$$\lim_{|n|\rightarrow\infty} c_n =+\infty,,
   \hspace{2cm}
    c_n\,\delta_n \geq 1\,.
$$

\noindent with $\delta_n = d_\ZM(n,n-1)$. The result is the following

\begin{lemma}
\label{met10.lem-dirZcomp}
The operator $D_\lambda=D_\ZM+\lambda X$ is selfadjoint with compact resolvent for any $\lambda \in \RM\setminus\{0\}$.
\end{lemma}

\noindent  {\bf Proof: } It should be noted that $X$ has compact resolvent, since $(z\id-X)^{-1}$ is a diagonal operator with eigenvalues 
$(z+\pm c_n c_{n+1})^{-1}$ which converges to zero as $|n|\rightarrow \infty$. It will be proved that $(z-\lambda X)^{-1}D_\ZM$ is compact as well. If so, the resolvent equation implies

$$(z-\lambda X-D_\ZM)^{-1} =
   (z-\lambda X)^{-1}+ (z-\lambda X)^{-1}D_\ZM(z-\lambda X-D_\ZM)^{-1}\,,
$$

\noindent leading to the result. Now, for $z=\imath$

$$\left((\imath-\lambda X)^{-1}D_\ZM f\right)_n=\left[
      \begin{array}{c}
      (\imath-c_n c_{n+1})^{-1} \delta_n^{-1} (f_n^--f_{n-1}^-)\\
      (\imath+c_n c_{n+1})^{-1} 
       (\delta_n^{-1}f_n^+-\delta_{n+1}^{-1}f_{n+1}^+)
     \end{array}
    \right]\,,
$$

\noindent where $f_n^{\pm}$ denote the two components of $f_n\in\CM^2$. The hypothesis made on $c_n$ shows that 

$$|\imath-c_n c_{n+1}|^{-2} \delta_n^{-2} =
   \frac{1}{c_n^2\delta_n^2 c_{n+1}^2 +\delta_n^2}\;\leq\; 
    c_{n+1}^{-2} \stackrel{|n|\uparrow \infty}{\longrightarrow} 0
$$

\noindent A similar estimate holds for the other component. This shows that 
$(z-\lambda X-D_\ZM)^{-1}$ is compact indeed.
\hfill $\Box$

\vspace{.2cm}

\noindent The Connes metric is now defined by

$$d_C(n,m) =
   \sup\{|a_n-a_m|\,;\, \|[D_\lambda,a]\|\leq 1\}
$$

\noindent Following Gromov \cite{Gr99} the path metric associated with a metric $d_\ZM$ is given by 

$$d_p(m,n) = \sum_{k=m+1}^n d_\ZM(k,k-1)\,,
  \hspace{2cm}
    \mbox{\rm for}\;\; m<n\,.
$$

\begin{lemma}
\label{met10.lem-Connes}
Let $d_\ZM$ be a metric on $\ZM$. Then, the Connes metric associated with $D_\lambda$ coincides with the path metric $d_p$ associated with $d_\ZM$.
\end{lemma}

\noindent {\bf Proof: } An elementary calculation shows that

$$[D_\lambda,a] =\left[
      \begin{array}{cc}
      0 & [\nabla,a]\\
      -[\nabla, a^{\ast}]^{\ast} & 0\\
     \end{array}
    \right]\,.
$$

\noindent It follows from then that $\|[D_\lambda,a]\|= \max\{\|[\nabla,a]\|\,,\,\|[\nabla,a^{\ast}]\| \}$ so that

\begin{equation}
\label{met10.eq-C1c0}
\|[D_\lambda,a]\|=
   \sup_n\left(\frac{|a_n-a_{n+1}|}{d_\ZM (n,n+1)}\right)
\end{equation}

\noindent Therefore $\|[D_\lambda,a]\|\leq 1$ implies $|a_n-a_m| \leq \sum_{j=m+1}^n d_\ZM(j,j-1)$ whenever $m<n$. Let now $m<n$ and define $b\in c_0(\ZM)$ by $a_j=0$ if $j\leq m$, $b_{m+k}= \sum_{j=1}^k d_\ZM(m+j,m+j-1)$ for $1\leq k\leq n-m$, then $b_{n+k}= b_n -\sum_{j=n-m-k}^{n-m} d_\ZM(m+j,m+j-1)$ if $1\leq k\leq n-m$, $b_l=0$ if $l\geq 2n-m$. Then, by construction $\|[D_\lambda,b]\| \leq 1$ and $|b_m-b_n|= \sum_{j=1}^{n-m} d_\ZM(m+j,m+j-1)$.
\hfill $\Box$

\vspace{.2cm}

\noindent At this point the eq.~\ref{met10.eq-C1c0} gives an important clue for the spectral triple to be a spectral metric space. In view of the Lemma~\ref{met10.lem-compc0}, the Lipschitz ball is compact only if $\lim_{n\rightarrow \infty} d_\ZM(n,n+1) =0$. More precisely

\begin{lemma}
\label{met10.lem-smsZM}
Let $d_\ZM$ be a metric on $\ZM$. Then the triple $(c_0(\ZM),\ell^2(\ZM)\otimes \CM^2, D_\lambda)$ is a spectral metric space if and only if $\sum_{n\in\ZM} d_\ZM(n,n-1) <\infty$. In the latter case $\Blip$ is norm compact.
\end{lemma}

\noindent {\bf Proof: } (i) Clearly, the representation of $c_0(\ZM)$ is non degenerate and the $D_\lambda$-commutant is trivial. Thus the only condition that ought to be checked is the total $\wu$-boundedness of the Lipschitz ball. It is sufficient to show that $\Blip$ is bounded if and only if $C=\sum_{n\in\ZM} d_\ZM(n,n-1) <\infty$, in which case it is norm compact. 

\vspace{.1cm}

\noindent (ii) Since $X$ commutes with the elements of $c_0(\ZM)$, it follows that $[D_\lambda,a]=[D_\ZM,a]$. The Lipschitz ball is therefore the set of $a\in c_0(\ZM)$ such that $|a_n-a_{n-1}|\leq d_\ZM(n,n-1)$ for all $n$. If the metric is summable, then, for $m\geq n$ and $a\in\Blip$ it follows that $|a_m-a_n|\leq \sum_{j=n+1}^m d_\ZM(j,j-1)$. If $m\rightarrow +\infty$ this gives $|a_n| \leq \sum_{j=n+1}^{+\infty} d_\ZM(j,j-1)\leq C$. In particular $\|a\|\leq C$. In much the same way, letting $n\rightarrow -\infty$ gives $|a_m| \leq \sum_{-\infty}^{m} d_\ZM(j,j-1) \leq C$. This shows that $\Blip$ is bounded and equicontinuous at infinity in the norm topology. Thanks to Lemma~\ref{met10.lem-compc0}, $\Blip$ is norm compact.

\vspace{.1cm}

\noindent (iii) If now $\sum_{j\in\ZM} d_\ZM(j,j-1)$ diverges, for any integer $N$ let $a^N$ be the sequence defined as follows: (a) $a_n^N = 0$ if $ n\leq 0$, (b) $a_n^N= \sum_{j=1}^n d_\ZM(j,j-1)$ if $1\leq n\leq N$ and (c) $a_n^N= \max\left\{\left(\sum_{j=1}^N d_\ZM(j,j-1) - \sum_{j=N+1}^n d_\ZM(j,j-1)\right), 0\right\}$ if $n\geq N+1$. By construction $a^N\in \Blip$ and $\|a^N\|= \sum_{j=1}^N d_\ZM(j,j-1)$ wich diverges as $N\rightarrow \infty$. Hence $\Blip$ is unbounded, and the spectral triple cannot be a spectral metric space.
\hfill $\Box$

\begin{rem}
\label{met10.rem-oneptcomp}
{\em The latter result shows that such a spectral triple can only be defined as a spectral metric space if $c_0(\ZM)$ can be replaced by its one point compactification. 
}
\hfill $\Box$
\end{rem}

\noindent The latter result suggest that, in order to recover the weak$^\ast$ topology from a spectral triple without requiring a one-point compactification of the \CS, it is necessary to use a different Dirac operator. One possibility consists in using a field of such operators, instead of only one, in the spirit of what was done in \cite{PB09}. A candidate is the following: let $\nabla$ be replaced by 

$$(\nabla_r f)_n =
    \frac{f_n-f_{n-r}}{d_\ZM(n,n-r)}\,,
$$

\noindent and let $D_r$ be the corresponding Dirac operator. Changing $c_{n+1}$ into $c_{n+r}$ in the definition of $X$ (see eq.~(\ref{met10.eq-X})), leads to a diagonal operator $X_r$, so that $D_{r,\lambda}:= D_r+\lambda X_r$ has compact resolvent for all $r\in\NM$ and $\lambda\in\RM\setminus\{0\}$. Then the Connes distance associated with the sequence $D_\lambda=\left(D_{r,\lambda}\right)_{r\geq 1}$ will now be defined as follows

$$d_C(\rho,\omega) =
   \sup\{|\rho(a)-\omega(a)| \,;\, 
    \sup_{r\geq 1}\|[D_{r,\lambda},a]\|\leq 1
       \}
$$

\noindent A convenient way to represent this family in a unique spectral triple consists in looking at the direct sum. This is achieved by changing the Hilbert space $\ell^2(\ZM)$ into $\ell^2(\ZM\times \NMo)$ where $\NMo$ is the set of nonzero natural integers. In order to get a Dirac operator, it will be necessary to increase the set of Dirac matrices. In this specific case it seems that at least four of them, denoted by $\gamma_1, \cdots, \gamma_4$, are needed. This gives the following construction: let $\Ee$ be a $4$-dimensional Hilbert space on which the Dirac matrices are represented and let $\Kk$ be the Hilbert space $\Kk=\ell^2(\ZM\times \NMo)\otimes \Ee$. A vector in $\Kk$ can be seen as a double sequence $f=(f_{n,r})_{(n,r)\in\ZM\times \NMo}$ with $f_{n,r}\in\Ee$. Then the algebra $c_0(\ZM)$ is represented by 

$$(af)_{n,r} = a_n f_{n,r}\,.
$$

\noindent The Dirac operator will be now given as before in two steps. First the part giving the various finite differences with the family of operator $\nabla_r$ and then a diagonal operator to insure that the resolvent will be compact. A solution to this problem consists in defining the following operators

$$(\nabla f)_{n,r} = \frac{f_{n,r}-f_{n-r,r}}{d_\ZM(n,n-r)}\,,
   \hspace{1cm}
    (Xf)_{n,r} = c_n c_{n+r} \gamma_3 f_{n,r}\,,
   \hspace{1cm}
     (Rf)_{n,r} = r \gamma_4 f_{n,r}\,.
$$

\noindent Then the Dirac operator on $\Kk$ will be defined by

\begin{equation}
\label{met10.eq-DirZM}
D_{\Kk} = 
   \frac{\gamma_1+\imath \gamma_2}{2} \;\nabla +
    \frac{\gamma_1-\imath \gamma_2}{2} \;\nabla^\ast +
     \lambda (X+R)\,. 
\end{equation}

\noindent It is tedious but straightforward to check that $D_{\Kk}$ has compact resolvent and that 

$$\|[D_{\Kk},a]\| =
   \sup_{n\neq m} \frac{|a_m-a_n|}{d_\ZM(n,m)}\,.
$$

\noindent The main result is the following

\begin{proposi}
\label{met10.prop-LipField}
Let $d_\ZM$ be a metric on $\ZM$. Then $(c_0(\ZM),\ell^2(\ZM\times\NMo)\otimes \Ee, D_{\Kk})$ is a spectral triple. It is a spectral metric space if and only if the metric $d_\ZM$ is bounded.
\end{proposi}

\noindent {\bf Proof: } (i) Checking that the representation of $c_0$ is non degenerate is clear. Checking that $\Cc^1$ is dense is also easy. The only technical point consists in checking that $D_\Kk$ has compact resolvent. Since the proof follows the lines of the proof of Lemma~\ref{met10.lem-dirZcomp}, it will be left to the reader.

\vspace{.1cm}

\noindent (ii) The second part consists in checking whether the Lipschitz ball is $\wu$-totally bounded or not. The same argument as in the proof of Lemma~\ref{met10.lem-smsZM} shows that, if $d_\ZM$ is not bounded, there is a sequence $(a^N)_{N\in\NM}$ in $\Blip$ such that $\|a^N\|\rightarrow \infty$ so that $\Blip$ is not bounded. Indeed it is sufficient to take (a) $a_n^N = 0$ if $ n\leq 0$, (b) $a_n^N= d_\ZM(0,n)$ if $1\leq n\leq N$ and (c) $a_n^N= \max\{(d_\ZM(0,N) - d_\ZM(N,n)), 0\}$ if $n\geq N+1$.

\vspace{.1cm}

\noindent (iii) Conversely if $d_\ZM$ is bounded then, any $a\in\Blip$ satisfies $|a_n|\leq \limsup_{m\rightarrow \infty} |a_n-a_m| \leq \limsup_{m\rightarrow \infty} d_\ZM(n,m) <\infty$. Therefore $\Blip $ is bounded in norm, so that, thanks to Corollary~\ref{met10.cor-unBalwucompact}, $\Blip$ is $\wu$ totally bounded.
\hfill $\Box$

\begin{rem}
\label{met10.rem-bdmetric}
{\em It is easy to define on $\ZM$ bounded metrics that are shift invariant. An example is $d_\ZM(n,m) = \tanh{(|n-m|)}$. More generally if $d_1 > d_2 > \cdots d_n > 0$ is a decreasing sequence of positive numbers such that $\sum_j d_j$ converges, then $d_\ZM(n,m) = \sum_{j=1}^{|n-m|} d_j$ if $n\neq m$, defines a shift invariant bounded metric on $\ZM$. 
}
\hfill $\Box$
\end{rem}

\noindent The main question is now to understand what happens if the metric is unbounded. The answer to this question is provided by the following result about the Connes metric

\begin{lemma}
\label{met10.lem-Connes2}
Let $d_\ZM$ be a bounded metric on $\ZM$. Then the Connes distance associated with the spectral metric space $(c_0(\ZM),\ell^2(\ZM\rtimes\NMo)\otimes \Ee, D_{\Kk})$ is exactly the Wasserstein distance of order $1$ on $\MG_1(\ZM)$ associated with $d_\ZM$. 
\end{lemma}

\noindent  {\bf Proof: } (i) Given two probabilities $\rho, \omega \in\MG_1(\ZM)$, let $\MG(\rho,\omega)$ be the set of probabilities on $\ZM^2$ with marginals given by $\rho$ and by $\omega$ respectively. Practically, $\mu\in\MG(\rho,\omega)$ is a double sequence $(\mu_{m,n})_{(m,n)\in\ZM^2}$ with $\mu_{m,n}\geq 0$, $\sum_{n\in\ZM}\mu_{m,n} =\rho_m$ and $\sum_{m\in\ZM}\mu_{m,n}=\omega_n$ for all $m,n$'s. The Wasserstein distance of order $p\geq 1$ is defined as

$$W_p(\rho,\omega)^p =
   \inf_{\mu\in\MG(\rho,\omega)}
     \sum_{(m,n)\in\ZM^2} \mu_{m,n}\; d_\ZM(m,n)^p\,.
$$

\noindent For $a\in c_0(\ZM)$, the difference $\rho(a)-\omega(a)$ can be written as $\rho(a)-\omega(a) = \sum_{(m,n)\in\ZM^2} \mu_{m,n}(a_m-a_n)$. Since $a$ belongs to the Lipschitz ball if and only if $|a_m-a_n| \leq d_\ZM(n,m)$ for all indices, then $|\rho(a)-\omega(a)|\leq \sum_{m,n} \mu_{m,n}d_\ZM(m,n)$. Minimizing over the measures $\mu$ gives $|\rho(a)-\omega(a)|\leq W_1(\rho,\omega)$. In particular the Connes distance is dominated by the Wasserstein distance $W_1$. 

\vspace{.1cm}

\noindent The opposite inequality is the content of a Theorem by Kantorovich and Rubinstein \cite{KR57}, showing that if $\rho$ and $\omega$ have finite support, then their Connes distance coincides with their Wasserstein distance.
\hfill $\Box$

\vspace{.2cm}

\noindent The Lemma~\ref{met10.lem-Connes2} gives a clue for what to do if the distance is unbounded. For indeed the Wasserstein distance is defined on a dense subset of measures, namely the probability measures $\omega\in\MG_1(\ZM)$ such that $\sum_{n\in\ZM} \omega_n \;d_\ZM(n_0,n)<\infty$ for some $n_0$. When $d_\ZM$ is unbounded this subset is weakly dense and is a countable union of weakly compact sets, by Prokhorov's Theorem. This indicates that the Connes metric might be defined for a spectral triple that is not a spectral metric space, provided the state space is replaced by a suitable closed subset. The following set of results is directly coming from \cite{Do70}. 

\vspace{.1cm}

\noindent Let $\MG_D(\ZM)$ be defined as the set of states $\omega\in\MG_1(\ZM)$ such that $\sum_{n\in\ZM} |n|\omega_n <\infty$. Clearly, $\MG_D(\ZM)$ is weakly dense in $\MG_1(\ZM)$. A subset $F\subset \MG_D(\ZM)$ will be called {\em $D$-tight} whenever for every $\epsilon >0$ there is $N\in\NMo$ such that 

$$\omega\in F 
   \hspace{1cm}\Rightarrow \hspace{1cm}
    \sum_{|n|>N} |n|\omega_n \leq \epsilon\,.
$$

\noindent The main result in \cite{Do70} is the following

\begin{proposi}[Dobrushin \cite{Do70}]
\label{met10.prop-Dob}
(i) The Wasserstein metric $W_1$ associated with the usual metric $d(m,n)=|m-n|$ on $\ZM$ is well defined on $\MG_D(\ZM)$ and makes it a complete metric space.

\noindent (ii) If $(\rho^j)_{j\in\NMo}$ is a sequence in $\MG_D(\ZM)$ that $W_1$-converges to $\rho$ then it converges weakly to $\rho$.

\noindent (iii) Any $D$-tight set is weak$^\ast$-compact. On any $D$-tight subset of $\MG_D(\ZM)$ the $W_1$-topology coincides with the weak$^\ast$ one. 
\end{proposi}

\noindent {\bf Proof: } (i) Let $\rho,\omega$ be a pair of elements in $\MG_D(\ZM)$. Then if $\mu\in \MG(\rho,\omega)$, it follows that 

$$\sum_{m,n} \mu_{m,n} |m-n| \leq 
   \sum_{m,n} \mu_{m,n} (|m|+|n|) = 
    \sum_{m} (\rho_{m}+\omega_m) |m| <\infty\,.
$$

\noindent This implies that $W_1(\rho,\omega) <\infty$. In particular, if $a\in c_0(\ZM)$ then, since $m\neq n\Rightarrow 1\leq |m-n|$, $|\rho(a)-\omega(a)| \leq \sum_{(m\neq n)\in\ZM^2} \mu_{m,n}|a_m-a_n| \leq 2\|a\| \sum_{(m\neq n)\in\ZM^2} \mu_{m,n} |m-n|$, it follows that

\begin{equation}
\label{met10.eq-weakW1}
|\rho(a)-\omega(a)| \leq 
   2\|a\| \;W_1(\rho,\omega)
\end{equation}

\noindent In particular, this proves (ii), namely that a sequence converging in the $W_1$-topology converges weakly. In addition to eq.~(\ref{met10.eq-weakW1}) the following estimate holds

\begin{equation}
\label{met10.eq-MDest}
|\sum_{n\in\ZM} (\rho_n-\omega_n)|n|\leq 
  W_1(\rho,\omega)\,.
\end{equation}

\noindent For indeed, $|\sum_{n\in\ZM} (\rho_n-\omega_n)|n|| = |\sum_{(m\neq n)\in\ZM^2} \mu_{m,n} (|m|-|n|)|\leq \sum_{(m\neq n)\in\ZM^2} \mu_{m,n}|m-n|$. Minimizing over the choice of $\mu$ gives the estimate.

\vspace{.1cm}

\noindent (ii) Let $(\rho^j)_{j\in\NMo}$ is a $W_1$-Cauchy sequence in $\MG_D(\ZM)$. It follows from eq.~(\ref{met10.eq-weakW1}) that it converges weakly to some state $\rho$. Thanks to eq.~(\ref{met10.eq-MDest}) it follows that the sequence of positive real numbers $Z_j= \sum_{n\in\ZM} \rho_n^j\,|n|$ is Cauchy and thus converges to $Z \geq 0$. Since $1 -\rho_0^j= \sum_{n\neq 0} \rho_n^j\leq \sum_{n\in\ZM} \rho_n^j\,|n|$, it follows that $Z\geq 1-\rho_0$. In addition, thanks to the weak$^\ast$-convergence, it follows that, for all $N\in\NMo$,

$$\sum_{|n|\leq N} \rho_n |n| =
   \lim_{j\rightarrow \infty}
    \sum_{|n|\leq N} \rho_n^j |n| \leq 
     \limsup_{j\rightarrow \infty} Z_j =Z\,.
$$

\noindent this proves that $\rho\in\MG_D(\ZM)$. It remains to prove that $W_1(\rho^j,\rho)$ converges to zero as well. Since $(\rho^j)_{j\in\NMo}$ is a $W_1$-Cauchy sequence, for any $\epsilon >0$, there is $J\in\NMo$ such that for $i,j\geq J, \;\;W_1(\rho^i,\rho^j)\leq \epsilon$. Therefore, for any $i,j\leq J$ let $\mu^{i,j}\in\MG(\rho^i,\rho^j)$ be such that $\sum_{m,n}\rho_{m,n}^{i,j}|m-n| \leq 2\epsilon$. Let $\mu_{m,n}^i$ be any limit point as $j\rightarrow \infty$. By construction $\sum_{m}\mu_{m,n}^{i,j}=\rho_n^j$ so that taking the limit gives $\mu^i\in\MG(\rho^i,\rho)$. Moreover $ \sum_{m,n}\rho_{m,n}^{i}|m-n| \leq 2\epsilon$ as well, showing that $W_1(\rho^i,\rho)\leq 2\epsilon)$ for $i\leq J$. Hence $\rho$ is the $W_1$-limit of the sequence. Therefore $\MG_D(\ZM)$ is complete.

\vspace{.1cm}

\noindent (iii) Let $F$ be $D$-tight. Then let $N(\epsilon)$ be the smallest integer $N$ such that $\sum_{|n|\geq N} |n|\omega_n \leq \epsilon$. Then for $l\in\NMo$ let $F_l = \{\omega\in \MG_1(\ZM)\,;\, \sum_{|n|\geq N(1/l)} |n|\omega_n \leq 1/l\}$. Clearly this set is weakly closed since $\sum_{|n|\geq N} |n|\omega_n =\sup_{M}\sum_{M\geq |n|\geq N} |n|\omega_n$. Thus $F=\bigcap_{l}F_l$ is also weakly closed. Hence $F$ is weakly compact. Let $(\rho^j)_{j\in\NMo}$ be a weakly convergent sequence in $F$ and let $\rho$ be the limit. Given $\epsilon>0$ let $N\in\NMo$ be such that $\sum_{|n| >N}\rho_n \leq \epsilon$. Then let $r_m^j= \min\{\rho_m^j,\rho_m\}$ if $|m|\leq N$ and $r_m^j=0$ for $|m|>N$. It is convenient to define the following sequences 

$$\omega_n^{j,+} = \rho_n^j - r_n^j\,,
\hspace{2cm}
  \omega_n^{j,-} = \rho_n - r_n^j\,.
$$

\noindent It follows that $\omega_n^{j,\pm} \geq 0$ and that, for $|n|\leq N$, $\omega_n^{j,\pm}=|\rho_n^j-\rho_n|_\pm= \max\{\pm(\rho_n^j-\rho_n), 0\}$. In particular, $\sum_{n\in\ZM}\omega_n^{j,+}=\sum_{n\in\ZM}\omega_n^{j,-}=C_j$. In addition, there is some $J$ such that for $j\geq J$, then $C_j =\sum_{|n|>N} \rho_n + \sum_{|n|\leq N} |\rho_n^j-\rho_n|_- \leq 2\epsilon$. Let $\mu^j\in\MG_1(\ZM^2)$ be defined as follows

$$\mu_{m,n}^j =
   r_n^j\delta_{m,n} +
    \frac{\omega_m^{j,+}\;\omega_n^{j,-} }{C_j}\,.
$$

\noindent It follows that $\mu^j\in\MG(\rho^j,\rho)$ as can be checked easily. Moreover

$$\sum_{m\neq n}\mu_{m,n}^j |m-n| =
   \frac{\sum_{m\neq n}\omega_m^{j,+}\omega_n^{j,-} |m-n|
         }{C_j}\,.
$$

\noindent The upper bound $|m-n|\leq |m|+|n|$ gives $\sum_{m\neq n}\mu_{m,n}^j |m-n| =\leq \sum_{n} (\rho_n^j+\rho_n-2r_n^j)|n|\leq \sum_{n} |\rho_n^j-\rho_n| |n|$. Since both measures belong to $F$, $\sum_{|n|>N(\epsilon)} |\rho_n^j-\rho_n| |n|\leq 2\epsilon$. On the other hand, there is $J_2>0$ such that for $j>J_2$, $\sum_{|n|\leq N(\epsilon)} |\rho_n^j-\rho_n| |n|\leq \epsilon$. Hence, for $j>J_2$, $W_1(\rho^j,\rho)\leq 3\epsilon$, showing the $\lim_{j\rightarrow \infty} W_1(\rho^j,\rho)=0$.
\hfill $\Box$

\vspace{.3cm}

 \subsection{The Metric Bundle}
 \label{met10.ssect-metbund}

\noindent Let $(A,\Hh,D)$ be a spectral triple and let $\alpha\in\Aut(A)$. If $\alpha$ is not equicontinuous, the Connes-Moscovici construction suggests to replace $\Aa$ by the \Cs of functions over the metric bundle vanishing at infinity. However, for the purpose of $\alpha$ it is sufficient to restrict this bundle to the orbit of the original metric under $\alpha$. This orbit can be identified with $\ZM$ so that the new algebra should simply be $\Bb= \Aa\otimes c_0(\ZM)$, namely the algebra of sequences $b=(b_n)_{n\in\ZM}$ with $b_n\in\Aa$ and $\lim_{n\rightarrow\infty} b_n=0$. The product is pointwise, namely $(bc)_n=b_n c_n$. Remarkably $A$ acts as a multiplier algebra if $(ab)_n= a b_n$ and $(ba)_n= b_n a$. Thanks to the detailed study made in the previous Section~\ref{met10.ssect-zm}, the following result hold

\begin{proposi}
\label{met10.prop-BbCompa}
Let $\Aa$ be a separable \Cs  admitting a strictly positive element and let $\Bb$ denote $\Aa\otimes c_0(\ZM)$. Let $e_n$ denote the evaluation map $e_n: b\in\Bb\mapsto b_n\in\Aa$.

\vspace{.1cm}

\noindent 1)- A closed subset $F\subset \Bb$ is compact if and only if the following three conditions hold

\noindent \hspace{.5cm} (a) $F$ is bounded,

\noindent \hspace{.5cm} (b) for each $n\in\ZM$, $e_n(F)\subset \Aa$ is compact,

\noindent \hspace{.5cm} (c) it is equicontinuous at infinty, namely, for each $\epsilon >0$ there is $N\in\NM$ such that

\noindent \hspace{.5cm} $b\in F\;\Rightarrow\; \|b_n\| \leq \epsilon$ for $|n| >N$.

\vspace{.1cm}

\noindent 2)- A closed subset $F\subset \Bb$ is $\wu$-totally bounded if and only if the following two conditions hold

\noindent \hspace{.5cm} (a) $F$ is bounded,

\noindent \hspace{.5cm} (b) for each $n\in\ZM$, $e_n(F)\subset \Aa$ is $\wu$-totally bounded.

\end{proposi}

\noindent {\bf Proof: } (i) Thanks to Result~\ref{met10.res-latr}, $F$ is $\wu$-totally bonded if and only if there is a strictly positive element $h\in\Bb$ such that $hFh$ is norm compact. A strictly positive element can be chosen to be given by $h_n= \delta_n h$ where $h$ is a strictly positive element of $\Aa$ and $\delta_n >0$ is a sequence of positive real numbers such that $\lim_{|n|\rightarrow \infty} \delta_n=0$. Therefore the part 2)- of the proposition results from the characterization of compact subsets.

\vspace{.1cm}

\noindent (ii) Let then $F$ be norm compact. Then the projection $F_n=e_n(F)\subset \Aa$ is compact since the projection is a $\ast$-homomorphism, thus continuous. In addition, given $\epsilon >0$, there is a finite family $\{b^1,\cdots,b^s\}\subset F$ such that the open balls of radius $\epsilon/2$ centered at the $b^j$'s cover $F$. Therefore, like in the proof of Lemma~\ref{met10.lem-compc0}, it implies that (a) $F$ is bounded, (b) for each $\epsilon >0$ there is $N\in\NM$ such that $b\in F\;\Rightarrow\; \|b_n\| \leq \epsilon$ for $|n| >N$.

\vspace{.1cm}

\noindent (iii) Conversely let $F\subset \Bb$ be closed, bounded, with $e_n(F)=F_n$ compact for all $n\in\ZM$ and such that for each $\epsilon >0$ there is $N\in\NM$ such that $b\in F\;\Rightarrow\; \|b_n\| \leq \epsilon$ for $|n| >N$. Let then $(b^j)_{j\in\NM}$ be a sequence in $F$. Since the product space $\prod_{n\in\ZM} F_n$ is compact for the product topology,  without loss of generality, it can be assumed that this sequence converges pointwise to some $b$, such that $b_n =\lim_{j} b_n^j$ for all $n$'s. The same argument as in the proof of Lemma~\ref{met10.lem-compc0} shows that $\|b-b^j\|\rightarrow 0$ so that, since $F$ is closed, $b\in F$. Hence every sequence admits a convergent subsequence, showing that $F$ is compact.
\hfill $\Box$

\vspace{.2cm}

\noindent Let now $X=(\Aa,\Hh,D)$ is a compact spectral metric space, namely $\Aa$ is unital. Let $\alpha\in\Aut(\Aa)$ be quasi-isometric. Whenever the corresponding $\ZM$-action  is not equicontinuous, the following construction overcomes the difficulty coming form the fact described in teorem~\ref{met10.th-qiCross}. First replace $\Aa$ by $\Bb=\Aa\otimes c_0(\ZM)$. It is worth remarking that $n$ labels the points in the orbit of the metric under $\alpha$. Therefore, the automorphism $\alpha$ induces on $\Bb$ an automorphim $\alpha_\ast$ defined by 

$$\left(\alpha_\ast (b)\right)_n =
    \alpha(b_{n-1})\,,
$$

\noindent which coincides with the interpretation that $\alpha$ changes also the metric, namely it moves it by one step in the orbit space indexed by $n$. Then $\alpha_\ast$ induces an automorphism $\alpha_{\ast\ast}$ on the multiplier algebra $M(\Bb)$ of $\Bb$. Since $A$ can be seen as a subalgebra of $M(\Bb)$ it is easy to check that $\alpha_{\ast\ast}(a)=\alpha(a)$ for $a\in A$.

\vspace{.1cm}

\noindent Then, using the results of Section~\ref{met10.ssect-zm} it is possible to define a spectral triple based on $\Bb$ as follows: 

\begin{enumerate}

\item the Hilbert space is the tensor product $\Kk=\Hh\otimes \ell^2(\ZM\times\NMo)\otimes \Ee$ where now, $\Ee$ support a representation of five Dirac matrices instead of four. In particular a vector $f\in\Kk$ can be seen as a double sequence $(f_{n,r})_{(n,r)\in\ZM\times \NMo}$ where now $f_{n,r}\in\Hh\otimes \Ee$.

\item the algebra $\Bb$ is represented as follows

$$ (bf)_{n,r} = \alpha^{-n}(b_n)\,f_{n,r}\,,
$$

\item a bounded shift invariant distance will be chosen; so that $d_\ZM(m,n) = d_{|m-n|}$ for some sequence $(d_j)_{j\in\NM}$ satisfying $d_0=0$ and $d_j\leq d_{j'}+d_{j"}$ for all $j$ such that $j\leq j'+j"$; for simplicity there is no loss of generality in choosing this distance so that $\sup_j d_j =1$;

\item the new Dirac operator is, using the same notation as in eq.~(\ref{met10.eq-DirZM}),

$$D_{\Bb} = 
   \frac{\gamma_1+\imath \gamma_2}{2} \;\nabla +
    \frac{\gamma_1-\imath \gamma_2}{2} \;\nabla^\ast +
     \lambda X + \gamma_5 D\,. 
$$

\noindent where now it is sufficient to take for $X$ the operator

$$(Xf)_{n,r} =
   (\gamma_3 \;n + \gamma_4 \;\frac{1}{d_r^2})f_{n,r}\,,
$$

\item the automorphism $\alpha_\ast$ is implemented by the unitary operator $u$ defined by

$$(uf)_{n,r} =
    f_{n-1,r}\,,
$$

\noindent implying that $ubu^{-1}=\alpha_\ast(b)$. 

\end{enumerate}

\noindent It is straightforward to check that the $D_{\Bb}$-commutant is trivial. In addition the Lipschitz ball is made of elements $b\in\Bb$ such that (i) $\|\alpha^{-m}b_m-\alpha^{-n}b_n\|\leq d_{|n,m|}$ for all $m,n$, and (ii) $\sup_n|[D,\alpha^{-n}b_n]\|\leq 1$. Let $\Blip$ be the Lipschitz ball for $X=(\Aa,\Hh,D)$ and let $\Blip^1$ the subset of $\Blip$ of elements of norm less that . These two conditions imply that $b_n \in \alpha^n(\Blip^1)$ for all $n\in\ZM$. Since $X$ is a compact spectral metric space, it follows that $\Blip^1$ is norm compact and therefore, thanks to Proposition~\ref{met10.prop-BbCompa}, it follows that the Lipschitz ball of the spectral triple $Y=(\Bb, \Kk, D_\Bb)$ is $\wu$-totally bounded.  The following result summarizes the properties of this triple and the proof will be left to the reader

\begin{proposi}
\label{met10.prop-BbSms}
The triple $Y=(\Bb,\Kk,D_\Bb)$ as defined above is a spectral metric space, which will be called the metric bundle over $X$. In addition, $u^{-1}[D_\Bb,u]$ is bounded and commutes with the elements of $\Bb$.
\end{proposi}

\vspace{.3cm}

 \subsection{Crossed Product Metric Bundle: Proof of Theorem~\ref{met10.th-bundle}}
 \label{met10.ssect-crossp}

\noindent This section is devoted to the proof of Theorem~\ref{met10.th-bundle}. The basic object will be the metric bundle $Y$ over $X$. That $\alpha_\ast$ defined an equicontinuous action on $Y$ is a consequence of the last result of Proposition~\ref{met10.prop-BbSms}. For indeed, like in Theorem~\ref{met10.th-qiCross}, the property $u^{-1}[D_\Bb,u]$ is bounded implies that $u\in \Cc^1$. Moreover since it commutes with the elements of $\Bb$, $u[D_\Bb,b]u^{-1}= [D_\Bb,\alpha_\ast(b)]$ showing that $\alpha$ acts equicontinuously. In particular, the algebra $\Bb$ and the operator $u$ generate a $\ast$-representation of the crossed product algebra $\Bb\rtimes_{\alpha_\ast}\ZM$ and $Y\rtimes_{\alpha_\ast}\ZM= (\Bb\rtimes_{\alpha_\ast}\ZM,\Kk,D_\Bb)$ is a spectral triple. It remains to show that the Lipschitz ball of this triple is $\wu$-totally bounded.

\vspace{.1cm}

\noindent First it will be convenient to represent also the dual action. Namely, for $k\in\TM$ let $v_k$ be the unitary operator given by

$$(v_k\;f)_{n,r} = e^{\imath kn}\; f_{n,r}\,.
$$

\noindent It is elementary to check that $v_k bv_k^{-1}= b$ is $b\in\Bb$. Moreover, $v_k u v_k^{-1} = e^{\imath k}\;u$. Hence $v_k$ implements the dual action on $\Bb\rtimes_{\alpha_\ast}\ZM$. However, this action is not equicontinuous. While $v_k$ commutes with $X$ and with $D$, it does not commute to $\nabla$. For if $\nabla_k= v_k\,\nabla \,v_k^{-1}$ then 

$$(\nabla_k\,f)_{n,r}=
   \frac{f_{n,r}-e^{\imath kr}f_{n-r,r}}{d_r}\,,
$$

\noindent so that

$$\left([\nabla_k,b]\;f\right)_{n,r}= e^{\imath kr}\;
   \frac{\alpha^{-n}(b_n)- \alpha^{-n+r}(b_{n-r})}{d_r}\;
    f_{n-r,r}
$$

\noindent This expression shows that $k\in\TM\mapsto v_k [\nabla,b] v_k^{-1}\in \Bb(\Kk)$ is strongly continuous but not norm continuous because of the dependence in $r$. If, however, $h=(h_n)_{n\in\ZM}$ is a strictly positive element of $c_0(\ZM)$, then 

$$\left(h [\nabla_k-\nabla_{k'},b] h\;f \right)_{n,r}=
    h_n\,(e^{\imath kr}-e^{\imath k'r}) h_{n-r}\;
   \frac{\alpha^{-n}(b_n)- \alpha^{-n+r}(b_{n-r})}{d_r}\;
    f_{n-r,r}\,,
$$

\noindent showing that the map $k\in\TM\mapsto h\,v_k [D_\Bb,b] v_k^{-1}\,h\,\in \Bb(\Kk)$ is continuous in norm.

\vspace{.1cm}

\noindent Let now $c= \sum_{l\in\ZM}c_l\,u^l\,\in \Clip$ where $\Clip$ denotes the Lipschitz ball of $Y\rtimes_{\alpha_\ast}\ZM$. Then $c_l\in\Bb$ can be written as a sequence $(c_{n,l})_{n\in\ZM}$ such that $c_{n,l}\in\Aa$. Moreover, by assumption $\|[D_\Bb,c]\|\leq 1$. Then applying the dual action and integrating over $k\in\TM$ gives

\begin{eqnarray*}
\left(u^{-l} \int_{\TM} 
   v_k[D_\Bb,c]v_k^{-1} \; e^{-\imath kl}\,\frac{dk}{2\pi} \;f
  \right)_{n,r}&=&
   \frac{\gamma_1+\imath\gamma_2}{2} \;
    \frac{\alpha^{-n}(c_{n,l-r})-\alpha^{-n+r}(c_{n-r,l-r})}{d_r}\;
     f_{n-l,r}\\
&&+
   \frac{\gamma_1-\imath\gamma_2}{2} \;
    \frac{\alpha^{-n}(c_{n,l+r})-\alpha^{-n-r}(c_{n+r,l+r})}{d_r}\;
     f_{n-l,r}\\
&&+
  \gamma_3 l(c_lu^lf)_{n,r}
   +\gamma_5 ([D,c_l]u^l f)_{n,r}\,.
\end{eqnarray*}

\noindent Taking a partial trace over the various Dirac matrices, leads to the following necessary conditions for $c$ to be in the Lipschitz ball $\Clip$

$$\|\alpha^{-n}(c_{n,l})-\alpha^{-n+r}(c_{n-r,l})\|\leq d_r\,,
   \hspace{1cm}
    \|[D,\alpha^{-n}(c_{n,l})]\| \leq 1\,,
     \hspace{1cm}
      \|\partial c\| \leq 1\,.
$$

\noindent The first inequality implies $\|c_{n,l}\|\leq \sup_r d_r\leq 1$. Hence with the second inequality, it means that $c_{n,l}\in\Blip^1$ for all $n,l$. In particular, since $\Blip^1$ is norm compact in $\Aa$, if $h$ is a strictly positive element of $c_0(\ZM)$ it follows that $hch$ satisfies the conditions of Proposition~\ref{met10.prop-comp} in Section~\ref{met10.sect-Xprod}. In particular, $h\Clip h$ is norm compact.

\vspace{.5cm}

\section{Noncommutative Tori with Real Multiplication}
\label{met10.sect-NCtori}

\noindent We give here an example of a situation where the problem of
extending a spectral triple to a crossed product by the action of the
integers arises naturally in a context motivated by the arithmetic geometry
of real quadratic fields. The noncommutative space described by the
original spectral triple is a noncommutative torus with real multiplication, while the action of the integers is induced by the action of the group of units of the real quadratic field. The geometry that arises in this case, in fact, turns out to be very similar to the prototype example of the Arnold cat map that we described earlier in the paper.

\vspace{.1cm}

\noindent In the '70s Hirzebruch formulated a conjecture expressing the signature defects of Hilbert modular surfaces in terms of the Shimizu $L$-function of real quadratic fields \cite{Hirz}. This conjecture was proved by Atiyah--Donnelly--Singer \cite{ADS} using the Atiyah--Patodi--Singer index formula, together with a detailed analysis of the induced operator on the solvmanifold that gives the link of a cusp singularity of the Hilbert modular surface. The main part of the argument of \cite{ADS} consists of separating out in the eta function of this operator a part that recovers the Shimizu $L$-function, whose value at zero gives the signature defect, and a remaining part that vanishes at zero and is not of arithmetic nature. It was recently shown in \cite{Mar-solv} that the Shimizu L-function of a real quadratic field is obtained from a spectral triple on a noncommutative torus with real multiplication, as an adiabatic limit of the Dirac operator on a 3-dimensional solvmanifold. The Dirac operator on this 3-dimensional geometry gives, via the Connes-Landi isospectral deformations \cite{CoLa}, a spectral triple for the noncommutative tori obtained by deforming the fiber tori to noncommutative spaces. It is also shown in \cite{Mar-solv} that the 3-dimensional solvmanifold of \cite{ADS} is the homotopy quotient in the sense of Baum--Connes of the noncommutative space obtained as the crossed product of the noncommutative torus by the action of the units of the real quadratic field. This noncommutative space is identified with the twisted group $C^\ast$-algebra of the fundamental group of the 3-manifold. The twisting can be interpreted as the cocycle arising from a magnetic field, as in the theory of the quantum Hall effect. The resulting interpretation of the Shimizu $L$-function in terms of the noncommutative geometry of noncommutative tori with real multiplication, in the sense of \cite{Man}, is interesting in view of the conjectured role of noncommutative tori in the Stark conjecture for real quadratic fields proposed in \cite{Man}. In fact, as explained also in \cite{Mar-ICM}, the result of \cite{Mar-solv} deals with a special case of the zeta functions defining the Stark numbers, which are the conjectural generators of abelian extensions of real quadratic fields, and it shows how to relate these to the geometry of noncommutative tori.

\vspace{.3cm}

 \subsection{Lattices and Noncommutative Tori}
 \label{met10.ssect-latNCtor}

\noindent To fix notation, we let $L\subset \KM$ be a lattice in a real quadratic field $\KM=\QM(\sqrt{d})$, with $U_L^+$ the group of totally positive units preserving $L$,

\begin{equation}
\label{met10.eq-ULambdaplus}
U_L^+ =
 \{ u\in O_\KM^* \,|\, u L \subset L, \ \alpha_i(u)\in \RM^*_+ \},
\end{equation}

\noindent with $\alpha_i$, $i=1,2$, the two embeddings of $\KM$ in $\RM$.
We let $V$ denote a finite index subgroup of $U_L^+$. We denote by $\epsilon$ a generator, so that $V=\epsilon^\ZM$.

\vspace{.1cm}

\noindent Let $\Lambda$ be the lattice in $\RM^2$ obtained by embedding $L$ in $\RM^2$ by

\begin{equation}
\label{met10.eq-R2lattice}
\Lambda=\{ (\alpha_1(\ell),\alpha_2(\ell))\,|\, \ell\in L \}.
\end{equation}

\noindent The group $V$ acts on $\Lambda$ by $\lambda=(\alpha_1(\ell),\alpha_2(\ell))\mapsto (\epsilon \alpha_1(\ell), \epsilon^{-1} \alpha_2(\ell))$. The following solvable Lie group will be considered

\begin{equation}
\label{met10.eq-solvLie}
S(\RM^2,\RM,\epsilon)=\RM^2\rtimes_\epsilon \RM,
\end{equation}

\noindent where $\RM$ acts on $\RM^2$ as the one-parameter subgroup $\{ e^{t\log A_\epsilon} \}_{ t\in \RM}$ of $\SL_2(\RM)$ with

\begin{equation}
\label{met10.eq-Aepsilon}
A_\epsilon = 
  \left[
   \begin{array}{cc}
     \epsilon & 0 \\
     0 & \epsilon'
   \end{array}
  \right]\; \in\; \SL_2(\RM).
\end{equation}

\noindent As in \cite{ADS}, the following discrete subgroup will be needed

\begin{equation}
\label{met10.eq-SLambdaV}
S(\Lambda,V)=\Lambda\rtimes_{A_\epsilon} \ZM 
\end{equation}

\noindent together with the quotient

\begin{equation}
\label{met10.eq-Xmanifold}
X_\epsilon =S(\Lambda,V)\backslash S(\RM^2,\RM,\epsilon).
\end{equation}

\noindent The noncommutative torus of irrational modulus $\theta$ is topologically the irrational rotation $C^\ast$-algebra $\Aa_\theta = C(S^1)\rtimes_\theta \ZM$. It can be equivalently described as a twisted group $C^\ast$-algebra $C^\ast_r(\ZM^2,\sigma)$, where $\sigma$ is a 2-cocycle on $\Gamma=\ZM^2$, namely a multiplier $\sigma: \Gamma\times \Gamma \to U(1)$ 
satisfying 

\begin{equation}
\label{met10.eq-sigmacocycle}
\sigma(\gamma_1,\gamma_2)\sigma(\gamma_1\gamma_2,\gamma_3)=
\sigma(\gamma_1,\gamma_2\gamma_3)\sigma(\gamma_2,\gamma_3),
\end{equation}

\noindent with $\sigma(\gamma,1)=\sigma(1,\gamma)=1$.  For the noncommutative torus the cocycle $\sigma$ can be taken of the form

\begin{equation}
\label{met10.eq-sigmaZ2}
\sigma((n,m),(n',m')):= 
 \exp(-2\pi i(\xi_1 nm' + \xi_2 mn')), 
\hspace{2cm}
 \theta =\xi_2-\xi_1.
\end{equation}

\noindent In addition, upon choosing $\xi_2=-\xi_1 =\theta/2$, the following property holds

\begin{equation}
\label{met10.eq-SL2inv}
\sigma((n,m),(n',m'))=
 \sigma((n,m)\varphi,(n',m')\varphi), 
\hspace{2cm} 
\forall \varphi \in \SL_2(\ZM) .
\end{equation}

\noindent We consider the noncommutative torus $C^\ast_r(\Lambda,\sigma)$, where the lattice $\Lambda$ is as in eq.~(\ref{met10.eq-R2lattice}) of the form $\Lambda=\{ (n+m\theta,n+m\theta')\,|\, (n,m)\in \ZM^2\}$, with $\theta$ and $\theta'$ Galois conjugates. The algebra is the norm closure of the twisted group ring $\CM(\Lambda,\sigma)$ generated by

\begin{equation}
\label{met10.eq-Rsigmalambda}
R^\sigma_\eta e_{\lambda} = 
 \sigma(\eta,\lambda) e_{\lambda+\eta},
\end{equation}

\noindent where $e_\lambda$ is the canonical basis of $\ell^2(\Lambda)$ and
$\sigma(\eta,\lambda)=\sigma((n,m),(u,v))$, for $\eta=(n+m\theta,n+m\theta')$ and $\lambda=(u+v\theta,u+v\theta')$. These satisfy the relation $R^\sigma_\eta R^\sigma_\lambda =\sigma(\eta,\lambda) R^\sigma_{\eta+\lambda}$.

\vspace{.1cm}

\noindent The cocycle $\sigma$ as in eq.~\ref{met10.eq-sigmaZ2} with eq.~\ref{met10.eq-SL2inv}, which we can write as

\begin{equation}
\label{met10.eq-sigmaLambda}
\sigma(\eta,\lambda)=
 \exp(-\pi i\theta \eta\wedge\lambda),
\end{equation}

\noindent extends to a cocycle $\tilde\sigma$ on the group $S(\Lambda,V)$, of the form

\begin{equation}
\label{met10.eq-tildesigma}
\tilde\sigma((\eta,\epsilon^k),(\lambda,\epsilon^r))=
 \sigma(\eta,A_\epsilon^k(\lambda))=
  \sigma((n,m),(u,v)\varphi_\epsilon^k),
\end{equation}

\noindent where $\varphi_\epsilon\in \SL_2(\ZM)$ is given by

\begin{equation}
\label{met10.eq-varphi}
\varphi_\epsilon=
 \left[\begin{matrix}
      a & b\\ 
      c & d 
 \end{matrix}\right]
\hspace{1cm}\text{ with } \hspace{1cm}
\epsilon = a + b\theta
\hspace{1cm}\text{ and }\hspace{1cm}
 \epsilon\theta = c + d \theta.
\end{equation}

\noindent so that $\theta$ and $\theta'$ are the two fixed points of the transpose of $\varphi_\epsilon$.

\vspace{.1cm}

\noindent The group $S(\Lambda,V)$ is amenable, so that the maximal and reduced (twisted) group $C^\ast$-algebras coincide. The group $S(\Lambda,V)$ satisfies the Baum--Connes conjecture (with coefficients),  since it is isomorphic to the crossed product $\ZM^2\rtimes_{\varphi_\epsilon}\ZM$, and the 3-manifold $X_\epsilon$ is the homotopy quotient, in the sense of  Baum--Connes, of the noncommutative space $C^\ast(S(\Lambda_\epsilon,V),\tilde\sigma)$.

\vspace{.1cm}

\noindent In particular, the algebra $C^\ast(S(\Lambda_\epsilon,V),\tilde\sigma) =C^\ast(\Lambda,\sigma)\rtimes V$ is the noncommutative space describing the quotient of the noncommutative torus $\TM_\theta$ with real multiplication, described by the algebra $C^\ast(\Lambda,\sigma)$ by the action of the group of units $V$. Thus, it provides a geometric analog, in the case of real quadratic fields, of the quotient $E/\Aut(E)=\PM^1$ of an elliptic curve with complex multiplication in the case of imaginary quadratic fields.

\vspace{.1cm}

\noindent In \cite{Mar-solv} one relates the Shimizu $L$-function to the geometry of the noncommutative torus $\TM_\theta$ by considering the
Dirac operator on the 3-dimensional solvmanifold $X_\epsilon$ as in \cite{ADS} and its isospectral deformations to a fibration of noncommutative tori, as in \cite{CoLa}. From the point of view we  are considering here, one can proceed in a different way, working with spectral triples on the noncommutative torus $\TM_\theta$ with algebra $C^\ast(\Lambda,\sigma)$, and the induced spectral triples, via the metric bundle construction, on the crossed product by $\ZM$ that gives the noncommutative space $C^\ast(S(\Lambda_\epsilon,V),\tilde\sigma)$.

\vspace{.3cm}

 \subsection{A Spectral Metric Bundle}
 \label{met10.ssect-stNCtori}

\noindent As in \S 7.4 of \cite{Mar-solv}, we consider on the noncommutative torus $\TM_\theta$ a spectral triple with representation

\begin{equation}
\label{met10.eq-Sp3Ttheta}
R^\sigma_{(r,k)} \psi_{n,m} =
 \sigma((r,k),(n,m)) \psi_{(n,m)+(r,k)}
\end{equation}

\noindent and Dirac operator

\begin{equation}
\label{met10.eq-Sp3TthetaD}
D_{\theta,\theta'} =
\left[
 \begin{array}{cc}
 0 & \delta_{\theta'}-i \delta_\theta \\
\delta_{\theta'} + i \delta_\theta & 0 
 \end{array}
\right],
\end{equation}

\noindent where $\theta'$ is the Galois conjugate of the quadratic irrationality $\theta$ and the derivations $\delta_\theta$ and $\delta_{\theta'}$ are given by

\begin{equation}
\label{met10.eq-deltatheta}
\delta_\theta: 
  \psi_{n,m} \mapsto (n+m\theta) \psi_{n,m}, 
\hspace{2cm}
\delta_{\theta'}: 
  \psi_{n,m} \mapsto (n+m\theta') \psi_{n,m}.
\end{equation}

\noindent This is equivalent to $D_{\theta,\theta'}: \psi_\lambda \mapsto (\lambda_1 \sigma_1+\lambda_2\sigma_2)\psi_\lambda$, with $\sigma_k$ the Pauli matrices (see eq.~\ref{met10.eq-Pauli}) and $\lambda=(\lambda_1,\lambda_2)$, which is the Dirac operator induced on the fiber noncommutative tori by the isospectral deformation of the 3-manifold $X_\epsilon$. The commutators of the Dirac operator with the generators of the algebra are bounded

$$[ D_{\theta,\theta'}, R^\sigma_{(r,k)}] =
  \left[
   \begin{array}{cc} 
    0 & (r+k\theta')-i(r+k\theta) \\ 
    (r+k\theta') + i(r+k\theta) & 0
   \end{array} 
  \right]. 
$$

\noindent Since $V$ acts in an hyperbolic way, similarly to the Arnold cat map acting on the commutative torus, it is required to construct the noncommutative analog of the metric bundle according to the general description given in Section~\ref{met10.ssect-crossp} in order to describe the metric space associated with the crossed product. This leads to consider the sequences $(b_n)_{n\in \ZM}$ of elements $b_n \in C^\ast(\Lambda,\sigma)$,

$$ b_n =
   \sum_{\lambda} c_{n,\lambda} R^\sigma_\lambda 
$$

\noindent with $\lim_n \| b_n \| =0$. They define elements $b\in \Bb_{\Lambda,\sigma}:= C^\ast(\Lambda,\sigma) \otimes c_0(\ZM)$. The action $\alpha_\epsilon$ of the group $V$ of units on $C^\ast(\Lambda,\sigma)$ by 

$$ \alpha_\epsilon^k :
    R^\sigma_\lambda \mapsto R^\sigma_{A_\epsilon^k(\lambda)} 
$$

\noindent or equivalently $R^\sigma_{(n,m)} \mapsto R^\sigma_{(n,m)\varphi_\epsilon^k}$, extends to an action on $\Bb_{\Lambda,\sigma}$ by setting

\begin{equation}
\label{met10.eq-phistaract}
((\alpha_\epsilon)_* (b))_n= 
  \alpha_\epsilon(b_{n-1})=
   \sum_\lambda c_{n-1,\lambda}
    R^\sigma_{A_\epsilon(\lambda)} =
    \sum_\lambda c_{n-1,A_\epsilon^{-1}(\lambda)}
      R^\sigma_\lambda. 
\end{equation}

\noindent In this case it is then convenient to take as Hilbert space $\Hh \otimes \ell^2(\ZM)$. This has a representation of the crossed product algebra  $\Bb_{\Lambda,\sigma}\rtimes_{(\alpha_\epsilon)_*}V$ by

$$ b \psi_{\lambda,n} = 
   \alpha_\epsilon^{-n}(b_n) \psi_{\lambda,n}
  = \sum_{\lambda'} c_{n,\lambda'} 
     R^\sigma_{A^{-n}_\epsilon(\lambda')} 
      \psi_{\lambda,n} ,
$$

\noindent where we write the latter as

$$ \sum_{\lambda'} 
   c_{n,A_\epsilon^n(\lambda')} 
    R^\sigma_{\lambda'} \psi_{\lambda,n} =
     \sum_{\lambda'} c_{n,A_\epsilon^n(\lambda')} 
      \psi_{\lambda+\lambda',n}, 
$$

\noindent and the unitary $\upsilon_\epsilon$ implementing the action of the generator of $V$ is given by

$$ \upsilon_\epsilon \, \psi_{\lambda,n} = 
    \psi_{\lambda, n-1}, 
$$

\noindent so that $\upsilon_\epsilon\, b\, \upsilon_\epsilon^{-1} =(\alpha_\epsilon)_*(b)$. We consider on this Hilbert space the Dirac operator

\begin{equation}
\label{met10.eq-DiracLambdaV}
\hat D \psi_{\lambda,n} = 
  \sigma_3 n \log(\epsilon) \psi_{\lambda,n} + 
   D_{\theta,\theta'} \psi_{\lambda,n},
\end{equation}

\noindent where $\sigma_3$ is the third Pauli matrix (see eq.~\ref{met10.eq-Pauli}), and $\log(\epsilon)$ is the length of the unit $\epsilon$, which is the length of the base circle in the fibration $T^2 \to X_\epsilon \to S^1$.

\vspace{.1cm}

\noindent As in \S 4 of \cite{Mar-solv}, we write the generators of the crossed product algebra $C^\ast(S(\Lambda,V),\tilde\sigma)$ as $R^{\tilde\sigma}_{(n,m,k)}$ or $R^{\tilde\sigma}_{(\lambda,k)}$ and we identify them with the elements $R^\sigma_\lambda \upsilon_\epsilon^k$ since we have

$$ R^{\tilde\sigma}_{(n,m,k)} 
    R^{\tilde\sigma}_{(n',m',k')} =
     \tilde\sigma((n,m,k),(n',m',k')) \, 
      R^{\tilde\sigma}_{(n,m,k)(n',m',k')} 
$$
$$ = \sigma((n,m),(n',m')\varphi_\epsilon^k) 
      R^\sigma_{(n,m)+(n',m')\varphi_\epsilon^k} \, 
      \upsilon_\epsilon^{k+k'} = 
       R^\sigma_{(n,m)}\upsilon_\epsilon^k 
       R^\sigma_{(n',m')} \upsilon_\epsilon^{k'}. 
$$

\noindent One has in this way also a representation of $C^\ast(S(\Lambda,V),\tilde\sigma)$ on the Hilbert space $\Hh \otimes \ell^2(\ZM)$ by identifying the $R^\sigma_\lambda$ with the constant sequence $b_n =R^\sigma_\lambda$ for all $n\in \ZM$, that is, $c_{n,\lambda'}=\delta_{\lambda',\lambda}$ for all $n\in \ZM$. 

\vspace{.1cm}

\noindent Then the commutator

\begin{equation}
\label{met10.eq-hatDu}
[\hat D, \upsilon_\epsilon]= 
  -\sigma_3 \log(\epsilon) \upsilon_\epsilon ,
\end{equation}

\noindent is clearly bounded, and

\begin{equation}
\label{met10.eq-hatDb}
[\hat D, b] \psi_{\lambda,n} = 
  [D_{\theta,\theta'}, \alpha_\epsilon^{-n}(b_n)] \psi_{\lambda,n},
\end{equation}

\noindent is bounded for all sequences $b=(b_n)$ with finite support. Thus, one obtains a spectral triple for the algebra $\Bb_{\Lambda,\sigma}\rtimes_{\alpha_\epsilon} V$ with bounded commutators on the dense subalgebra $\Bb_{c,\Lambda,\sigma}\rtimes_{\alpha_\epsilon} V$, with $\Bb_{c,\Lambda,\sigma}$ consisting of the sequences with finite support. The algebra $C^\ast(S(\Lambda,V),\tilde\sigma)$ acts as multipliers with commutators with
$\hat D$ that are unbounded multipliers of $\Bb_{c,\Lambda,\sigma}\rtimes_{\alpha_\epsilon} V$.

\vspace{.1cm}

\noindent The Dirac operator $\hat D$ of eq.~(\ref{met10.eq-DiracLambdaV}) is in fact given by the Fourier modes of the Dirac operator on the 3-manifold $X_\epsilon$, namely

\begin{equation}
\label{met10.eq-DiracXepsilon}
D_{X_\epsilon} = 
 \left[ 
  \begin{array}{cc} 
   \frac{\partial}{\partial t} & 
  e^{-t}\frac{\partial}{\partial y} -i e^t \frac{\partial}{\partial x} \\
   e^{-t} \frac{\partial}{\partial y} + i e^t \frac{\partial}{\partial x} &
    -\frac{\partial}{\partial t} 
  \end{array}
 \right],
\end{equation}

\noindent when the spinors $\psi((x,y),t)$ on $X_\epsilon$ are expanded in the form $\psi = \sum_\lambda \psi_\lambda E_\lambda$, with $E_\lambda = \exp(2\pi i \langle\Theta_{-t}(x,y),\lambda\rangle)$, with $\langle (a,b),\lambda\rangle = a\lambda_1 + b \lambda_2$ and $\Theta_{-t}(x,y)=(e^{-t}x, e^t y)$, so that

$$ D_{X_\epsilon} E_\lambda = 
   (\frac{\partial}{\partial t} \sigma_3 + 
     2\pi i \lambda_1 \sigma_1 + 2\pi i \lambda_2 \sigma_2) 
      E_\lambda. 
$$

\noindent This example is in essence very similar to the case of the Arnold cat map, with the matrix $A_\epsilon \in \SL_2(\RM)$ acting on the fiber torus $T^2$ of the mapping torus $X_\epsilon$. 

\vspace{.1cm}

\noindent Up to a unitary transformation, as in \cite{ADS}, \cite{Mar-solv}, the Dirac operator $D_{\theta,\theta'}$ on the noncommutative torus ca be written as

\begin{equation}
\label{met10.eq-Dtheta0}
D_{\theta,\theta', 0} =
 \sum_{\mu \in (\Lambda \smallsetminus\{0\})/V} 
   D^\mu_{\theta,\theta', 0} ,
\end{equation}

\noindent where $D_{\theta,\theta', 0}$ is the restriction to the complement of the zero modes $\lambda=0$, and $D^\mu_{\theta,\theta', 0}$ is given by

\begin{equation}
\label{met10.eq-Dmutheta}
D^\mu_{\theta,\theta', 0} \psi_{A^k_\epsilon(\mu)}= 
 {\rm sign}(N(\mu))\, |N(\mu)|^{1/2} \, 
  (\epsilon^k \sigma_1 + \epsilon^{-k} \sigma_2) \, 
  \psi_{A^k_\epsilon(\mu)},
\end{equation}

\noindent where a fundamental domain $\Ff_V$ for the action of $V$ on $\Lambda \smallsetminus\{0\}$ is chosen and $\lambda\in \Lambda$, if $\lambda\neq 0$, is written as $\lambda =A^k_\epsilon(\mu)$ for some $\mu \in \Ff_V$ and some $k\in \ZM$. This formulation shows clearly that the action of $\ZM$ that implements the crossed product changes the metric on the noncommutative torus by mapping ${\rm sign}(N(\mu))\, |N(\mu)|^{1/2} \, 
(\epsilon^k \sigma_1 + \epsilon^{-k} \sigma_2) \mapsto  {\rm sign}(N(\mu))\, |N(\mu)|^{1/2} \, (\epsilon^{k+1} \sigma_1 + \epsilon^{-(k+1)} \sigma_2)$.

\vspace{.2cm}

\noindent This example, just as in the case of the Arnold cat map, gives a clear illustration of the subtleties involved in the problem of the choice of the metric on $\ZM$ that we discussed in the previous sections. In fact, the natural Dirac operator to consider here, which in this case is dictated by the requirement that it recovers the 3-manifold geometry of \cite{ADS}, corresponds to the usual metric on $\ZM$, which, however, as we have seen, requires the use of the smaller state space characterized by the finite first moment condition. While one can make sense of this condition in the specific case at hand, a general treatment of such first moment conditions in the general setting of noncommutative geometry is missing and will be an interesting topic of investigation.

\vspace{.5cm}

\section{Cuntz--Krieger Algebras and AF-Algebras}
\label{met10.sect-CKAF}

\noindent One of the earlier results on the interaction between noncommutative and arithmetic geometry was the use of spectral triples techniques to model the ``fiber at infinity" in the Arakelov geometry of arithmetic surfaces \cite{ConsMa2}, in terms of the choice of a Schottky uniformization as in \cite{Man3d}, and of the action of the Schottky group on its limit set. This point of view gave rise to the construction of a $\theta$-summable, not finitely summable spectral triple, where the algebra representation was related to the archimedean cohomology at infinity, see also \cite{ConsMa3}. Based on the analogies between the degenerate fiber at the archimedean place in the Arakelov geometry of arithmetic surfaces, and the case of Mumford curves with $p$-adic uniformization, the same technique were adapted in \cite{ConsMa} to describe a noncommutative space associated to the action of a p-adic Schottky group on its limit set. The construction
was later refined in \cite{CMRV} and \cite{CaMaRe}.

\vspace{.1cm}

\noindent The original motivation for considering spectral triples on crossed products by the integers in \cite{CMRV} was to improve the $\theta$-summable spectral triple on the crossed product algebra $C(\Lambda_\Gamma)\rtimes \Gamma$ of a Kleinian or $p$-adic Schottky group on its limit set, considered in \cite{ConsMa}. The idea was to replace the algebra $C(\Lambda_\Gamma)\rtimes \Gamma$, which can be described in terms of a Cuntz--Krieger algebra $O_A$, with its stabilization $\overline{O}_A =\overline{\Ff}_A \rtimes_T \ZM$.  One then considers the problem of
constructing a spectral triple on the non-unital AF algebra $\Aa=\overline{\Ff}_A$ and extending it to the crossed product $\overline{\Ff}_A \rtimes_T \ZM$. 

\vspace{.1cm}

\noindent The case of $\overline{\Ff}_A$ can be handled in the way proposed by Christensen and Ivan \cite{ChI06}, namely by choosing eigenvalues with $|\lambda_n| \geq (\dim \Aa_n)^q$, for $\Aa_n \subset \Aa_{n+1}$ the filtration of the AF algebra and with $q> 2/p$. This determines a $p$-summable Dirac operator on a Hilbert space $\Hh$, which is the GNS space of a state $\varphi$, with quotient map $\pi_\phi: \Aa \to \Hh$ and $\Hh_n=\pi_\phi(\Aa_n)$, with $\dim\Hh_n \leq \dim \Aa_n$. 

\vspace{.1cm}

\noindent In general, if the restrictive quasi-invariance assumption of \cite{CMRV} is removed, the metric bundle construction leads to a spectral triple on the noncommutative metric bundle algebra $\Bb \rtimes_{T_*}\ZM$,
with $\Bb = \Aa \otimes c_0(\ZM)$, on which $\Aa\rtimes_T\ZM$ acts as multipliers, instead of a spectral triple for the crossed product $\overline{\Ff}_A \rtimes_T \ZM$.

\vspace{.1cm}

\noindent In this case, the metric on $\ZM$ and the form of the Dirac operator on the crossed product by the integers is not strictly dictated by the geometry, so that one can modify the construction by choosing a metric on $\ZM$ that satisfies the conditions of \S \ref{met10.ssect-zm} and obtain a spectral metric space with the desired properties.

\vspace{.5cm}

\section{Algebro-Geometric Codes and Spectral Triples}
\label{met10.sect-codes}

\noindent In the recent work \cite{ManMar}, the asymptotic bound problem in the theory of algebro-geometric codes was reformulated in terms of quantum statistical mechanical systems and operator algebras associated to coding maps. In particular, it is shown in \cite{ManMar} that the code parameters, rate and relative minimum distance, can be recovered from the Hausdorff dimensions of certain Cantor-like fractals associated with the code.

\vspace{.1cm}

\noindent Just as the prototype example of the Arnold cat map provides a model for the application to the case of noncommutative tori with real multiplication described above in Section~\ref{met10.sect-NCtori}, the other prototype example, of the bilateral shift and Cantor sets described in the Example~\ref{met10.exam-shift} (Section~\ref{met10.ssect-neqact}) serves as a model for the application to algebro-geometric codes.

\vspace{.3cm}

 \subsection{The Spectral Metric Space of Codes}
 \label{met10.ssect-smsCodes}

\noindent The basic terminology of \cite{ManMar} is given here. Let $\AG$ be a finite alphabet of cardinality $\# \AG=q$. The main application here we will focus on the case where $\AG$ is a finite field $\FM_q$, which is set theoretically identified with the set $\{ 0,\ldots, q-1 \}$ of $q$-ary digits. For some given $n\in \NM$ let $\AG^n$ be the $n$-hypercube. The Hamming metric on $\AG^n$ is defined as $d(x,y) =\# \{ k\,|\, x_k \neq y_k \}$, for $x=(x_1,\ldots,x_n)$ and $y=(y_1,\ldots,y_n)$ in $\AG^n$. This
satisfies the bound $d(x,y)\leq n$.

\begin{defini}
\label{met10.def-ncodes}
An $[n,k,d]_q$-code is a subset $C\subset \AG^n$ with

$$ k=\log_q (\# C), 
\hspace{1cm} \text{ and } \hspace{1cm}
d=\min\{ d(x,y)\,|\, x,y \in C, \, x\neq y \}. 
$$

\noindent The rate of the code is $R=k/n$ and the relative minimum
distance is $\delta=d/n$.
\end{defini}

\noindent In particular, if the alphabet has $q$ letters with $q=p^r$, it will be identified with the elements of a finite field $\FM_q$. However, the codes $C: \FM_q^k \hookrightarrow \FM_q^n$ will not be required to be $\FM_q$-linear maps, that is, nonlinear codes are included in the present description.

\vspace{.1cm}

\noindent Let $Q^n =[0,1]^n$ be the standard unit cube in $\RM^n$. A point in $Q^n$ will be written as $x=(x_1,\ldots,x_n)$ using the q-adic expansion of points in $[0,1]$. Then, a fractal $\bar S_C$ in $Q^n$ will be associated to a code $C$, by first subdividing the unit cube $Q^n$ into $q^n$ smaller cubes, each of volume $q^{-n}$: each of these smaller cubes consist of the points $x\in Q^n$ such that in the $q$-adic expansion
the first digits $(x_{11},\ldots,x_{n1})$ are equal to a given element 
$(a_1,\ldots, a_n)$ in $\AG^n$. Among these one then only those for which $(x_{11},\ldots,x_{n1})$ is an element of $C\subset \AG^n$ are kept and the others are deleted. The same procedure is then iterated on each of the remaining cubes. The set obtained in this way is the Sierpinski fractal $\bar S_C$ consisting of all points $x\in Q^n$ such that all the points $(x_{1k},\ldots,x_{nk}) \in \AG^n$, for all $k\geq 1$, belong to the code $C\subset \AG^n$.  The Hausdorff dimension $\dim_H(\bar S_C)$ of this fractal
is $\log_q(\# C)$; when normalized to the ambient space dimension it equals the rate of the code,

$$ \frac{\log(\# C)}{n \log q} = 
    \frac{k}{n} = R. 
$$

\noindent Similarly, for a given $[n,k,d]_q$-code $C\subset \AG^n$, the set 
$\Xi_C$ of doubly infinite sequences $x=(x_m)_{m\in\ZM}$ of code words $x_m \in C$ will be considered. Let also $\Ww$ denote the set of all words of finite length in the code language, that is, all finite sequences of elements of $C$. As recalled in \S 5 of \cite{ManMar} (see \cite{Eilen}) the {\em entropy} of a language is defined through the following 
procedure. The {\em structure function} of the code language is defined by counting words of a given length

$$ s_C(N)=\# \{ w\in \Ww\,:\, |w|=N \}. 
$$

\noindent Here the length of a word $w\in \Ww$ is defined as the number of letters in the underlying alphabet $\AG$, rather than the number of code words $w$ consists of. Namely if $w= w_1 \cdots w_m$ with $w_i \in C$ this gives $|w|=n m$, where $C\subset \AG^n$. The generating function for the $s_C(N)$ is given by

\begin{equation}
\label{met10.eq-GCseries}
 G_C(t) = \sum_N s_C(N) t^N ,
\end{equation}

\noindent and the {\em entropy} of the code language is defined as

\begin{equation}
\label{met10.eq-Centropy}
 \Ee_C = - \log_{\# \AG} \rho(G_C) ,
\end{equation}

\noindent which is the logarithm of the radius of convergence $\rho(G_C)$.
It turns out that the series in eq.~(\ref{met10.eq-GCseries}) and the entropy in eq.~(\ref{met10.eq-Centropy}) have a natural interpretation in terms of the spectral triples of \cite{PB09}. In fact, the spectral metric spaces $X_{\tau,C}=(\Cc(\Xi_C), \Hh_C, \pi_\tau, D)$ is defined as in Example~\ref{met10.exam-shift}, where  $\Hh_C = \ell^2(\Ww_C)\otimes \CM^2$, with representations

$$ \pi_\tau (f)\psi(w) = 
    \left[ 
     \begin{array}{cc} 
      f(x_w) & 0 \\
       0 & f(y_w) 
     \end{array}
    \right] \psi(w), 
$$

\noindent for $\tau: \Ww_C \to \Xi_C\times \Xi_C$ a choice map, and with Dirac operator

$$ D \psi(w) = q^{|w|} \, 
    \left[ 
     \begin{array}{cc} 
     0&1 \\
     1&0 
     \end{array}
    \right] \psi(w). 
$$

\noindent Then the zeta function of the spectral triple is given by (\cite{PB09} \S 6)

$$\zeta_D(s) = 
   \sum_{w\in \Ww_C} q^{-|w|s} = 
    \sum_m q^{km} q^{-s nm} =
     \sum_m q^{(R-s)n m} = 
      (1-q^{(R-s)})^{-1}, 
$$

\noindent with convergence for $\Re(s)>R$. This gives
\begin{equation}\label{zetaDseriesC}
 \zeta_D(s) = G_C(q^{-s}) ,
\end{equation} 

\begin{equation}
\label{met10.eq-zetaDseriesC}
 \zeta_D(s) = G_C(q^{-s}) ,
\end{equation}

\noindent which recovers the structure functions for the code language, and the degree of summability of the spectral triple is the code rate, which is also equal to the entropy of the code language. In fact, $s_C(N)=0$ if $N\neq nm$ and $s_C(nm)=q^{km}$ for an $[n,k,d]_q$-code.

\vspace{.1cm}

\noindent A similar type of fractal construction can be done for the code parameter $d$ instead of $k$, as explained in \cite{ManMar}. Namely, 
the property that $C$ has minimum distance $d$ means that any pair of distinct points $x\neq y$ in $C$ must have at least $d$ coordinates that do not coincide, since $d(x,y)=\#\{ i\,|\, x_i\neq y_i \}$. Thus, in particular, this means that no two points of the code lie on the same $\pi$, for any $\pi$ as above of dimension $\ell \leq d-1$, while there exists at least one $\pi$ in $\Pi_d$ which contains at least two points of $C$. In terms of the iterative construction of the fractal $\bar S_C$, this means the following. For a given $\pi\in \Pi_\ell$ with $\ell\leq d-1$, if the intersection $C\cap \pi$ is non-empty it must consist of a single point. Thus, another fractal can be constructed, corresponding to a choice of a linear space $\pi\in \Pi_\ell$ with $\ell\leq d-1$, where at the first step the single cube $Q^\ell =Q^n\cap \pi$ is replaced with a scaled cube of volume $q^{-\ell}$, successively iterating the same procedure. This produces a family of nested cubes of volumes $q^{-\ell N}$ with
intersection a single vertex point. The Hausdorff dimension is clearly zero. When $\ell =d$ there exists a choice of $\pi\in \Pi_d$ for which 
$C\cap \pi$ contains at least two points. Then the same inductive construction, which starts by replacing the cube $Q^d =Q^n\cap \pi$ with $\# (C\cap \pi)$ copies of the same cube scaled down to have volume $q^{-d}$ produces a genuinely fractal object of positive Hausdorff dimension. (This is the same argument as in Proposition 3.3.1 of \cite{ManMar}, except that it is formulated here in terms of $\bar S_C$, as in \S 4.1.3 of \cite{ManMar}, rather than in terms of $S_C$ as defined in \S 3.1 of \cite{ManMar}.) The spectral metric space construction of \cite{PB09} can also be applied to the resulting fractal $S_{C\cap \pi}$ for $\pi\in \Pi_d$.

\vspace{.3cm}

 \subsection{The Shift Action on a Code Space}
 \label{met10.ssect-shCode}

\noindent The action of the bilateral shift on $\Xi_C$ or on $\Xi_{C\cap\pi}$ then produces a non-equicontinuous action. Geometrically, in such cases the resulting crossed product construction corresponds to considering code languages up to the natural equivalence relation given by the 
action of the bilateral shift. This gives rise to the type of crossed product construction we have been investigating. Once again we see that the issue of the choice of metric structure over $\ZM$ arises and determines the type of topological and metric properties that one expects to find on the state space.

\vspace{.1cm}

\noindent In this case the resulting noncommutative metric geometries have interesting global symmetries that come from Galois actions.

\vspace{.2cm}

\noindent Let $\FM_{q^m}$ be a field extension of $\FM_q$. Then a $[n,k,d]_q$-code $C: \FM_q^k \hookrightarrow \FM_q^n$ determines a $[n,k,d']_{q^m}$-code $C_{(m)}$ with $d'\leq d$, by setting $C_{(m)}: \FM_{q^m}^k \to \FM_{q^m}^n$ 

$$ C_{(m)}(u) =
   (C(u_1), \cdots, C(u_m)), 
$$

\noindent where the vector spaces $\FM_{q^m}$ and $\FM_q^m$ are identified. Hence $u\in \FM_{q^m}^k$ can be written as $u=(u_1,\ldots,u_m)$ with $u_i\in \FM^k$. In particular $\# C_{(m)} = (\# C)^m = (q^m)^k$, so that the parameter $k$ remains the same. To see that $d'\leq d$ let $\pi_{d,q}$ denote the $\FM_q$-linear space such that $\# (C\cap \pi_{d,q})\geq 2$
and such that all spaces $\pi$ of smaller dimension have $\# (C\cap \pi)\leq 1$. Then the $\FM_{q^m}$-linear space $\pi_{d,q^m}=\pi_{d,q}\otimes_{\FM_q} \FM_{q^m}$ satisfies the same property $\#( C_{(m)} \cap \pi_{d,q^m} ) \geq 2$. Thus, a code $C$ over $\FM_q$ defines a family of codes $C_{(m)}$ over all the fields extensions, with inclusions $C_{(\ell)} \subset C_{(m)}$ for $\ell |m$ induced by the corresponding inclusions of fields $\FM_{q^m} \subset \FM_{q^\ell}$. All the codes $C_{(m)}$ have the same rate $R$ and non-increasing $\delta$'s.

\vspace{.1cm}

\noindent Consider then the Galois action of $Gal(\FM_{q^m}/\FM_q)$ acting on both the source and target spaces $\FM_{q^m}^k$ and $\FM_{q^m}^n$ of a code $C: \FM_{q^m}^k \hookrightarrow \FM_{q^m}^n$, obtained as above, so that the coding map is equivariant with respect to this action. The Galois action of $Gal(\FM_{q^m}/\FM_q)$ then induces a homeomorphism of the fractal sets $\bar S_C$, induced by the Frobenius $\phi$. This in turn determines an automorphism $\phi$ of the algebra $\Cc(\Xi_C)$ and a unitary transformation of the Hilbert space $\Hh_C$, so that the representations $\pi_\tau$ are equivariant. The metric structure given by the spectral triple is also preserved by this action, since the length of words in the code language is preserved. This gives a global symmetries of the spectral geometry.

\vspace{.5cm}

\section{Appendix: Proof of the Proposition~\ref{met10.prop-lipcomp}}
\label{met10.sect-Xprod}

\noindent Let $\Aa$ be a unital \CS. Let $\alpha\in\Aut(\Aa)$ be an automorphism of $\Aa$. Then $\as$ will denote the crossed product algebra $\Aa\rtimes_\alpha \ZM$. Then $\as$ is generated by the elements of $\Aa$ and by a unitary element $u$ such that $uau^{-1}=\alpha(a)$. Hence $\as$ is the norm closure of $\as_c$ made of elements of the form $b=\sum_{l=-L}^L b_l\,u^l$ for some integer $L$, where $b_l\in\Aa$ for all $l$. The dual action is the family of automorphisms of $\as$ defined by $\eta_k(a) = a$ for $a\in\Aa$ and by $\eta_k(u) = e^{\imath k}u$, where $k\in\TM$. Then

\begin{enumerate}
  \item If $b\in \as_c$ then $\eta_k(b) = \sum_{l\in\ZM} b_l\,e^{\imath kl}\,u^l$.

  \item The $\ast$ derivation $\partial$ generating the dual group action is defined by $\partial b= d\eta_k(b)/dk\upharpoonright_{k=0}$ giving $\partial b= \sum_{l\in\ZM}\imath l\, b_l\,\,u^l$.

  \item The conditional expectation $\EM:\as\mapsto \Aa$ is defined by $\EM(b)= \int_{\TM} \eta_k(b) \, dk/2\pi = b_0$.\\
 In particular $\EM(bb^\ast)= \sum_l b_l\,b_l^\ast$.

  \item The equation $b_l = \EM(b\,u^{-l})$, valid on $\as_c$ extends to $\as$.

\end{enumerate}

\noindent The first important result of this Section is

\begin{proposi}
\label{met10.prop-comp}
For each $l\in\ZM$, let $K_l\subset \Aa$ be a compact set for the norm topology. Let $B(K)$ be the set of $b\in\as=\Aa\rtimes_\alpha\ZM$ such that (i) $b_0=\EM(b)=0$, (ii) $b_l=\EM(bu^{-l})\in K_l$ for all $l\in\ZM$, (iii) $\partial b\in \as$ and $\|\partial b\|\leq 1$. Then $B(K)$ is compact.
\end{proposi}

\noindent The proof of this result requires several intermediate steps.

\vspace{.3cm}

 \subsection{Algebraic Bounds}
 \label{met10.ssect-algb} 

\noindent The first step is given by

\begin{lemma}
\label{met10.lem-sobolev}
If $b\in\as$ is such that $\partial b\in\as$, then 

$$\|b-\EM(b)\|^2 \leq 
   \frac{\pi^2}{3}
    \|\EM(\partial b \,\partial b^\ast)\|
$$ 
\end{lemma}

\noindent  {\bf Proof: } Let $\pi$ be any representation of $\as$ and let $\Hh$ be the corresponding Hilbert space. Then, for $f,g\in\Hh$ and $b\in \as$ the following holds

$$\langle f |\pi(b-\EM(b) )\,g\rangle =
   \int_{-\pi}^\pi \frac{dk}{2\pi}\;
    \langle f |\pi(b-\eta_k(b) )\,g\rangle= -
     \int_{-\pi}^\pi \frac{dk}{2\pi}\;
      \int_{0}^k dp\;
       \langle f |\pi(\eta_p(\partial b) )\,g\rangle\,.
$$

\noindent Exchanging the order of integration gives

$$\langle f |\pi(b-\EM(b) )\,g\rangle = -
   \int_{-\pi}^\pi \frac{dp}{2\pi}\;
    (\pi-|p|)\,\sign(p)\,
    \langle f |\pi(\eta_p(\partial b) )\,g\rangle\,.
$$

\noindent Applying a Cauchy-Schwartz inequality leads to 

$$|\langle f |\pi(b-\EM(b) )\,g\rangle|^2 \leq
   \int_{-\pi}^\pi \frac{dp}{2\pi}\;
    (\pi-|p|)^2\;
   \int_{-\pi}^\pi \frac{dk}{2\pi}\;
    |\langle f |\pi(\eta_k(\partial b) )\,g\rangle|^2
$$

\noindent The first term in the {\em r.h.s.} gives $\pi^2/3$. Using again the Cauchy-Schwartz inequality, the other terms can be bounded by 

$$\int_{-\pi}^\pi \frac{dk}{2\pi}\;
  |\langle f |\pi(\eta_k(\partial b) )\,g\rangle|^2 \leq
   \|g\|^2 \int_{-\pi}^\pi \frac{dk}{2\pi}\;
    \langle f |\pi(\eta_k(\partial b\,\partial b^\ast) )\,f\rangle\leq
     \|g\|^2 \|f\|^2 \|\EM(\partial b\,\partial b^\ast)\| \,.
$$

\noindent Since $\pi$ is an arbitrary representation, this last estimates leads directly to the result.
\hfill $\Box$

\begin{coro}
\label{met10.cor-BKbnd1}
If $b\in B(K)$ then $\|b\|\leq \pi/\sqrt{3}$. 
\end{coro}

\noindent {\bf Proof: } This follows immediately from Lemma~\ref{met10.lem-sobolev} and from the definition of $B(K)$.
\hfill $\Box$

\vspace{.3cm}

 \subsection{The Fejer Kernel and Approximation Estimates}
 \label{met10.ssect-fejer}

\noindent The Fejer kernel is the function $F_N$ defined on $\TM$ by

\begin{equation}
\label{met10.eq-fejer}
F_N(k) =
   \frac{1}{N} \frac{\sin^2\{(N)k/2\}}{\sin^2{(k/2)}}=
    \sum_{n=N+1}^{N-1}
     \left(1-\frac{|n|}{N} \right) \;e^{\imath nk}
\end{equation}

\noindent It has the following properties: (i) it defines a probability distribution on the torus

$$F_N(k) \geq 0\,,
\hspace{2cm}
\int_{-\pi}^\pi \frac{dk}{2\pi}\; F_N(k) =1\,.
$$

\noindent (ii) it is concentrated in the vicinity of the origin

$$\Prob\left\{
    |k|\geq \frac{\pi}{\sqrt{N}}
       \right\} =
   \int_{|k|\geq \pi/\sqrt{N}}\;
    \frac{dk}{2\pi}\; F_N(k) \leq
     \frac{1}{\sqrt{N}}
$$

\noindent The Fejer approximant of order $N$ of $b\in\as$ is defined by

$$b^{(N)} = 
   \int_{-\pi}^\pi \frac{dk}{2\pi}\;
    F_N(k) \,\eta_k(b)
$$

\noindent This definition shows that $b^{(N)}\in\as$ and that $\|b^{(N)}\|\leq \|b\|$. In particular, the map $b\in\as\mapsto b^{(N)}\in\as$ is continuous. Moreover if $b\in\as_c$ it follows from eq.~(\ref{met10.eq-fejer}) that

$$b^{(N)}= 
   \sum_{|l| <N}
    \left(1-\frac{|l|}{N} \right) \; b_l\,u^l
$$

\noindent By continuity this equation is still valid for $b\in\as$ showing that $b^{(N)}$ is always an element of $\as_c$. 

\begin{lemma}
\label{met10.lem-fejapprox}
For any $b\in\as$, its Fejer approximants satisfy $\lim_{N\rightarrow \infty} \|b-b^{(N)}\| =0$. Moreover, if $\partial b\in\as$ then

$$\|b-b^{(N)}\| \leq
   \frac{\pi}{\sqrt{N}} \,\|\partial b\| +
    \frac{2}{\sqrt{N}}\, \|b\|\,.
$$
\end{lemma}

\noindent  {\bf Proof: } By definition 

$$b-b^{(N)} =  
   \int_{-\pi}^\pi \frac{dk}{2\pi}\;
    F_N(k) \,\left(b-\eta_k(b)\right)\,.
$$

\noindent The integral on the {\em r.h.s} will be decomposed into $\int_{|k|\geq \pi/\sqrt{N}} + \int_{|k|\leq \pi/\sqrt{N}}= \mbox{\rm (I)}+\mbox{\rm (II)}$. The first term can be estimated by

$$\mbox{\rm (I)} \leq 2\|b\| 
   \Prob\left\{
    |k|\geq \frac{\pi}{\sqrt{N}}
        \right\} \leq 
     \frac{2}{\sqrt{N}}\, \|b\|\,.
$$

\noindent To estimate the other term, let $\epsilon >0$ be chosen. Then, there is $N_\epsilon\in\NM$ such that, for $N\geq N_\epsilon$, then (i) $2\|b\|/\sqrt{N}\leq \epsilon/2$ and (ii) $\|b-\eta_k(b)\|\leq \epsilon/2$ as soon as $|k|\leq \pi/\sqrt{N}$. This gives the obvious bound 

$$N\geq N_\epsilon
 \hspace{1cm}\Rightarrow\hspace{1cm}
   \mbox{\rm (II)} \leq \frac{\epsilon}{2}\,,
\hspace{1cm}\Rightarrow\hspace{1cm}
    \|b-b^{(N)}\| \leq \epsilon\,.
$$

\noindent If now $\partial b\in\as$ the estimate of (II) can be improved. For indeed, by definition

\begin{equation}
\label{met10.eq-taylor}
b-\eta_k(b) = -
   \int_0^k dp \,\eta_p(\partial b)\,,
\hspace{1cm}\Rightarrow\hspace{1cm}
\| b-\eta_k(b)\| \leq |k|\|\partial b\|\,.
\end{equation}

\noindent Plugging this last inequality into (II), gives the result.
\hfill $\Box$

\begin{coro}
\label{met10.cor-BKbnd2}
If $b\in B(K)$ then $\|b-b^{(N)}\|\leq 2.2\pi/\sqrt{N}$. 
\end{coro}

\noindent  {\bf Proof: } Follows immediately from Lemma~\ref{met10.lem-fejapprox}, the definition of $B(K)$ and Corollary~\ref{met10.cor-BKbnd1}.
\hfill $\Box$

\vspace{.3cm}

 \subsection{Compactness Properties}
 \label{met10.ssect-comprop}

\begin{lemma}
\label{met10.lem-BKclosed}
Given a sequence $(K_l)_{l\in\ZM}$ of norm compact subsets of $\Aa$, the set $B(K)$ defined in Proposition~\ref{met10.prop-comp} is closed.
\end{lemma}

\noindent  {\bf Proof: } Let $b\in\overline{B(K)}$. Since $\as$ is a metric space, there is a sequence $(b_n)_{n\in\NM}$ contained in $B(K)$ converging in norm to $b$.

\vspace{.1cm}

\noindent (i) By definition $\EM(b_n)=0$ for all $n\in\NM$ so that, taking the limit $n\rightarrow\infty$ leads to $\EM(b)=0$. In addition $\EM(b_n\,u^{-l})\in K_l$. Since $K_l$ is compact, it is closed in norm so that, taking the limit $n\rightarrow\infty$ leads to $\EM(b\,u^{-l})\in K_l$ as well.

\vspace{.1cm}

\noindent (ii) Thanks to eq.~(\ref{met10.eq-taylor}), it follows that, for all $n\in\NM$

$$\left\| \frac{\eta_k(b_n)-b_n}{k}\right\| \leq 
   \|\partial b_n\|\leq 1
$$

\noindent Taking the limit as $n\rightarrow \infty$ implies 

$$\sup_{-\pi\leq k\leq \pi}
   \left\| \frac{\eta_k(b)-b}{k}\right\|
    \leq 1
$$

\noindent In particular, since $\eta_{k+k'}=\eta_k\circ\eta_{k'}$ (group property),  the map $k\in\TM\mapsto \eta_k(b)$ is Lipschitz continuous with Lipschitz constant $1$. Thanks to Lebesgue's Lemma, it follows that this map is almost surely differentiable. Let $k_0\in \TM$ be such that the differential exists at $k_0$ namely

$$\frac{d\eta_k(b)}{dk}\upharpoonright_{k=k_0}= c = 
   \lim_{h\rightarrow 0}
    \frac{\eta_{k_0+h}(b)-\eta_{k_0}(b)}{h}\,.
$$

\noindent Applying $\eta_{k_0}^{-1}=\eta_{-k_0}$ on both sides implies that  the map $k\in\TM\mapsto \eta_k(b)$ is differentiable at zero (actually everywhere) and that the derivative has a norm less than $1$. Therefore $b\in B(K)$
\hfill $\Box$

\vspace{.2cm}

\noindent {\bf Proof of Prop~\ref{met10.prop-comp}: } Let $(b_n)_{n\in\NM}$ be a sequence contained in $B(K)$. It will be proved that there is a convergent subsequence with limit in $B(K)$. Since $B(K)$ is metric, this is sufficient to prove the compactness. 

\vspace{.1cm}

\noindent (i) By construction the sequence $b_{n,l}= \EM(b_n\,u^{-l})$ is entirely contained in $K_l$. In particular $\hb_n= (b_{n,l})_{l\in\ZM\setminus\{0\}}\in \Omega$ if $\Omega$ denotes the product space $\prod_{l\neq 0}K_l$. Thanks to the Tychonov theorem, $\Omega$ is compact for the product topology. Therefore there is a subsequence $(\hb_{n_i})_{i\in\NM}$ which converges in $\Omega$. Replacing $(b_n)_{n\in\NM}$ by $(b_{n_i})_{i\in\NM}$, there is no loss of generality to assume that this is the sequence $(\hb_n)_{n\in\NM}$ itself. Hence for each $l\in\ZM$, the limit $\lim_{n\rightarrow \infty} b_{n,l}=b_l\in K_l$ exists.

\vspace{.1cm}

\noindent (ii) Let now $b^{(N)}$ be defined as

$$b^{(N)}= 
   \sum_{|l| <N}
    \left(1-\frac{|l|}{N} \right) \; b_l\,u^l =
     \lim_{n\rightarrow \infty} b_n^{(N)}\,,
$$

\noindent where $b_n^{(N)}$ denotes the Fejer approximant of order $N$ of $b_n$. Thanks to Corollary~\ref{met10.cor-BKbnd2} it follows that if $N<M$ then 

$$\|b^{(N)}-b^{(M)}\|\leq
   \limsup_{n\rightarrow \infty}
    \|b_n^{(N)}-b_n^{(M)}\|\leq
     \limsup_{n\rightarrow \infty} \left(
      \|b_n^{(N)}-b_n\|+ \|b_n-b_n^{(M)}\|
       \right)\leq
        \frac{4.4\pi}{\sqrt{N}}\,.
$$

\noindent It follows that the sequence $(b^{(N)})_{N\in\NM}$ is Cauchy and therefore it converges in $\as$ to a limit $b$.

\vspace{.1cm}

\noindent (iii) By construction $\partial b^{(N)}$ exists since $b^{(N)}\in\as_c$. For a fixed $N$ it satisfies 

$$\partial b^{(N)} =
   \sum_{l=-N+1}^{N-1} \imath l\,
    \left(1-\frac{|l|}{N}\right)\,
     b_l\,u^l =
      \lim_{n\rightarrow \infty} \partial b_n^{(N)}
$$

\noindent From then, it follows that $\|\partial b^{(N)}\|\leq \limsup_{n\rightarrow \infty} \|\partial b_n^{(N)}\|\leq \limsup_{n\rightarrow \infty} \|\partial b_n\| \leq 1$. In particular $b^{(N)}\in B(K)$. Since $B(K)$ is closed, thanks to Lemma~\ref{met10.lem-BKclosed}, it follows that $b\in B(K)$ as well.
\hfill $\Box$

\vspace{.3cm}

 \subsection{Quotient Spaces}
 \label{met10.ssect-quot}

\noindent Let $\Aa$ be a unital \Cs with unit denoted by $\id$. Let $\delta\Aa$ denote the Banach space $\Aa/\CM\id$ and let $[a]\in\delta\Aa$ denote the equivalence class of $a\in\Aa$ modulo $\CM\id$. The quotient norm is defined by 

$$\|[a]\|_1 = \inf_{\lambda\in\CM} \|a-\lambda\id\|\,.
$$

\noindent Let now $\omega$ be a state on $\Aa$ and let $\psi_\omega$ denote the map $\psi_\omega(a) = ([a],\omega(a))\in \delta\Aa\times \CM$. The space $\delta\Aa\times \CM$ is a complex vector space that will become a Banch space if endowed with the norm $(c,z)=\|c\|_1+|z|$. 

\begin{lemma}
\label{met10.lem-quotiso}
The map $\psi_\omega:\Aa\mapsto \Aa/\CM\id\times \CM$ is linear, bounded, invertible with bounded inverse.
\end{lemma}

\noindent  {\bf Proof: } The linearity and continuity of $\psi_\omega$ is obvious from the definition. If $a,b\in\Aa$ satisfy $\psi_\omega(a)=\psi_\omega(b)$, then $[a]=[b]$ implying that there is $\lambda\in\CM$ such that $b=a+\lambda\id$. Then the equation $\omega(a)=\omega(b)$ implies that $\lambda =0$ showing that $a=b$. Hence $\psi_\omega$ is one-to-one. To show that it is onto, let $(c,z)\in\delta\Aa\times \CM$. Then there is $a\in\Aa$ such that $c=[a]$. If $\lambda=z-\omega(a)\in\CM$ it follows that $[a+\lambda\id]=c$ and $\omega(a+\lambda\id)=z$. Hence $\psi_\omega(a+\lambda\id)=(c,z)$. Since the inverse is everywhere defined, by the closed graph theorem, it follows that the inverse is continuous as well. In particular, the inverse image of a compact set in $\delta\Aa\times \CM$ by $\psi_\omega$ is compact.
\hfill $\Box$

\vspace{.3cm}

 \subsection{Proof of Proposition ~\ref{met10.prop-lipcomp}}
 \label{met10.ssect-prooflipcomp}

\noindent Let $\Clip$ denote the Lipschitz ball of $\XP$. By definition, $b\in\Clip$ if and only if $\|[\hD,\hpi(b)]\|\leq 1$. 

\vspace{.1cm}

\noindent (i) This condition implies two properties: (a) for all $l\in\ZM$, $b_l=\EM(bu^{-l})$ belongs the the Lipschitz ball $\Blip$ of $\Aa$, (b) $\|\partial b\|\leq 1$. For indeed, if $f,g\in\Hh\otimes \ell^2(\ZM)$ and if $e,e'\in\CM^2$, the inner product $\langle g\otimes e|[\hD,\hpi(b)]\,f\otimes e'\rangle$ can be written as $\langle e|\left(\langle g|[\hD,\hpi(b)]\,f\rangle\right)| e'\rangle$, where now $\langle g|[\hD,\hpi(b)]\,f\rangle$ is a $2\times 2$ matrix equal to $\mbox{\rm (I)} + \mbox{\rm (II)}$ with

\begin{eqnarray*}
\mbox{\rm (I)} &=& \sum_{n,l\in\ZM}
    \langle g_n |[D,\pi\circ\alpha^{-n}(b_l)]\;f_{n-l}\rangle\;
     \sigma_1\\
\mbox{\rm (II)} &=& \sum_{n,l\in\ZM}
      \langle g_n |\,l\,\pi\circ\alpha^{-n}(b_l)]\;f_{n-l}\rangle\;
       \sigma_2=
 -\imath\,\langle g|\hpi\left(\partial b\right)\,f\rangle \;
       \sigma_2\
\end{eqnarray*}

\noindent Let $\sigma_0$ denote the unit $2\times 2$ matrix. Then, the Pauli matrices satisfy $\TR(\sigma_i\sigma_j)=2\delta_{i,j}$ for $i,j\in \{0,1,2,3\}$. It follows that

$$\frac{\imath}{2}\TR\left(
   \langle g|[\hD,\hpi(b)]\,f\rangle
                     \right) =
    \langle g|\hpi\left(\partial b\right)\,f\rangle\,,
$$

\noindent showing that $b\in\Clip$ implies $\|\partial b\|\leq 1$. In particular, $\imath lb_l=\EM(\partial b\;u^{-l})$ implies that $\|b_l\|\leq 1/|l|$, whenever $l\neq 0$. On the other hand the same type of formula with $\sigma_2$ replaced by $\sigma_1$ implies that $\|[D,\pi\circ\alpha^{-n}(b_l)]\|\leq 1$ for all $n,l\in\ZM$. In particular $b_l\in\Blip$ for all $l\in\ZM$.

\vspace{.1cm}

\noindent (ii) Let $K_l$ be the closure of the set of $a\in\Blip\subset \Aa$ such that $\|a\|\leq 1/|l|$ for $l \neq 0$. If $\omega$ is any state on $\Aa$ it follows that $a\in K_l\;\Rightarrow\; |\omega(a)|\leq 1/|l|$. Let $\psi_\omega$ denote the map $\psi_\omega(a) = ([a],\omega(a))\in \Aa/\CM\id\times \CM$ defined in Section~\ref{met10.ssect-quot}. Then the image $\psi_\omega(K_l)$ is included in the closure of $[\Blip]\times\{z\in\CM: |z|\leq 1/|l|\}$, which is compact. Consequently, thanks to Lemma~\ref{met10.lem-quotiso}, $K_l$ is compact. It follows that, the subset $C_0\subset \Clip$ of elements with $\EM(b)=0$ is included in $B(K)$. Thanks to the Proposition~\ref{met10.prop-comp}, this set has a compact closure. Since any element of $\Clip$ has the form $b= b_0 + b-\EM(b)$, and since $b_0\in\Blip$, it follows that the image of $\Clip$ inside $\XP/\CM\id$ has also a compact closure, because $b-\EM(b)\in B(K)$ and $[\Blip]$ has a compact closure.
\hfill $\Box$

\vspace{.5cm}


\begin{thebibliography}{9999}


\bibitem{ACh04} C.~Antonescu, E.~Christensen, ``Metrics on group $C^*$-algebras and a non-commutative Arzel\`a-Ascoli theorem'', {\em J. Funct. Anal.} {\bf 214} (2004), 247-259.

\bibitem{AA67} V.~I.~Arnold, A.~Avez, {\em Probl\`emes ergodiques de la m\'ecanique classique}, (French) Gauthier-Villars, Paris 1967. English translation: V.~I.~Arnold, A.~Avez, {\em Ergodic problems of classical mechanics}, W. A. Benjamin, Inc., New York-Amsterdam, 1968.

\bibitem{Ar98} W.~Arveson, {\em An invitation to C$^\ast$-algebras}, 2nd printing, Graduate Text in Mathematics, Springer, 1998.

\bibitem{At67} M.F.~Atiyah, {\em $K$-theory}, Lecture notes by D. W. 
Anderson, W. A. Benjamin, Inc., New York-Amsterdam, 1967.

\bibitem{ADS} M.F.~Atiyah, H.~Donnelly, I.M.~Singer, ``Eta invariants, signature defects of cusps, and values of $L$-functions'', {\em Ann. of Math.} {\bf 118} (1983), 131-177.

\bibitem{Be32} A.~S.~Besicovitch, {\em Almost periodic functions}, Cambridge Univ. Press, 1932.

\bibitem{Bi99} P.~Billingsley, {\em Convergence of Probability Measures}, New York, NY, John Wiley \& Sons, Inc., 1999.

\bibitem{BvN35} S.~Bochner, J.~von~Neumann, ``Almost Periodic Function in a Group II'', {\em Trans. Amer. Math. Soc.} {\bf 37} (1935), 21-50.

\bibitem{Bo47} H.~Bohr, {\em Almost-periodic functions}, Chelsea, reprint, 1947.

\bibitem{Bu99} S.~V.~Buyalo, ``Spectral geometries on a compact metric space'', {\em POMI Preprint}, vol. 19, (1999).

\bibitem{Bu00} S.~V.~Buyalo, ``Measurability of self-similar spectral geometries'', {\em Algebra i Analiz} {\bf 12} (2000), 1-39; translation in  {\em St. Petersburg Math. J.}  {\bf 12}  (2001),  no. 3, 353-377.

\bibitem{CaMaRe} A.~Carey, M.~Marcolli, A.~Rennie, 
``Modular index invariants of Mumford curves", {\tt arXiv:0905.3157} (2009).

\bibitem{ChI06} E.~Christensen, C.~Ivan, ``Spectral triples for AF \Css and metrics on the Cantor set'', {\em J. Operator Theory} {\bf 56} (2006), 17-46.

\bibitem{ChI07} E.~Christensen and C.~Ivan, ``Sums of Two Dimensional Spectral Triples'', {\em Math. Scand.} {\bf 100} (2007), 35-60.

\bibitem{CIL} E.~Christensen, C.~Ivan, M.~Lapidus, ``Dirac Operators and Spectral Triples for some Fractal Sets Built on Curves'', {\em Adv. Math.} {\bf 217} (2008), 42-78.

\bibitem{CIS} E.~Christensen, C.~Ivan, E.~Schrohe, ``On the choice of a spectral triple'', {\tt arXiv:1002.3081v1}, (2010).

\bibitem{Co73}  A.~Connes, ``Une classification des facteurs de type ${\rm III}$'', (French) {\em Ann. Sci. \'Ecole Norm. Sup.} (4) {\bf 6} (1973), 133-252.

\bibitem{Co88} A.~Connes, ``Trace de Dixmier, modules de Fredholm et g\'eom\'etrie riemannienne'', in {\em Conformal field theories and related topics (Annecy-le-Vieux, 1988)}, {\em Nuclear Phys. B Proc. Suppl.} {\bf 5B} (1988), 65-70. 

\bibitem{Co94} A.~Connes, {\em Noncommutative Geometry}, Academic Press, San Diego, 1994.

\bibitem{Co95} A.~Connes, ``Gravity coupled with matter and the foundation of non-commutative geometry'', {\em  Comm. Math. Phys.} {\bf 182} (1996), 155-176.

\bibitem{CM95} A.~Connes, H.~Moscovici, ``The local index formula in noncommutative geometry'', {\em Geom. Funct. Anal.} {\bf 5} (1995), 174-243.

\bibitem{CoLa} A.~Connes, G.~Landi, ``Noncommutative manifolds, the 
instanton algebra and isospectral deformations'', {\em Comm. Math. Phys.} {\bf 221} (2001), 141-159. 

\bibitem{Co08} A.~Connes, ``On the Spectral Characterization of Manifolds'', {\tt arXiv:0810.2088}, (2008).

\bibitem{ConsMa3} C.~Consani, M.~Marcolli, ``Archimedean cohomology revisited", {\em Noncommutative Geometry and Number Theory}, 109--140, Vieweg Verlag, 2006.

\bibitem{ConsMa2} C.~Consani, M.~Marcolli, ``Noncommutative geometry, dynamics, and $\infty$-adic Arakelov geometry", {\em Selecta Mathematica} {\bf 10} (2004), no. 2, 167--251.

\bibitem{ConsMa} C.~Consani, M.~Marcolli, ``Spectral triples from Mumford curves'', {\em International Mathematics Research Notices (IMRN)} {\bf 36} (2003), 1945--1972.

\bibitem{CMRV} G.~Cornelissen, M.~Marcolli, K.~Reihani, A.~Vdovina,
``Noncommutative geometry on trees and buildings'', in {\em Traces in Geometry, Number Theory, and Quantum Fields}, Aspects Math., E38, Friedr. Vieweg, Wiesbaden, 2008, pp. 73-98.

\bibitem{Do70} R.~L.~Dobrushin, ``Prescribing a system of random variables by conditional distributions'', {\em Theor. Probab. Appl.} {\bf 15} (1970), 458-486.

\bibitem{Eilen} S.~Eilenberg, {\em Automata, Languages, and Machines}, 
vol. A. Academic Press, London, 1974.

\bibitem{FS10} K.~J.~Falconer, T.~Samuel, ``Dixmier traces and coarse multifractal analysis'', {\em Ergod. Th. \& Dynam. Sys.}, First view 2010, pp. 1-13.

\bibitem{GBVF} J.~Garcia-Bond\'{\i}a, J.~V\'{a}rilly, H.~Figueroa, {\em Elements of noncommutative geometry}, Birkh\"{a}user Advanced Texts: Basler Lehrb\"{u}cher, Birkh\"{a}user Boston, Inc., Boston, MA, 2001.

\bibitem{GPS95} T.~Giordano, I.~F.~Putnam, C.~F.~Skau, ``Topological orbit equivalence and $C^\ast$-crossed products'', {\em J. Reine Angew. Math.} {\bf 469} (1995), 51-111.

\bibitem{GPS10} T.~Giordano, H.~Matui, I.~F.~Putnam, C.~F.~Skau, ``Orbit equivalence for Cantor minimal $\ZM^d$-systems'', {\em Invent. Math.} {\bf 179} (2010), 119-158.

\bibitem{GI01} D.~Guido, T.~Isola, ``Fractals in noncomutative geometry'', in {\em Mathematical Physics in Mathematics and Physics (Sienna 2000)}, vol. 30, {\em Fields Inst. Commun.}, Providence, RI, AMS, 2001, pp. 171-186.

\bibitem{GI03} D.~Guido, T.~Isola, ``Dimensions and Singular Traces for Spectral Triples, with Applications to Fractals'', {\em J. Funct. Anal.} {\bf 203} (2003), 362-400.

\bibitem{GI05} D.~Guido, T.~Isola, ``Dimension and Spectral Triples for Fractals in $\RM^N$''.  In {\em Advances in Operator Algebras and Mathematical Physics}, vol. 5, {\em Theta Ser. Adv. Math.}, Theta, Bucharest, 2005, pp. 89-108.

\bibitem{Gr99} M.~Gromov, {\em Metric structures for Riemannian and non-Riemannian spaces}, Birkh\"auser, 1999.

\bibitem{Hirz} F.~Hirzebruch,``Hilbert modular surfaces'',  {\em Enseignement Math.} (2) {\bf 19} (1973), 183-281.

\bibitem{JS09} A.~Julien, J.~Savinien, ``Transverse Laplacians for Substitution Tilings'', {\tt arXiv:0908.1095} (2009), to be published in {\em Commun. Math. Phys.}.

\bibitem{Ka40} L.~V.~Kantorovich, ``On one effective method of solving certain classes of extremal problems'', {\em Dokl. Akad. Nauk USSR} {\bf 28} (1940), 212-215.

\bibitem{KR57} L.~V.~Kantorovi\v{c}, G.~\v{S}.~Rubin\v{s}te\u{\i}n, ``On a functional space and certain extremum problems'', (Russian) {\em Dokl. Akad. Nauk SSSR (N.S.)} {\bf 115} (1957), 1058-1061; ``On a space of completely additive functions'', (Russian) {\em Vestnik Leningrad. Univ.} {\bf 13} (1958), 52-59.

\bibitem{KS07} L.~B.~Koralov, Ya.~G.~Sinai, {\em Theory of probability and random processes}, Springer, 2nd Ed., 2007.

\bibitem{La94} M.~L.~Lapidus, ``Analysis on fractals, Laplacians on Self-similar sets, noncommutative geometry and spectral dimensions'', {\em Topol. Methods Nonlinear Anal.} {\bf 40} (1994), 137-195.

\bibitem{La97} M.~L.~Lapidus, ``Towards a noncommutative fractal geometry? Laplacians and volume measures on fractals'', in {\em Harmonic analysis and nonlinear differential equations (Riverside, CA, 1995)}, vol. 208 of {\em Contemp. Math.}, Providence, RI, AMS, 1997, pp. 211-252.

\bibitem{La07} F.~Latr\'emoli\`ere, ``Bounded-Lipschitz distances on the state space of a \CS'', {\em Taiwanese J. of Math.} {\bf 11} (2007), 447-469.

\bibitem{Man} Yu.I.~Manin, ``Real multiplication and noncommutative geometry (ein Alterstraum)". {\em The legacy of Niels Henrik Abel}, Springer, Berlin, 2004,  pp. 685-727.

\bibitem{Man3d} Yu.I.~Manin, ``Three-dimensional hyperbolic geometry as $\infty$-adic Arakelov geometry", {\em Invent. Math.} {\bf 104} (1991), 223-244.

\bibitem{ManMar} Yu.I.~Manin, M.~Marcolli, ``Algebro-geometric codes and phase transitions'', to appear in {\em Mathematics in Computer Science}.

\bibitem{Mar-solv} M.~Marcolli, ``Solvmanifolds and noncommutative tori with real multiplication'', {\em Commun. Number Theory and Phys.} {\bf 2} (2008), no. 2, 423-479.

\bibitem{Mar-ICM} M.~Marcolli, ``Noncommutative Geometry and Arithmetic'', {\tt arXiv:1003.3662}, to appear in the {\em Proceedings of ICM-2010}.

\bibitem{vN34}  J.~von~Neumann, ``Almost Periodic Functions in a Group I'', {\em Trans. Amer. Math. Soc.} {\bf 36}
(1934), 445-492.

\bibitem{Oz-Ri} N.~Ozawa, M.~Rieffel, ``Hyperbolic group $C^*$-algebras and free-product $C^*$-algebras as compact quantum metric spaces", {\em Canad. J. Math.} {\bf 57} (2005), no. 5, 1056--1079. 

\bibitem{Pa98} B.~Pavlovi\'c, ``Defining metric spaces viaoperators from unital $C^\ast$-algebras'', {\em Pacific J. Math.} {\bf 186} (1998), 285-313.

\bibitem{Pa10} I.~Palmer, ``Noncommutative Geometry and Compact Metric Spaces'', PhD Thesis, May 2010, Georgia Institute of Technology.

\bibitem{PB09} J.~Pearson, J.~Bellissard, ``Noncommutative Riemannian Geometry and Diffusion on Ultrametric Cantor Sets'', {\em J.  Noncommut. Geo.} {\bf 3} (2009), 447-480.

\bibitem{Ped79} G.~K.~Pedersen, {\em $C^{\ast}$-Algebras and their Automorphism Groups}, Academic Press, London, 1979.

\bibitem{Pr56} Y.~V.~Prokhorov, ``Convergence of random processes and limit theorems in probability theory'', (in English translation) {\em Theory of Prob. and Appl.}  {\bf I 2} (1956), 157-214. 

\bibitem{Re01} A.~Rennie, ``Commutative geometries are spin manifolds'',  {\em Rev. Math. Phys.} {\bf 13} (2001), 409-464. 

\bibitem{RV06} A.~Rennie, J.~C.~Varilly, ``Reconstruction of manifolds in noncommutative geometry'', {\tt arXiv:math/0610418}, (2006).

\bibitem{Ri76} M.~Rieffel, ``Strong Morita equivalence of certain transformation group $C^\ast$-algebras'', {\em Math. Ann.} {\bf 222} (1976), 7-22.

\bibitem{Ri98} M.~Rieffel, ``Metrics on states from actions of compact groups'', {\em Doc. Math.} {\bf 3} (1998), 215-229.

\bibitem{Ri99} M.~Rieffel, ``Metrics on states spaces'', {\em Doc. Math.} {\bf 4} (1999), 559-600.

\bibitem{Ri04} M.~Rieffel, ``Compact quantum metrics spaces'', in {\em Operator Algebras, Quantization, and Noncommutatve Geometry: a Centennial Celebration Honoring John von Neumann and Marshall H. Stone}, ({\sc Doran, R.~S.} and {\sc Kadison, R.~V.}, eds.), vol. 365 of {\em Contemporary Mathematics}, AMS, 2004, pp. 315-330.

\bibitem{TaH74} H.~Takai, ``Dualit\'e dans les produits crois\'es de $C^{\ast} $-alg\`ebres'', (French), {\em C. R. Acad. Sci. Paris S\'er. A} {\bf 278} (1974), 1041-1043.

\bibitem{TaH75} H.~Takai, ``On a duality for crossed products of $C^{\ast} $-algebras, {\em J. Functional Analysis} {\bf 19} (1975), 25-39.

\bibitem{TaM73} M.~Takesaki, ``Duality for crossed products and the structure of von Neumann algebras of type III'', {\em Acta Math.} {\bf 131} (1973), 249-310.

\bibitem{Wa69} L.~N.~Wasserstein, ``Markov processes over denumerable products of spaces describing large systems of automata'', {\em Probl. Inform. Transmission} {\bf 5} (1969), 47-52.

\end{thebibliography}
\end{document}